\documentclass[dvips]{compositio}
\usepackage{amsthm,amsmath,amssymb,amscd,mathrsfs}
\input diagrams
\diagramstyle[PostScript=dvips]
\newcommand{\includefig}[1]{\raisebox{-3ex}{\resizebox{!}{7ex}{\includegraphics{pix/#1}}}}
\newcommand{\tree}{\includefig{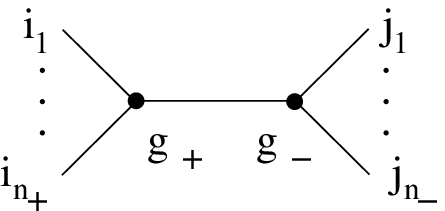}}
\newcommand{\treecut}{\includefig{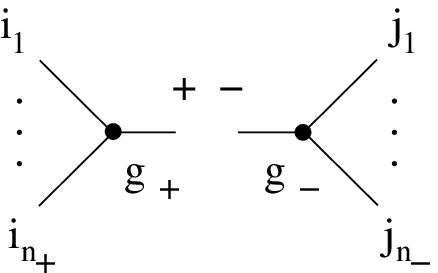}}

\newcommand{\treemb}{\includefig{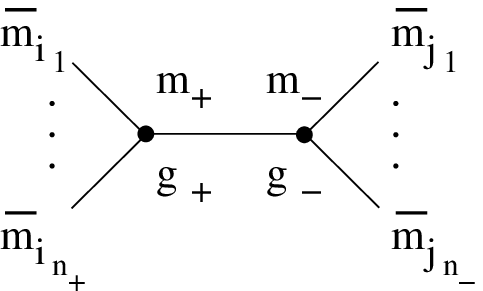}}
\newcommand{\treembcut}{\includefig{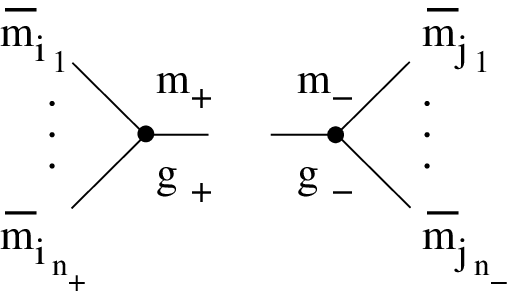}}
\newcommand{\treembb}{\includefig{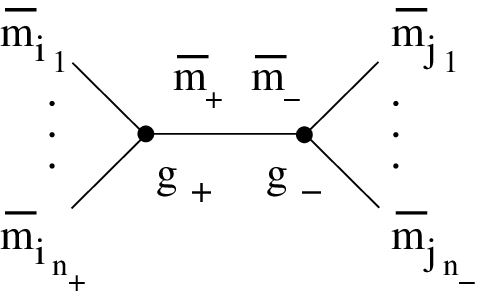}}
\newcommand{\treembbcut}{\includefig{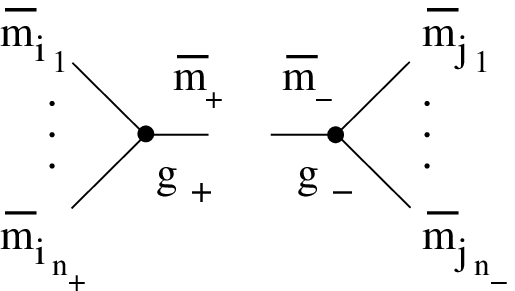}}

\newcommand{\loopp}{\includefig{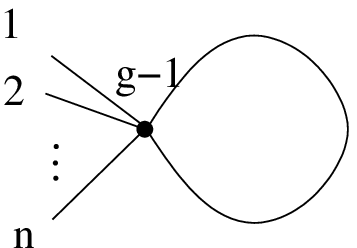}}
\newcommand{\loopcut}{\includefig{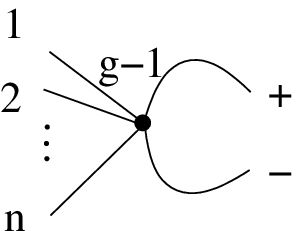}}

\newcommand{\loopmb}{\includefig{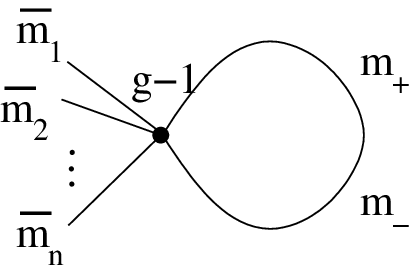}}
\newcommand{\loopmbcut}{\includefig{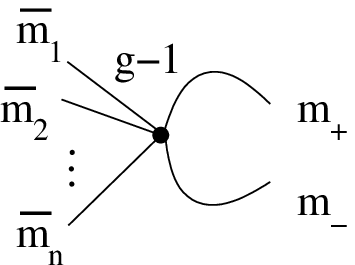}}
\newcommand{\loopmbb}{\includefig{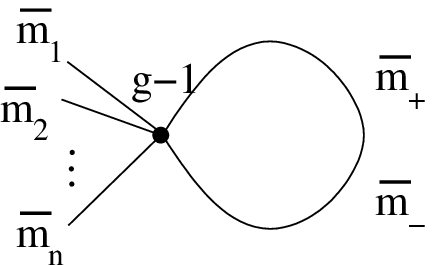}}
\newcommand{\loopmbbcut}{\includefig{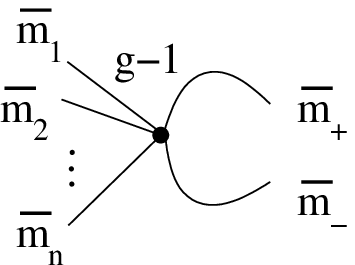}}

\newcommand{\ce}{{\mathscr E}}
\newcommand{\co}{{\mathscr O}}
\newcommand{\cc}{{\mathscr C}}
\newcommand{\ch}{{\mathscr H}}
\newcommand{\ct}{{\mathscr T}}
\newcommand{\bg}{{\boldsymbol{\gamma}}}
\newcommand{\stack}[1]{\mathscr {#1}}
\newcommand{\mgnbar}{\overline{\stack{M}}_{g,n}}

\newcommand{\nc}{\mathbb{C}}
\newcommand{\nq}{\mathbb{Q}}
\newcommand{\nr}{\mathbb{R}}
\newcommand{\xg}{[X/G]}

\newcommand{\Aut}{\operatorname{Aut}}
\newcommand{\be}{\mathbf{e}}

\newcommand{\bm}{\mathbf{m}}
\newcommand{\bmb}{{\overline{\mathbf{m}}}}

\newcommand{\cft}{CohFT}
\newcommand{\cfts}{CohFTs}
\newcommand{\chb}{\overline{\ch}}

\newcommand{\ctt}{\widetilde{c}}

\newcommand{\db}[1]{[\![{#1}]\!]} % renormalized cohomology classes
\newcommand{\Deltac}{\Delta_\mathrm{cut}}
\newcommand{\Deltat}{\widetilde{\Delta}}
\newcommand{\etab}{\overline{\eta}}
\newcommand{\etah}{{\widehat{\eta}}}

\newcommand{\gammat}{\widetilde{\gamma}}
\newcommand{\Gammah}{{\widehat{\Gamma}}}
\newcommand{\Gammat}{{\widetilde{\Gamma}}}
\newcommand{\Gammac}{{\Gamma_\mathrm{cut}}}
\newcommand{\Gammatc}{\Gammat_{\mathrm{cut}}}
\newcommand{\Gammahc}{\Gammah_{\mathrm{cut}}}
\newcommand{\Gb}{\overline{G}}
\newcommand{\gb}{\overline{\gamma}}
\newcommand{\Gcft}{$G$-\cft}
\newcommand{\Gcfts}{$G$-\cft s}
\newcommand{\Gpcft}{$G'$-\cft}

\newcommand{\Gppcft}{$G''$-\cft}

\newcommand{\ih}{\widehat{i}}
\newcommand{\im}{\operatorname{im}}  %%image of a map
\newcommand{\Ic}{\hat{I}}
\newcommand{\It}{\widetilde{I}}
\newcommand{\itt}{\widetilde{i}}
\newcommand{\irightarrow}{\rTo^{\sim}} %%isom arrow

\newcommand{\Lambdab}{\overline{\Lambda}}
\newcommand{\Lambdac}{\check{\Lambda}}
\newcommand{\Lambdah}{\widehat{\Lambda}}

\newcommand{\M}{\overline{\mathscr{M}}}
\newcommand{\MM}{\M^G} % small moduli spaces
\newcommand{\MMp}{\M^{G'}} % small moduli spaces
\newcommand{\MMpp}{\M^{G''}} % small moduli spaces
\newcommand{\MMq}{\mathscr{Q}} % Quotient stack
\newcommand{\mb}{{\overline{m}}}
\newcommand{\mub}{\overline{\mu}}
\newcommand{\muh}{\widehat{\mu}}
\newcommand{\mut}{\widetilde{\mu}}
\newcommand{\nz}{\mathbb{Z}}

\newcommand{\pit}{\widetilde{\pi}}

\newcommand{\pr}{\mathrm{pr}}
\newcommand{\prc}{\mathrm{pr}_\mathrm{cut}}
\newcommand{\prh}{\widehat{\mathrm{pr}}}
\newcommand{\prt}{\widetilde{\mathrm{pr}}}

\newcommand{\pt}{\widetilde{p}}
\newcommand{\rh}{\widehat{r}}
\newcommand{\rt}{\widetilde{r}}
\newcommand{\st}{\mathrm{st}}
\newcommand{\stc}{\check{\mathrm{st}}}
\newcommand{\stcp}{{\check{\mathrm{st}}'}}
\newcommand{\sth}{\widehat{\mathrm{st}}}
\newcommand{\stt}{\widetilde{\mathrm{st}}}
\newcommand{\vac}{\mathbf{1}}
\newcommand{\varrhot}{\widetilde{\varrho}}

\newcommand{\id}{\boldsymbol{1}}
\newcommand{\spec}{\operatorname{Spec}}

\newcommand{\BG}{{\mathscr{B}}G}
\newcommand{\CP}{\mathbb{CP}}
\newcommand{\bgamma}{\boldsymbol{\gamma}}
\newcommand{\ts}{\tau} %                          %forgetting tails map
\newcommand{\tst}{\widetilde{\tau}} %  %forgetting tails map on \MM
\newcommand{\tsh}{\widehat{\tau}} %    %forgetting tails map on \M(\BG)
\newcommand{\X}{\widehat{X}} % target of the eval morphism for X/G
\newcommand{\stab}{\operatorname{stab}}
\newcommand{\bd}{\mathbf{d}} %n-tuple of paths from p_0 to p_i
\newcommand{\str}{{{\widetilde{\mathrm{st}}_X}}} %stab map from \beta=0 G-stable maps to G-covers.
\newcommand{\stq}{{{\mathrm{st}}_{(X,\beta)}}} %stab map from all
\newcommand{\E}{\ce} %universal G-bundle

\newcommand{\vrta}{\varrhot_a}
\newcommand{\itimes}{\odot} %interior tensor product of cohfts
\newcommand{\etimes}{\otimes} %external tensor product of cohfts
\newcommand{\btimes}{\itimes} % tensor product of G-graded G-modules

\newtheorem{quest}{Question}

\newtheorem{thm}{Theorem}[section]
\newtheorem{lm}[thm]{Lemma}
\newtheorem{prop}[thm]{Proposition}
\newtheorem{crl}[thm]{Corollary}

\theoremstyle{definition}

\newtheorem{rem}[thm]{Remark}

\newtheorem{df}[thm]{Definition}
\newtheorem{ex}[thm]{Example}
\newtheorem{df-pr}[thm]{Definition-Proposition}

\theoremstyle{remark}
\newtheorem{nota}[thm]{Notation} \renewcommand{\thenota}{\kern-1ex}
\newenvironment{ack}{\noindent \emph{Acknowledgments}.}

\begin{document}

\title[$G$-equivariant Cohomological
Field Theories] {Pointed Admissible $G$-Covers and
$G$-equivariant Cohomological Field Theories}


\author[T. J. Jarvis]{Tyler J. Jarvis}\email{jarvis@math.byu.edu}
\address{Department of Mathematics, Brigham Young University, Provo, UT
84602, USA}

\author[R. Kaufmann]{Ralph Kaufmann}\email{kaufmann@math.uconn.edu}
\address{Department of Mathematics, University of Connecticut, 196 
Auditorium Road, Storrs, CT 06269-3009, USA}

\author[T. Kimura]{Takashi Kimura} \email{kimura@math.bu.edu}
\address{Department of Mathematics and Statistics; 111 Cummington
  Street, Boston University; Boston, MA 02215, USA and School of
  Mathematics; Institute for Advanced Study; 1 Einstein Dr.;
  Princeton, NJ 08540, USA } 

\thanks{ Research of the first
author was partially supported by NSF grant DMS-0105788.  Research of
the second author was partially supported by NSF grant DMS-0070681.
Research of the third author was partially supported by NSF grants
DMS-0204824 and DMS-9729992\\
Published in \emph{Compositio Mathematica}, \textbf{141}
  (2005) 926--978. doi:10.1112/S0010437X05001284.  \copyright
Foundation Compositio Mathematica, 2005.}

\date{\today}

\begin{abstract}
For any finite group $G$ we define the moduli space of pointed
admissible $G$-covers and the concept of a $G$-equivariant
cohomological field theory (\Gcft), which, when $G$ is the trivial
group, reduce to the moduli space of stable curves and a cohomological
field theory (\cft), respectively. We prove that taking the
``quotient'' by $G$ reduces a \Gcft\ to a \cft. We also prove that
a \Gcft\ contains a $G$-Frobenius algebra, a $G$-equivariant
generalization of a Frobenius algebra, and that the ``quotient'' by
$G$ agrees with the obvious Frobenius algebra structure on the space
of $G$-invariants, after rescaling the metric. We then introduce the
moduli space of $G$-stable maps into a smooth, projective variety $X$
with $G$ action.  Gromov-Witten-like invariants of these spaces
provide the primary source of examples of \Gcfts.  Finally, we explain
how these constructions generalize (and unify) the Chen-Ruan orbifold
Gromov-Witten invariants of $\xg $ as well as the ring
$H^{\bullet}(X,G)$ of Fantechi and G\"ottsche.
\end{abstract}

\maketitle

%\tableofcontents

\section{Introduction}

\ 

The purpose of this paper is to introduce a generalization of
Kontsevich and Manin's notion of a cohomological field theory (or
\cft) \cite{KoMa}, in the presence of a finite group $G$, which we
call a $G$-equivariant cohomological field theory (or \Gcft).
Examples of (usual) \cfts\ include the Gromov-Witten invariants of
a smooth, projective variety (cf. \cite{Ma}) and the $r$-spin
\cft\ \cite{JKV,PoVa,Po2}.  A \Gcft\ provides a framework for
studying the physical procedure of orbifolding \cite{Ka,Ka2,Mo}, as
well as a structure for understanding both Chen-Ruan orbifold
Gromov-Witten invariants of global quotients by a finite group
\cite{CR1,CR2,AGV} and the non-commutative ring structure of
Fantechi and G\"ottsche \cite{FG}.  We now describe in some detail
the motivation for studying \Gcfts.

The first motivation comes from topological field theory. Recall
that a Frobenius algebra $\ch$ is a finite-dimensional,
commutative, associative, unital algebra with an invariant metric.
It can be regarded as a two-dimensional topological field theory,
in the sense of Atiyah-Segal, associated to a cobordism category
of two (real) dimensional, compact, oriented surfaces with
boundary. A \cft\ is a generalization of the above, but where the
role of the cobordism category is replaced by $\{H_r(\M_{g,n})\}$
for all $r$, where $\M_{g,n}$ is the moduli space of stable curves
of genus $g$ with $n$ marked points.  By specializing to $r=0$,
one finds that the state space of the theory $\ch$ recovers the
structure of a Frobenius algebra.

For every finite group $G$ Turaev \cite{Tu} introduced a
$G$-equivariant topological field theory (which he called a
homotopy field theory) whose state space $\ch$ is a
(non-projective) $G$-Frobenius algebra associated to a cobordism
category of principal $G$-bundles over two (real) dimensional,
compact, oriented surfaces with boundary. A (non-projective)
$G$-Frobenius algebra (borrowing terminology from \cite{Ka,Ka2}) is a
finite-dimensional, $G$-graded $G$-module with a $G$-equivariant
associative multiplication, metric, and unit, and whose
multiplication is braided commutative, satisfying an additional
genus-one compatibility condition (called the trace axiom). By
\emph{braided commutative} we mean that the multiplication
commutes with the action of the generator of the braid group which
acts on tensor products of $G$-graded $G$-modules. If $G$ is the
trivial group, then a $G$-Frobenius algebra is a Frobenius algebra.
Furthermore, the space of $G$-invariants $\chb$ of a $G$-Frobenius
algebra inherits the structure of a Frobenius algebra graded by
$\Gb$, the set of conjugacy classes of $G$. Kaufmann \cite{Ka,Ka2}
considered a generalization of the above construction which
allowed for projective factors.

This procedure of restricting to the space of invariants can be
interpreted as a kind of orbifolding procedure from physics
\cite{Ka,Ka2,Mo} where the subspace of $\chb$ graded by $1$ in $\Gb$
is called the \emph{untwisted sector}, and the subspaces graded by
nontrivial elements in $\Gb$ are called \emph{twisted sectors}.

\begin{quest}\label{qu:one} Is there a generalization of
a \cft, called a \Gcft, where $\M_{g,n}$ is replaced by another
moduli space $\MM_{g,n}$, such that for all $r$, the collection
$\{ H_r(\MM_{g,n}) \}$ endows the state space $\ch$ of the theory
with an algebraic structure whose specialization to $r=0$ induces
the structure of a $G$-Frobenius algebra on $\ch$? A \Gcft\ should
also have the property that when $G$ is the trivial group, a
\Gcft\ reduces to a \cft. Furthermore, by performing the correct
``quotient'' by $G$, the space of $G$-invariants $\chb$ should
inherit the structure of a \cft\ graded by $\Gb$.
\end{quest}

\begin{figure}[hbp]
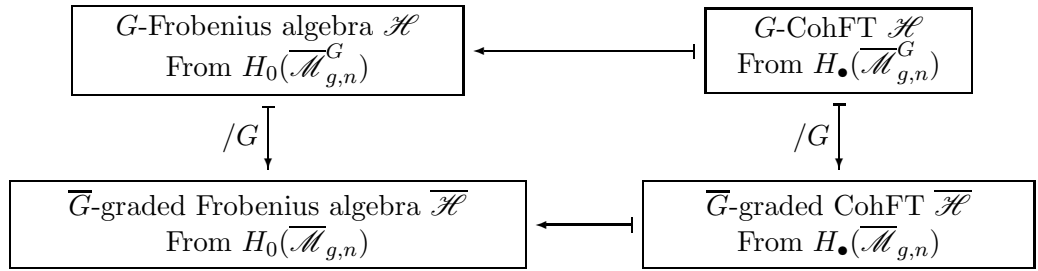

\[
\begin{diagram}
\fbox{\parbox{.3\linewidth}{\centerline{$G$-Frobenius algebra $\ch$}
\centerline{From $H_0(\MM_{g,n})$}}} & \lMapsto &
&
\fbox{\parbox{.2\linewidth}{\centerline{\Gcft\  $\ch$} \centerline{From
$H_\bullet(\MM_{g,n})$}}} \\
\dMapsto^{/G} &                                  &  & \dMapsto^{/G }
\\
\fbox{\parbox{.4\linewidth}{\centerline{$\Gb$-graded Frobenius algebra
$\chb$} \centerline{From $H_0(\M_{g,n})$}}} &
\lMapsto & &
\fbox{\parbox{.3\linewidth}{\centerline{$\Gb$-graded \cft\  
$\chb$} \centerline{From $H_\bullet(\M_{g,n})$}}}
\\
\end{diagram}
\]
\caption{Schematic of Question (\ref{qu:one}) where each box contains
  an algebraic structure and the responsible homology group, the
  horizontal arrows are restrictions, and the vertical arrows are
  ``quotients'' by $G$. When $G$ is the trivial group, the two rows
  coincide.}
\end{figure}

The second motivation for studying \Gcfts\ comes from orbifold
Gromov-Witten invariants and is about how to construct certain
examples of \Gcfts\ associated to a smooth, projective variety $X$
with an action of a finite group $G$. 

Consider the $G$-graded $G$-module $\ch(X) := \bigoplus_{m\in G}
H^\bullet(X^m)$, where $X^m$ denotes the fixed-point set in $X$ of
$m$, and let $\chb(X)$ denote its space of $G$-invariants. Chen and
Ruan \cite{CR1,CR2} introduced the notion of Gromov-Witten invariants
for orbifolds, which, when applied to the global quotient $\xg $, has
a state space isomorphic to $\chb(X)$.  An algebro-geometric version
of this theory was introduced by \cite{AGV}. The key geometric object
in these constructions was $\M_{g,n}(\xg )$, the moduli space of
orbifold stable maps into the quotient $\xg $. The state space
$\chb(X)$ of this theory is graded by $\Gb$, and the Gromov-Witten
invariants are expected to yield a \cft\ associated to each $\xg $.
An important special case arises by considering only those
contributions from $\M_{g,n}(\xg ,0)$, the moduli of orbifold stable
maps which have degree zero.  This endows $\chb(X)$ with the structure
of a Frobenius algebra graded by $\Gb$, called variously \emph{stringy
  orbifold cohomology}, \emph{Chen-Ruan cohomology}, or just
\emph{orbifold cohomology} of $\xg $.

When $G$ is a trivial group, $\M_{g,n}(\xg )$ reduces to the usual
moduli space $\M_{g,n}(X)$ of stable maps into $X$, and the
Gromov-Witten invariants of $X$ make $H^\bullet(X)$ into a
\cft. Restricting to  contributions from $\M_{g,n}(\xg ,0)$ alone, one
obtains the usual cohomology ring of $X$, which is a Frobenius
algebra. ``Forgetting'' the stable map yields a morphism
$\M_{g,n}(X)\rTo \M_{g,n}$ for all stable pairs $(g,n)$, which is an
isomorphism when $X$ is a point.

Fantechi and G\"ottsche \cite{FG} were able to obtain the
structure of the Chen-Ruan orbifold cohomology on $\chb(X)$ by
first introducing a certain ring structure with metric on $\ch(X)$ and
then taking $G$-invariants. In fact, their ring satisfies all of
the axioms of a $G$-Frobenius algebra except, possibly, the trace
axiom. However, their construction is not obviously part of a
larger structure and does not explicitly involve the moduli space
of orbifold stable maps.

\begin{quest}\label{qu:two}
For any smooth, projective variety $X$ with a $G$-action, does there
exist a moduli space $\MM_{g,n}(X)$ of a $G$-equivariant version of
stable maps such that ``forgetting'' the map yields a morphism
$\MM_{g,n}(X) \rTo \MM_{g,n}$ for stable pairs $(g,n)$?  This map
should be an isomorphism when $X$ is a point.

There should also exist $G$-equivariant Gromov-Witten invariants
associated to $\MM_{g,n}(X)$ which yield a \Gcft\ with state space
$\ch(X)$, generalizing the usual construction when $G$ is the trivial
group. Furthermore, by taking the appropriate ``quotient'' by $G$, one
should recover the orbifold Gromov-Witten invariants of $\xg $ as in
\cite{CR1,CR2,AGV} with associated state space $\chb(X)$.

Finally, by considering only those contributions from the moduli
$\MM_{g,n}(X,0)$ of stable maps of degree zero, one should be able to
recover the $G$-Frobenius algebra structure in \cite{FG} and
  prove that the trace axiom must hold.
\end{quest}
\begin{figure}[hbt]
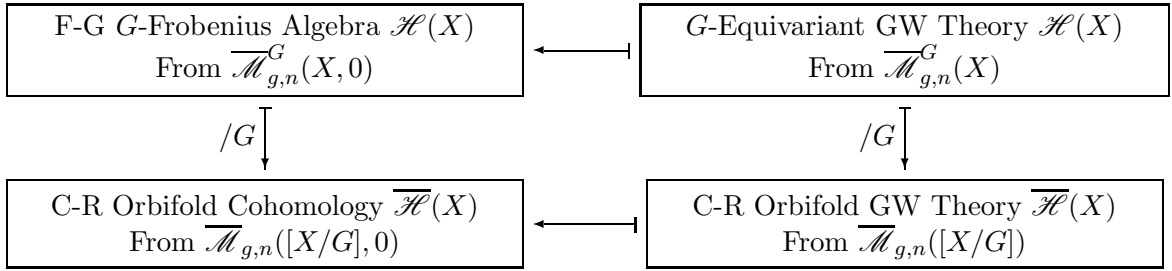

\[
\begin{diagram}
\fbox{\parbox{.40\linewidth}{\centerline{F-G $G$-Frobenius Algebra $\ch(X)$}
\centerline{From $\MM_{g,n}(X,0)$}}} & \lMapsto & &
\fbox{\parbox{.41\linewidth}{\centerline{$G$-Equivariant GW Theory $\ch(X)$}
\centerline{From $\MM_{g,n}(X)$}}} \\
\dMapsto^{/G} &                                  &  & \dMapsto^{/G}
\\
\fbox{\parbox{.40\linewidth}{\centerline{C-R Orbifold Cohomology $\chb(X)$}
\centerline{From $\M_{g,n}(\xg ,0)$}}} & \lMapsto & &
\fbox{\parbox{.40\linewidth}{\centerline{C-R Orbifold GW Theory $\chb(X)$}
\centerline{From $\M_{g,n}(\xg )$}}}
\\
\end{diagram}
\]
\caption{Schematic of Question (\ref{qu:two}) where each box contains an
algebraic structure and the relevant moduli space, F-G denotes
Fantechi-G\"ottsche, C-R denotes Chen-Ruan, the horizontal arrows denote
restriction, and the vertical arrows denote  taking ``quotients'' by $G$.}
\end{figure}

This paper provides affirmative answers to both of these questions.  

The first part of this paper is devoted to answering the first
question. We introduce $\MM_{g,n}$, the moduli space of $n$-pointed
admissible $G$-covers of genus $g$.  Roughly speaking, it consists of
a tuple $(E\rTo^\pi C;\pt_1,\ldots,\pt_n)$, where $E$ and $C$ are (at
worst, nodal) curves, $(C,p_1,\ldots,p_n)$ is a stable curve of genus
$g$, where $\pt_i$ are points in $E$ and $p_i := \pi(\pt_i)$, and
$\pi$ maps nodes of $E$ to nodes of $C$. Furthermore, we require that,
away from $\pi^{-1}(p_i)$ and nodes, $E$ is a principal $G$-bundle;
however, $E$ is allowed to have ramification over the marked points
and nodes. Our construction differs from the stack of admissible
covers in \cite{ACV}, as we require the additional data of $\pt_i$ in
$E$ over each marked point $p_i$ in $C$.

By forgetting the data associated to the $G$-cover, one obtains
a morphism $\st:\MM_{g,n}\rTo \M_{g,n}$, where $\M_{g,n}$ is the
moduli space of stable curves. We prove that $\MM_{g,n}$ is a
smooth, Deligne-Mumford stack, flat, proper, and quasi-finite (but not
representable) over
$\M_{g,n}$. Furthermore, $\MM_{g,n}$ has an action of the
symmetric group $S_n$ by permuting the ordering of the marked
points, and it has an action of $G^n$ by translation of the marked
points. In fact, $\MM_{g,n}$ admits an action of the braid group
$B_n$, which factors through the $S_n$ and $G^n$ actions.

The collection $\{ \MM_{g,n} \}$ possesses gluing morphisms,
provided that the monodromies of the two marked points to be glued
together are inverses of one another. These gluing morphisms are
equivariant under the action of $S_n$ and $G^n$. One may regard
the collection $\{ \MM_{g,n} \}$ as a $G$-equivariant colored
modular operad, where the coloring is by elements of $G$.
Furthermore, the morphism $\st$ respects the $S_n$ actions and the gluing
morphisms.

A \Gcft\ is defined analogously to a \cft, but where the role of
$\M_{g,n}$ is replaced by $\MM_{g,n}$, and where $G$-equivariance is
maintained throughout the construction.  We prove that there is an
external tensor product and a (usual) tensor product associated to
equivariant \cfts.  We then define the correct notion of taking a
``quotient'' by $G$ and prove that this procedure has the desired
properties.  The procedure of taking quotients involves an
intermediate step on the stack $\M_{g,n}(\BG)$ of stable maps into the
classifying stack $\BG$ (i.e., the stack of admissible covers without
the additional points $\pt_i$).  We show that in this intermediate
step the stack $\M_{g,n}(\BG)$ can be replaced by the quotient
$[\MM_{g,n}/G^n]$, but that the resulting ``quotient" \cfts\ are
isomorphic.

The last part of this paper will treat the second
question. We introduce the moduli space of $G$-stable maps
$\MM_{g,n}(X)$ and describe the $G$-equivariant Gromov-Witten
invariants. By restricting to contributions from
$\MM_{g,n}(X,0)$ alone, we prove that the state space $\ch(X)$ inherits
a $G$-Frobenius algebra structure which agrees with that from
\cite{FG}, and in particular that the trace axiom holds for their
ring. The proof consists of relating the virtual fundamental class
to an analogous cohomology class in their construction.

The details of the construction of a \Gcft\ for general equivariant
Gromov-Witten invariants, properties of potential functions, and
applications to higher spin curves will be explored elsewhere
\cite{IP}.

The Gromov-Witten invariants of orbifolds which are global quotients
of a variety by a finite group are particularly interesting in light
of the results of Costello \cite{Cos} which state that the Gromov-Witten
invariants of a smooth, projective variety $X$ of arbitrary genus are
determined by the genus zero Gromov-Witten invariants of the orbifolds
$[X^n/S_n]$ where $S_n$ is the symmetric group acting upon $X^n$  by
permuting its factors. We expect that our generalization of \cite{FG}
to higher degree stable maps will be useful in calculating these
invariants.

Finally, we observe that orbifolding plays an important role in mirror
symmetry, in certain Landau-Ginzburg theories (see, for example,
\cite{CoKa,Ma}), and in conformal field theory. In particular, there
are related notions of orbifolding which appear in the context of
vertex algebras (see, for example, \cite{Ki,FrSz}).  Furthermore, a
variant of our moduli spaces is used in the announcement \cite{Ki2} of
the construction of a modular functor associated to a finite group,
and this can be regarded as an example of an orbifold conformal field
theory.  It  would be enlightening to further clarify the relationship
between these notions.

The outline of this paper is as follows.  In Section \ref{2} we
describe the moduli spaces $\MM_{g,n}$, their associated forgetful
and gluing morphisms, group actions, and automorphism groups.  In
Section \ref{3} we briefly review important facts from the
category of $G$-graded $G$-modules, including the braid group
action and tensor products. In Section \ref{4} we define
\Gcfts\ and their tensor products. We prove that a
(non-projective) \Gcft\ always contains a $G$-Frobenius algebra.
In Section \ref{5} we define how to obtain a \cft\ from a \Gcft\
by taking the appropriate ``quotient.''  We prove that this is
consistent with the obvious notion of taking a quotient for a
$G$-Frobenius algebra, after rescaling the metric, and then work
out the example of the orbifold cohomology of $\BG$. In Section
\ref{6} we introduce the moduli space of $G$-stable maps and
equivariant Gromov-Witten invariants, reproduce the ring of
\cite{FG} as a special case, and prove that the trace axiom is
satisfied.

\begin{rem}
Unless otherwise specified, we assume that all cohomology rings
are over the ground ring $\nc$,
although all constructions here
are also valid over the rationals $\nq$.  

Also, unless otherwise specified, all groups which appear are finite and all
group actions are right group actions.
\end{rem}

\begin{nota}
  The stack quotient of a variety $X$ by $G$ will be denoted $\xg $
  and the coarse moduli space of this quotient will be denoted $X/G$.
\end{nota}

\begin{ack}
  We would like to thank Heidi Jarvis for help with typesetting and
  proofing.  We are grateful to Dan Abramovich for his careful reading
  of the manuscript and many helpful suggestions and corrections.  We
  also thank Jim Stasheff for his help with improving the exposition
  of the manuscript.
  
  The first author would like to thank Boston University for its
  hospitality.  The second author would like to thank the
  Max-Planck-Institut f\"ur Mathematik and the Institut des Hautes
  \'Etudes Scientifiques for their hospitality.  And the third author
  would like to thank the Institute for Advanced Study, where much of
  this work was done, for their financial support under National
  Science Foundation agreement No.~DMS-9729992.  He would also like to
  thank the Institut des Hautes \'Etudes Scientifiques for their
  hospitality where the finishing touches were applied to the paper.
  Any opinions, findings and conclusions or recommendations expressed
  in this material are those of the authors and do not necessarily
  reflect the views of the National Science Foundation.
\end{ack}
\section{The moduli spaces}\label{2}

\ 

Let $(C\rTo^{\varpi}T, p_1, \dots, p_n)$ be a stable curve over $T$ of
genus $g$, with marked points (sections) $p_1, \linebreak[2]\dots,
p_n$.  We want to study a variant of the space of admissible
$G$-covers of $C$, as defined in \cite[Def 4.3.1]{ACV}.  We recall the
definition here:
\begin{df}
A finite morphism $\pi: E \rTo C$ to an $n$-pointed, genus-$g$, stable
curve $C\rTo^{\varpi}T, p_1,\linebreak[2]\dots, p_n$ over $T$ is an
\emph{admissible $G$-cover} if
\begin{enumerate}
\item $E/T$ is itself a nodal curve (not necessarily connected).
\item Nodes of $E$ map to nodes of $C$.
\item There is a right action $\rho_E$ of $G$ on $E$ preserving $\pi$, and
such that
\item the restriction of $\pi$ to $C_{gen}$(the points of $C$
which are neither marked points nor nodes) is a principal
$G$-bundle.
\item At points of $E$ lying over nodes of $C$ the structure of the
maps $E\rTo^{\pi} C \rTo^{\varpi} T$ is locally the same as
(analytically isomorphic to) that of $$ \spec A[z,w]/(zw-t) \rTo \spec
A[x,y]/(xy-t^r) \rTo \spec A,$$ where we have $t\in A$, $x=z^r$ and
$y=w^r$, for some integer $r>0$.
\item At points of $E$ lying over marked points of $C$ the structure
of the maps $E\rTo^{\pi} C \rTo^{\varpi} T$ is locally the same as
(analytically isomorphic to) that of $$ \spec A[z] \rTo \spec A[x]
\rTo \spec A,$$ where $x=z^s$ for some integer $s>0$.
\item The action of the stabilizer $G_q\subseteq G$ at each node $q$
of $E$ is \emph{balanced}; that is, the eigenvalues of the action on
the tangent space at $q$ are multiplicative inverses of each other.
\end{enumerate}
\end{df}

Theorem 4.3.2 of \cite{ACV} shows that the stack of admissible
$G$-covers is isomorphic to the stack $\mgnbar(\BG)$ of balanced
twisted stable maps into the classifying stack of $G$.

\subsection{Definition, construction, and basic properties of $\MM_{g,n}$}

\

Given an admissible $G$-cover $(E \rTo^{\pi} C,\linebreak[0] p_1,
\dots, p_n)$, let ${\pt}_i \in \pi^{-1}(p_i)$ be a choice of a
point in the fiber over $p_i$ for all $i=1\dots n$.

\begin{df}
Let $\MM_{g,n}$ denote the stack of admissible $G$-covers $$(\pi:E
\rTo \linebreak[0] C, p_1,\dots,p_n, \pt_1, \dots, \pt_n)$$ of
$n$-pointed, genus-$g$, stable curves, together with a choice of $n$
marked points $\pt_i \in E$ such that $\pi(\pt_i) = p_i$ for all
$i=1,\ldots,n$. We call such objects \emph{$n$-pointed admissible
  $G$-covers}. A morphism of such objects is a $G$-equivariant fibered
diagram; that is, a morphism of the underlying stable curves, together
with a $G$-equivariant morphism of the induced admissible $G$-covers
preserving the points $\pt_i$.
\end{df}

Because the curve $C$ is oriented, a pointed admissible $G$-cover $(E
\rTo^{\pi} C, {\pt}_1, \dots, {\pt}_n)$ has a well-defined monodromy
$m_i$ at each marked point ${\pt}_i$; namely, $E$ induces a principal
$G$-bundle over $C-\{p_1,\dots,p_n\}$, and the orientation gives a
small loop in $C-\{p_1,\dots,p_n\}$ around each $p_i$, with a lift to
a path in a small neighborhood of ${\pt_i}$ in
$E-\{{\pt_1}, \dots, {\pt_n}\}$.  The lift is not uniquely
determined, but the difference between the starting and ending sheets
of the lifted path is given by a well-defined element $m_i\in G$.

Since the points $\pt_i$ are determined up to a discrete choice by
$C$, $\pi$, and the points $p_i$, the monodromy is invariant under
deformation of the curve $C$, the cover $E$, and the points $p_i$.
Also note that, while the action $\rho_E$ acts on the points $\pt_i$
by right multiplication, it acts on the holonomies by conjugation.

Let $G_A$ denote the set $G$, considered as a right $G$-space
under conjugation.  Associated to any object $(E \rTo^{\pi} C,
\pt_1, \dots, \pt_n)$ there exists an element $\bm=(m_1, \dots,
m_n) \in G^n_A$; namely, $m_i$ is the monodromy of $E$ at the
point $\pt_i$.

\begin{df}~\label{df:e}
Denote the canonical morphism we have just described by
\begin{equation}\label{eq:mm}
\be:\MM_{g,n}\rTo G_A^n, \end{equation} and let

\begin{equation}\label{eq:MMeval}
\MM_{g,n} (\bm):=\be^{-1}(\bm)
\end{equation}
denote the substack of objects in $\MM_{g,n}$ that map to $\bm$.
\end{df}
Since $\be$ is locally constant, we may write $$\MM_{g,n}=
\coprod_{\bm \in G^n_A}\MM_{g,n}(\bm).$$

The stack $\MM_{g,n}$ and the substacks $\MM_{g,n}(\bm)$ can be
explicitly constructed as follows.
\begin{thm}\label{thm:MMDM}
The stack $\MM_{g,n}$ and the substacks $\MM_{g,n}(\bm)$ are
smooth Deligne-Mumford stacks, flat, proper, and quasi-finite over
$\mgnbar$.
\end{thm}

\begin{proof}
Let $Adm_{g,n}^G$ be the stack of admissible $G$-covers of
$n$-pointed, genus-$g$ curves, and let $E\rTo^{\pi} C \rTo^{\varpi}
Adm_{g,n}^G$ be the universal $G$-cover and stable curve, with
universal gerbe markings $\mathscr{S}_i \rTo \cc:=[E/G]$. Let $$E_i:=E
\times_{\cc} \mathscr{S}_i$$ be the fibered product of $E$ with
$\mathscr{S}_i$. Let $$W:=E_1 \times_{Adm_{g,n}^G} E_2
\times_{Adm_{g,n}^G}\cdots \times_{Adm_{g,n}^G} E_n$$ be the fibered
product of the $E_i$. It is straightforward to see that $W$ is the
stack of admissible $G$-covers, together with explicit choices of
sections $\pt_i \in E$ lying over the sections $p_i$; that is,
$W=\MM_{g,n}$.

Theorems 3.0.2 and 4.3.2 of \cite{ACV} show that the space
$Adm_{g,n}^G$ is isomorphic to $\M_{g,n}(\BG)$, the stack of balanced
twisted stable maps to the classifying stack $\BG$, and is a smooth DM
stack, flat, proper, and quasi-finite over $\mgnbar$.  Since the
$\mathscr{S}_i$ are \'etale over $Adm_{g,n}^G$ and $E$ is \'etale over
$\cc$, these properties are preserved by the above-listed fibered
products.  Thus the theorem follows for $\MM_{g,n}$. The substacks
$\MM_{g,n}(\bm)$ are finite disjoint unions of connected components of
$\MM_{g,n}$, so the theorem also holds for them.
\end{proof}

\begin{rem}
The above construction of the moduli stack requires the use of the
gerbe sections $\mathscr{S}_i$ rather than the coarse sections
$A_i:=\im(p_i)$ in the coarse curve $C$.  This is due to the fact that
the fibered product of the $A_i$ with $E$ over $C$ does not
necessarily represent \emph{reduced} points of $E$---which is what we
really mean when we say a \emph{point}.
\end{rem}

\subsection{Morphisms and group actions on $\MM_{g,n}$}

\ 

There are several obvious morphisms on $\MM_{g,n}$.  First, there are
the forgetful morphisms $$\MM_{g,n}\rTo^{\stt} Adm^G_{g,n} \cong
\mgnbar (\BG) \rTo^{\sth} \mgnbar ,$$
which were shown in Theorem
~\ref{thm:MMDM} to be proper, flat, and quasi-finite. We denote the
composition by $$\st:=\sth\circ\stt.$$
We also have the evaluation
morphism (\ref{eq:MMeval}): $$\be:\MM_{g,n} \rTo G^n_A.$$

Recall (see \cite{JK}) that while $\M_{g,n}(\BG)$ cannot be written as
a disjoint union of substacks indexed by $\bm \in G^n_A$, it does have
a decomposition indexed by conjugacy classes of $G$.

\begin{df}
We denote the set of conjugacy classes of $G$ by $\Gb$ and the
conjugacy class of $m \in G$ by $\mb$.  Similarly, we denote by $\bmb
\in \Gb^n$ the $n$-tuple of conjugacy classes determined by $\bm \in
G^n$.
\end{df}

As described in \cite{JK}, we have $$\M_{g,n}(\BG)= \coprod _{\bmb \in
\Gb^n} \M_{g,n} (\BG; \bmb),$$ where some of the substacks may be
empty.

\begin{df}
We let $\MM_{g,n}(\bmb)$ denote the preimage
$\stt^{-1}(\M_{g,n}(\BG;\bmb)),$ which is easily seen to be
$$\MM_{g,n}(\bmb)=\coprod_{\bm' \in \bmb} \MM_{g,n}(\bm').$$
\end{df}

The stack $\MM_{g,n} (\bmb)$ has a right $G^n$ action
\begin{equation}\rho(\gamma_1,\ldots,\gamma_n):\MM_{g,n}(\bm)\rTo
\MM_{g,n}(\gamma_1^{-1} m_1\gamma_1,\ldots,\gamma_n^{-1}
m_n\gamma_n),\end{equation} which acts by right multiplication on the
$n$ marked points $(\pt_1,\ldots,\pt_n)\mapsto
(\pt_1\cdot\gamma_1,\ldots, \pt_n\cdot\gamma_n)$.  We sometimes write
$\varrho_i$ for the action on the $i$th factor: $\varrho_i(\gamma) =
\varrho(1,\dots, \gamma,\dots 1)$.

Together with the action of the symmetric group $S_n$ on $\MM_{g,n}$,
which reorders the marked points, $\MM_{g,n}$ has the action of the
semi-direct product group $G^n \rtimes S_n$, called the \textsl{wreath
  product},  where $S_n$ acts on $G^n$ by permuting its
factors. One consequence is that $\MM_{g,n}$ has the action of the
braid group $B_n$.

\begin{df}
Let $B_1$ be the trivial group. If $n\geq 2$, let $B_n$ be the group
with generators $\{ b_1,\ldots,b_{n-1}\}$ subject to the relations
\begin{equation}\label{eq:braid1}
b_i b_{i+1} b_i = b_{i+1} b_i b_{i+1}
\end{equation}
for all $i=1,\ldots, n-1$ and
\begin{equation}\label{eq:braid2}
b_i b_j =b_j b_i
\end{equation}
if $|i-j|>1$.
$B_n$ is called the \emph{braid group on $n$-letters}.
\end{df}

For each generator $b_i$ of the braid group $B_n$, there is an isomorphism
\begin{equation}\label{eq:braidact}
\MM_{g,n}(m_1,\ldots,m_i,m_{i+1},\ldots,m_n)\rTo^{b_i}
\MM_{g,n}(m_1,\ldots,m_i m_{i+1} m_{i}^{-1},m_i,\ldots,m_n).
\end{equation}
  These are given by $b_i := \rho_i(m_i^{-1})\circ s_i$, where $s_i$
is the element $(i,i+1)$ in $S_n$ which transposes $i$ and $i+1$, and
$\rho_i$ is the group action on $\MM_{g,n}$ obtained by right
translation of the $i$-th marked point.  It is straightforward to
check that the induced isomorphisms satisfy the braid relations
(\ref{eq:braid1}) and (\ref{eq:braid2}), thus they induce an action of
$B_n$ on $\MM_{g,n}$.

Finally, there are the three fundamental morphisms:
\emph{forgetting tails}, \emph{gluing trees}, and \emph{gluing
loops}.

\begin{description}
\item[Forgetting Tails]

Let $\bm$ be any $n$-tuple $$\bm:=(m_1, \dots, m_n) \in G^n_A,$$
and let $1$ be the identity in $G$. Whenever the pair $(g,n)$ is
stable (i.e., $2g-2+n> 0$) there is a natural
\emph{forgetting tails} morphism $\tst:\MM_{g,n+1}(\bm, 1) \rTo
\MM_{g,n}(\bm)$ defined as follows.

First, simply forgetting the data associated to the $(n+1)$st
marked point usually yields an object of $\MM_{g,n}(\bm)$, but if
the resulting curve is unstable, then we need to contract the
unstable component to a point $p$.  In those cases it is true, but
not immediately obvious, that we can produce a suitable $G$-cover
$E$ on the new curve, and where necessary, assign a point $\pt$ in
$E$ over $p$. We now describe how this works.

We have two cases:  first, when the resulting unstable component
$D$ is a (genus-zero) $-1$-curve with one marked point $p_i$, and
one node $q$; and second, when the unstable component $D$ is a
$-2$-curve with two nodes $q$ and $q'$ and no marked points.

In either case, the unstable component $D$ is a genus-zero curve
with two special points (call them $q$ and $q'$ for
simplicity of notation).    It is straightforward to see that for
any $\widetilde{q}' \in E$ over $q'$ with monodromy $m$, the connected
component $\widetilde{D}$ of $E$ containing $\widetilde{q}'$ is a finite
cover of $D$ with automorphism group
$\operatorname{\Aut}_D\widetilde{D}$ generated by $m$, which acts
faithfully on all points but $q$ and $q'$. In particular, it is
fully ramified over $q$ and $q'$, and unramified at all other
points. Thus there is a canonical $G$-equivariant isomorphism
$\varphi:E|_q \rTo^{\sim} E|_{q'}$. This shows in the first case,
where $p_i=q'$, that there is a canonical choice of $\widetilde{q} \in
E|_q$ (namely, $\widetilde{q}=\varphi^{-1}({\pt}_i)$), and thus a
well-defined point of $\MM_{g,n}(\bm)$.

In the second case, the isomorphism $\varphi$ allows the
construction of a principal $G$-bundle on the curve with the
unstable component $D$ contracted.  In this case, we need no point
$\widetilde{q}$---the data we already have will give a point of
$\MM_{g,n}$. Thus in every case the forgetting tails morphism
exists.

\item[Gluing Trees]
Given any $\bm \in G^{n_1}_A$ and $\bm' \in G^{n_2}_A$, as well as
an additional element $\mu \in G_A$, let $g:=g_1+g_2$ and
$n:=n_1+n_2$. We have the \emph{gluing trees} morphism:
\begin{equation}\label{eq:gluetree}
\varrho_{tree}:\MM_{g_1,n_1+1}(\bm,\mu) \times
\MM_{g_2,n_2+1}(\mu^{-1},\bm') \rTo
\MM_{g,n}(\bm,\bm')
\end{equation}
given by attaching the universal $G$-covers $E\rTo^{\pi}C'
\rTo^{\varpi}\MM_{g_1,n_1+1}(\bm,\mu)$ and
$E'\rTo^{\pi'}C\linebreak[0]\rTo^{\varpi'}\MM_{g_2,n_2+1}(\mu^{-1},\bm')$
along the sections $\rho(\gamma)\pt_{n_1+1}\in E$ and
$\rho(\gamma)\pt'_1\in E'$ for all $\gamma$ in $G$, and attaching the
universal curves $C$ and $C'$ along the sections $p_{n_1+1}$ and
$p'_1$.  It is straightforward to see that, because the monodromies
$\mu$ and $\mu^{-1}$ are inverses, the induced cover is indeed an
admissible $G$-cover of the resulting stable curve, and thus gives an
object in $\MM_{g,n}(\bm,\bm')$.

More generally, let  $I = \{i_1, \dots, i_{n_1} \}$ and $J= \{j_1,
\dots, j_{n_2} \}$  be any disjoint subsets of  $\{1,\dots,n\}$
such that $I \sqcup J=\{1,\dots,n\}$. For any integers $s,t$ with
$i \leq s \leq n_1,$ $1 \leq t \leq n_2$ there is a morphism
\begin{equation}\label{eq:glue-perm}
 \MM_{g_1,n_1+1}(m_{i_1},\dots, m_{i_{s-1}}, \mu, m_{i_{s}},
\dots, m_{i_{n_1}})
\end{equation}
$$\times
\MM_{g_2,n_2+1}(m_{j_1},\dots, m_{j_{t-1}},\mu^{-1},
m_{j_{t}},\dots, m_{j_{n_2}})
 \rTo $$  $$\MM_{g,n}(m_1, \dots, m_n).$$

\item[Gluing Loops]
Given any $\bm \in G^n_A$ and $\mu \in G_A$ we have the
\emph{gluing loops} morphism:
\begin{equation}\label{eq:glueloop}
\varrho_{loop}:\MM_{g-1,n+2}(\bm,\mu,\mu^{-1}) \rTo
\MM_{g,n}(\bm),
\end{equation}
defined in a manner similar to the gluing trees morphism; namely,
one attaches the universal $G$-cover $E$ to itself along the two
sections $\pt_{n+1}$ and $\pt_{n+2}$, and the universal curve $C$
to itself along the sections $p_{n+1}$ and $p_{n+2}$.

As with gluing trees, the gluing loops morphism can be defined
more generally for any two sections $\pt_i$, and $\pt_{j}$,
provided they have inverse monodromies.

\end{description}

\begin{rem} Even more generally, if two points do not have inverse monodromies,
  the braid group action may still allow one to glue them.  For
  example, for any $i_1<i_2$ with $i_1,i_2 \in \{0,\dots,n+1\}$ and
  $\sigma = m^{-1}_{i_1+1}m^{-1}_{i_1+2}\dots m^{-1}_{i_2}
  \mu^{-1}m_{i_2}\dots m_{i_1+1}$, we have a morphism
\begin{align*} & \MM_{g,n+2}(m_1,\dots, m_{i_1}, \mu, m_{i_1+1}, 
\dots, m_{i_2}, \sigma, m_{i_2+1}, \dots, m_{n+2}) \\ &
\rTo^{b_{i_1+1} \circ b_{i_1+2}\circ \dots \circ b_{i_2}}
\MM_{g,n+2}(m_1,\dots, m_{i_1}, \mu, \mu^{-1}, m_{i_1+1}, \dots, m_n)
\\ & \rTo^{\varrho_{loop}} \MM_{g+1,n} (m_1, \dots, m_n).
\end{align*}
\end{rem}

\begin{rem}
Since the collection $\{ \MM_{g,n} \}$ has gluing morphisms which
are equivariant under the actions of $G^n$ and $S_n$, one may
regard $\{ \MM_{g,n} \}$ as a $G$-equivariant colored modular
operad where the set of colors is the $G$-set $G_A$.  Since the
action of the braid group $B_n$ (see Equation~(\ref{eq:braidact}))
on $\MM_{g,n}$ is constructed from the $G^n$ and $S_n$ actions,
one may also regard $\{ \MM_{g,n} \}$ as a colored modular operad,
but where the role of the permutation group is replaced by the
braid group.
\end{rem}

\begin{rem}
  It is worth pointing out that the stack $\MM_{g,n+1}(\bm, 1)$ is not
  the universal curve or orbicurve over $\MM_{g,n}(\bm)$ nor is it the
  universal admissible $G$-cover.  On the one hand, generic locations
  of $\pt_{n+1}$ will have no automorphisms, since they must fix the
  point $\pt_{n+1}$.  On the other hand, when $p_{n+1}$ ``collides''
  with another marked point (i.e., they bubble off a genus-zero
  component), then the point $\pt_{n+1}$ only prevents the existence
  of non-trivial automorphisms of $E$ over the new component, but
  automorphisms over the remainder of the curve need only fix the
  fiber of $E$ over the new node.  \end{rem}

\subsection{Holonomy and other tools for studying $G$-covers}

\ 

Let $G_R$ denote $G$ considered as a right $G$-module.  Note that
the automorphism group $\Aut^G(G_R)$ of $G_R$ is again $G$, acting
by \emph{left} multiplication. Given a pointed admissible cover
$(E\rTo^{\pi} C, \pt_1, \dots,\pt_n)$ and given any point $\pt_0
\in E_{gen}:=\pi^{-1}(C_{gen})$ lying over $p_0 \in C_{gen}$, we
have an isomorphism of  right $G$-modules $\nu_{\pt_0} :E|_{p_0}
\irightarrow G_R$, given by $$\nu_{\pt_0}(\pt_0\gamma) :=
\gamma.$$

Changing the base point $\pt_0$ to $\pt_0\alpha$ changes the map
$\nu_{\pt_0}$ by left multiplication by $\alpha^{-1}$.
$$\nu_{\pt_0\alpha}= \alpha^{-1}\nu_{\pt_0}.$$

\begin{df}\label{df:holonomy}
  The choice of $\pt_0 \in E_{gen}$ gives a homomorphism from the
  fundamental group to $G$: $$\chi_{\pt_0}:\pi_1(C_{gen},p_0) \rTo
  G,$$
  which we call \emph{holonomy}.  One way to see this
  homomorphism explicitly is to pull $E_{gen}$ back to the trivial
  admissible cover $\widetilde{E}$ of the universal cover $U$ of
  $C_{gen}$.  Automorphisms of $U$ are precisely $\pi_1(C_{gen},p_0)$,
  and they induce automorphisms of $\widetilde{E} \cong U \times G_R$,
  and therefore of $G_R$: $$
  \pi_1(C_{gen},p_0) \cong \Aut_{C_{gen}}U
  \rTo^{\chi} \Aut^G G_R =G.$$
  
  Conversely, given any homomorphism $\chi:\pi_1(C_{gen},p_0) \rTo G$,
  it is easy to see that we get a uniquely determined admissible
  $G$-cover of $(C, p_1, \dots, p_n)$ and a distinguished point
  $\pt_{0,\chi}$ over $p_0$.  This $G$-cover is given by first taking
  the quotient of $U \times G_R$ by the action of $$\pi_1(C_{gen},p_0)
  \rTo^{\mathrm{id},\chi} \Aut_{C_{gen}} U \times \Aut^G G_R$$
  and
  then extending it to all of $C$. Such an extension is uniquely
  determined by the $G$-cover on $C_{gen}$.  The point
  ${\pt}_{0,\chi}$ is the image of $(p_0,1) \in U \times G_R$ under
  this quotient.  We call this cover \emph{the admissible $G$-cover of
    $C$ induced by $\chi$} and $p_0$, and we denote it $E_{\chi,p_0}$,
  or $E_{\chi}$ if $p_0$ is clear from context.
\end{df}

The following proposition is an immediate consequence of well-known
corresponding results for principal $G$-bundles (see, for example,
\cite[Chapters 13--14]{Fu}) and is straightforward to check.
\begin{prop}
Let $C_{gen}$ be connected.  For any homomorphism
$\chi:\pi_1(C_{gen},p_0) \rTo G$, the induced $E$ and $\pt_{0,\chi}$
have holonomy $\chi_{\pt_0}$ equal to $\chi$, and conversely, given an
$E$ and $\pt_0$ the bundle $E_{\chi_{\pt_0}}$ is canonically
isomorphic to $E$, via an isomorphism identifying
$\pt_{0,\chi_{\pt_0}}$ with the original $\pt_0$.
Thus the data of $E,\pt_0$ is equivalent to a choice of
homomorphism $\chi:\pi_1(C_{gen},p_0) \linebreak[0] \rTo G$.

A different choice of point $\pt_0$, say $\pt_0 \alpha$, changes
$\chi$ by conjugation $\chi_{\pt_0 \alpha}=\alpha^{-1}\chi
\alpha$.  Furthermore, given a path $\gamma$ from $p_0$ to $p'_0$
in $C_{gen}$ and the corresponding unique lift $\widetilde{\gamma}$ of
$\gamma$ from $\pt_0$ to $\pt'_0 \in E|_{p'_0}$, the holonomy
$\chi_{\pt'_0}$ is induced from $\chi_{\pt_0}$ by conjugation with
$\gamma$.

\begin{diagram}
\pi_1(C_{gen},p_0) & \rTo^{\chi_{\pt_0}} & G\\ \dTo^{ad(\gamma)} &
\ruTo_{\chi_{\pt'_0}}&\\ \pi_1(C_{gen},p_0)\\
\end{diagram}

And conversely, given any $\chi':\pi_1(C_{gen},p_0) \rTo G$ determined
from $\chi$ by conjugation by $\gamma$, the induced $G$-cover
$E_{\chi'}$ is canonically isomorphic to $E_{\chi}$, and the induced
point $\pt'_{0,\chi'}$ is that obtained by parallel transporting
$\pt_0$ along $\gamma$ (i.e., $\pt'_{0,\chi}$ is the endpoint of
$\widetilde{\gamma}$).
\end{prop}

\begin{df}\label{df:path-loop}
For any path $d$ from $p_0$ to $p_i$ in $C_{gen}$ (that is, a path in
$C$ such that the image of $(0,1)$ lies in $C_{gen}$ and $d(1)=p_i$
and $d(0)=p_0$), we have an induced element $\sigma_{d}$ of
$\pi_1(C_{gen},p_0)$ defined by following $d$ from $p_0$ to a little
loop around $p_i$, tracing the loop out once counterclockwise, and
then returning along $d$ (or rather $d^{-1}$) to $p_0$.

Moreover, for any admissible $G$-cover with point $\pt_0 \in
E|_{p_0}$, the path $d$ determines a point $\pt(d) \in E|_{p_i}$,
which is the endpoint of the unique lift $\widetilde{d}$ of $d$ in $E$
that begins at $\pt_0$.

Finally, given $\pt_0$ and a path $d$ from $p_0$ to $p_i$,
holonomy and the map $\nu_{\pt_0}$ induce an isomorphism of right
$G$-modules $\overline{\nu}_{d,\pt_0}:E|_{p_i}\irightarrow \langle
m_i\rangle\backslash G_R$, where $m_i:=\chi_{\pt_0}(\sigma_d)$,
and $\bar{\nu}_{d,\pt_0}$ maps the point $\pt(d)$  to the coset
$\langle m_i\rangle$.
\end{df}

\begin{df} \label{df:zeta}
  In genus zero, a choice of paths $d_i$ from the point $p_0$ to the
  point $p_i$ for each $i\in\{1,\dots,n\}$ induces loops
  $\sigma_{d_i}$ that generate $\pi_1(C_{gen},p_0)$.  Not every choice
  of monodromy $\bm \in G^n$ satisfies the same relations that the
  generators $\sigma_i$ do, and thus not every choice of monodromy
  defines a holonomy $\chi$, but for those $\bm$ that do, there is a
  uniquely determined pointed admissible $G$-cover $$\zeta (d_1,
  \dots, d_n; \bm):=(E_{\chi}\rTo \CP^1, \pt_1,\dots,\pt_n)$$
  by
  defining the holonomy $\chi$ to be given by the monodromy
  $$\chi(\sigma_{d_i}) = m_i, \quad \mbox{ for $i\in \{1,\dots,n\}$}$$
  and letting the points $\pt_i:=\pt(d_i)$ be the points induced as in
  Definition~\ref{df:path-loop}.  Since the loops $\sigma_{d_i}$
  generate the fundamental group of $\CP^1-\{p_1,\dots,p_n\}$, this
  construction gives a well-defined pointed admissible $G$-cover.
\end{df}

It is clear from our discussion so far that every smooth, genus-zero,
$n$-pointed, admissible $G$-cover $(E\rTo \CP^1, \pt_1, \dots, \pt_n)$
that has all of its points $\pt_i$ in the same connected component of
$E$ must be of the form $\zeta(d_1, \dots, d_n; \bm)$ for some choice
of $p_0$, paths $(d_1, \dots, d_n)$, and $\bm\in G^n$.  Assume that
the points $p_0, \dots, p_n \in C$ are given.  Of special interest is
the case where the induced generators of the fundamental group have
product equal to $1$.  We denote the subset of such $n$-tuples of
paths by
\begin{equation}
P_C:=\{(d_1, \dots, d_n)| d_i \text{ a path from } p_0 \text{ to }
p_i, \text{ and } \prod_{i=1}^n \sigma_{d_i} = 1 \},
\end{equation} and the corresponding pointed admissible $G$-covers of $C$  by
 \begin{equation}
\zeta_C:= \{ \zeta(\bd;\bm)| \bd\in P_C, \ \bm\in G^n, \
\prod_{i=1}^n m_i = 1\}.
\end{equation}

\begin{df-pr}\label{df-prp:trnsBact}
Given a choice of $p_0,\dots, p_n \in C$, with
$\operatorname{genus}(C)=0$, there is a transitive action of the
braid group on the set $P_C$, where
\begin{align}
b_i d_i  &= d_{i+1} \\
b_id_{i+1} &= \sigma_{d_i+1} d_i\\
b_i d_j &= d_j \text{ if } j \neq i,i+1 .
\end{align}
This action of $B_n$ on the set $P_C$ is compatible with the usual
braid action on $\pi_1(C_{gen},p_0)$; that is, for each $i$ we have
$\sigma_{d_i} \in \pi_1(C_{gen},p_0)$,
and \begin{equation}\label{eq:braidtwo} \sigma_{bd_i} = b\sigma_{d_i}.
\end{equation}
Consequently, the braid action on $P_C$ induces an action of the braid
group on $\zeta_C$, distinct from the braid action on all of
$\MM_{g,n}$ that we defined earlier.  To distinguish the two, we will
denote this new action by $\beta:B_n \rTo \Aut(\zeta_C)$.
\end{df-pr}
\begin{proof}
The fact that the given equation defines an action and that the action
is compatible with the usual action on the fundamental group is a
straightforward calculation.  That the action is transitive follows
from the classical fact that the braid group generates all
outer automorphisms of the fundamental group that preserve the property of
the product of generators being trivial.

Since the product of generators and the product of
monodromies are both  trivial, the induced holonomy
$b\chi:\sigma_{(b d_i)}\mapsto m_i$ is still a well-defined
homomorphism of groups.  Thus for each admissible cover
$\zeta(\bd;\bm)\in \zeta_C$ and for each  $b\in B_n$ we can define
\begin{equation}\beta(b)\zeta(\bd;\bm):=\zeta(b\bd;\bm).\end{equation}
 \end{proof}

\subsection{Automorphisms, isomorphisms, and fibers}

\ 

\begin{df}
Let $\Aut^G_C E$ denote the group of $G$-equivariant automorphisms
of $E$ over $C$.  Any $\varphi \in \Aut^G_C E$ must induce a
$G$-equivariant automorphism $\varphi':G_R \rTo G_R$ of right
$G$-modules by $\varphi'=\nu_{\pt_0} \circ \varphi \circ
\nu^{-1}_{\pt_0}$.  It is easy to see that if $\varphi (\pt)=\pt_0
g$, then $\varphi'$ is simply left multiplication by $g$.  This
gives a homomorphism $$\Psi_{\pt_0}:\Aut^G_C E\rTo G.$$
\end{df}

\begin{prop}\label{prop:commute}
The homomorphism $\Psi_{\pt_0}:\Aut^G_CE \rTo G$ commutes with
every element of $\im(\chi_{\pt_0})$, and depends only upon the
(path-)component of $E_{gen}$ in which $\pt_0$ lies. Moreover, if
$C$ is irreducible, then $\Psi_{\pt_0}$ is an isomorphism to the
centralizer of (i.e., the subgroup of $G$ which commutes with every
element of) the image of $\chi_{\pt_0}$: $$\Psi_{\pt_0}:\Aut^G_CE
\irightarrow C(\im \chi_{\pt_0}).$$
\end{prop}

\begin{proof}
It is straightforward to check that a change of base point from
$\pt_0$ to $\pt'_0=\pt_0 \gamma$ changes $\Psi_{\pt_0}$ by
conjugation.  $$\Psi_{\pt_0
\gamma}=\gamma^{-1}\Psi_{\pt_0}\gamma.$$

On the other hand, given a path $\sigma:[0,1]\rTo E_{gen}$ from
$\pt_0$ to another point $\widetilde{q}_0$ we may parallel transport
any point $\pt_0 \gamma$ of the fiber $E|_{p_0}$ to the point
$\widetilde{q}_0 \gamma$ in the fiber $E|_{q_0}$, thus giving an
isomorphism of right $G$-sets $\sigma_*:E|_{p_0}\irightarrow E|_{q_0}$,
and one can check that the induced homomorphisms $\Psi_{\pt_0}$
and $\Psi_{\widetilde{q}_0}$ are the same:
$$\Psi_{\pt_0}=\Psi_{\widetilde{q}_0}:\Aut^G_CE \rTo G.$$ The first
two claims of the proposition follow.

It is straightforward to check that if $C_{gen}$ is path
connected, then $\Psi_{\pt_0}$ is injective, and surjectivity can
be seen by uniformizing $C_{gen}$, pulling $E$ back to a trivial
bundle on the uniformizer, and checking that left multiplication
by any element of $G$ which commutes with holonomy descends to a
$G$-equivariant automorphism of $E$ over $C$.
\end{proof}

We now turn our attention to automorphisms of \emph{pointed} admissible
$G$-covers. For a pointed admissible $G$-cover ($E\rTo^{\pi} C,
\pt_1, \dots, \pt_n$), we denote the group of $G$-equivariant
automorphisms of $E$ over $C$ which fix the points $\pt_i$ by
$\Aut^G_C(E, \pt_1,\dots, \pt_n)$.

\begin{prop}\label{prop:mono}
If $C$ is an irreducible curve, and if $m_1,\dots, m_n \in G_A$
are the monodromies of the admissible $G$-cover E at $\pt_1,
\dots, \pt_n$, respectively, then for any elements
$\gamma_1,\dots, \gamma_n$ such that $\pt_0 \in E_{gen}$ lies in
the same connected component of $E_{gen}$ as $\pt_1
\gamma_1,\dots, \pt_n \gamma_n$, the map $\Psi_{\pt_0}$ induces an
isomorphism $$\Psi_{\pt_0}:\Aut^G_C(E,\pt_1,\dots,\pt_n) \irightarrow
\langle\gamma^{-1}_1 m_1 \gamma_1\rangle \cap \dots \cap
\langle\gamma^{-1}_n m_n \gamma_n\rangle \cap C(\im(\chi_{{\pt}_0})),$$
where $C(\im(\chi_{{\pt}_0}))$ denotes the centralizer of the image of
$\chi_{{\pt}_0}$.
\end{prop}

\begin{proof}
If $\pt_i \gamma_i$ is in the same component of $E_{gen}$ as
$\pt_0$, then there is a path $d$ in $C_{gen}$ from $p_0$ to $p_i$
which lifts to a path $\widetilde{d}$ from $\pt_0$ to $\pt_i$, and we
have an isomorphism $\overline{\nu}_{d,\pt_0}:E|_{p_i}
\irightarrow \langle\gamma^{-1}_i m_i \gamma_i\rangle\backslash
G_R $ of right $G$-sets taking $\pt_i \gamma_i$ to the coset
$\langle\gamma^{-1}_i m_i \gamma_i\rangle$. An automorphism
$\varphi \in \Aut^G_C E$ with $\Psi_{\pt_0}(\varphi) =g$ takes the
coset $\langle\gamma^{-1}_i m_i\gamma_i\rangle$ to itself if and
only if $g \in \langle\gamma^{-1}_i m_i\gamma_i\rangle$.  Thus
$\varphi$ fixes the points $\pt_i \gamma_i$ and also $\pt_i$ if
and only if $\Psi_{\pt_0} (\varphi) \in \langle\gamma^{-1}_i m_i
\gamma_i\rangle$ for every $i$.
\end{proof}

Of course, if $\pt_i \gamma$ is in the same connected component of
$E_{gen}$ as $\pt_i \alpha$, then, since $\Psi_{\pt_0}(\varphi)$
commutes with holonomy, including $\gamma^{-1} \alpha$, the
condition $\Psi_{\pt_0}(\varphi) \in \langle\gamma^{-1} m_i
\gamma\rangle$ is the same as the condition $\Psi_{\pt_0}(\varphi)
\in \langle\alpha^{'-1} m_i \alpha\rangle$.

\begin{prop}\label{prop:central}
For any smooth pointed curve $(C,p_1,\dots, p_n)$ having \emph{no
  non-trivial automorphisms}, choose an admissible cover $(E
\rTo^{\pi}C, p_1,\dots, p_n) \in Adm^G_{g,n}$. For any ${\pt}_0
\in E_{gen}$ and for any choice of paths $d_i$ in $C_{gen}$ from $p_0
= \pi(\pt_0)$ to $p_i$, let $\sigma_i=\sigma_{d_i}$ be the
corresponding element of $\pi_1(C_{gen},p_0)$.  We can describe the
fiber $(\stt)^{-1}([E \rTo^{\pi}C,p_1,\dots,p_n])$ of the forgetful
map $\stt: \MM_{g,n}\rTo Adm^G_{g,n}$ as the quotient
stack $$(\stt)^{-1} (E\rTo^{\pi}C,p_1,\dots,p_n)=
\left[\left(\prod^n_{i=1} \langle\chi_{\pt_0} (\sigma_i)
  \rangle\backslash
  G_R\right)/C(\chi_{\pt_0})\right]=\coprod_{I_{\pt_0}}\mathscr{B}H_{\pt_0},$$
where $C(\chi_{\pt_0})$ is the centralizer of the image of
$\chi_{\pt_0}$, acting diagonally on the product, the index set
$I_{\pt_0}$ is
$\left(\prod^n_{i=1}\langle\chi_{\pt_0}(\sigma_i)\rangle\backslash
G_R\right)/\left(C (\chi_{\pt_0})/H_{\pt_0}\right)$, and the group
$H_{\pt_0}$ is the image under $\Psi_{\pt_0}$ of the automorphism
group of any pointing $(\pt_1,\dots,\pt_n)$ of $E$:
$$H_{\pt_0}=\Psi_{\pt_0}(\Aut^G_C(E,\pt_1,\dots,\pt_n))=C(\chi_{\pt_0})
\cap \langle\chi_{\pt_0}(\sigma_1)\rangle\cap \dots \cap
\langle\chi_{\pt_0}(\sigma_n)\rangle.$$
\end{prop}

\begin{proof}
A choice of pointing $\pt_1,\dots, \pt_n \in E$ is equivalent to a
choice
$\overline{\nu}_{\pt_0}(\pt_i)\in\langle\chi_{\pt_0}(\sigma_i)\rangle
\backslash G_R$ for each $i \in \{1,\dots, n\}$, and any
isomorphism between two pointings $(E,\pt_1,\dots,\pt_n)$ and
$(E,\pt'_1,\dots,\pt'_n)$ induces an automorphism of $E$.
Conversely, the automorphisms of $E$ act on the set of all
pointings, thus Proposition ~\ref{prop:commute} gives the first
equality.  For any pointing, the homomorphism $\Psi_{\pt_0}$ takes
the automorphism group $\Aut^G_C (E, \pt_1,\dots, \pt_n)$ to
$H_{\pt_0}:=C(\chi_{\pt_0}) \cap
\langle\chi_{\pt_0}(\sigma_1)\rangle\cap \dots, \cap
\langle\chi_{\pt_0}(\sigma_n)\rangle$ by Proposition
~\ref{prop:mono}.  The second equality follows.
\end{proof}

\begin{prop}
Let $C=C^1 \cup C^2$ be the union of two irreducible curves joined
at a single node $q$.  Choose points $\pt^1_0, \pt^2_0 \in
E_{gen}$ lying over $C^1_{gen}$ and $C^2_{gen}$, respectively, and
such that $\pt^1_0$ and $\pt^2_0$ lie in the same connected
component of $E$.  Also, choose a point $\widetilde{q} \in E|_q$ of the
fiber over $q$ which lies in the same connected component of $E$
as $\pt^1_0$ and $\pt^2_0$.  Let $\mu$ and $\mu^{-1}$ be the
monodromy of $E$ at $q$ with respect to the orientations of $C^1$
and $C^2$ respectively.  We have an injective homomorphism
$$\Psi:=(\Psi_{\pt^1_0}, \Psi_{\pt^2_0}):\Aut^G_C(E,
\pt_1,\dots,\pt_n) \hookrightarrow G \times G,$$ which depends
only on the connected component of $E$ in which $\pt^1_0$ and
$\pt^2_0$ lie, and $$\im \Psi=\{(g_1,g_2) \in
\im\Psi_{\pt^1_0}\times \im \Psi_{\pt^2_0}| g_1g^{-1}_2 \in
\langle \mu_1\rangle\}.$$
\end{prop}

\begin{proof}
The injectivity follows from arguments similar to the irreducible
case.  The condition on the elements $(g_1,g_2) \in \im
\Psi_{\pt^1_0} \times \im \Psi_{\pt^2_0}$ comes from the fact that
any automorphism of $E$ must take both ``sides'' of the node
$\widetilde{q}$ to the same point:  $\widetilde{q}g_1 =\widetilde{q}g_2$, but
$\widetilde{q}g_i$ is only determined up to a (left) coset of
$\langle \mu \rangle$.
\end{proof}

Let $C$ be an irreducible curve with one node $q$ obtained by
attaching $2$ points $q_+$ and $q_-$ of the normalized curve
$C^{\nu}$.  An admissible $G$-cover $E$ of $C$ is obtained by
attaching two points $\widetilde{q}_+ \in E^{\nu}|_{q_+}$ and
$\widetilde{q}_- \in E^{\nu}|_{q_-}$ of an admissible $G$-cover
$E^{\nu}$ on $C^{\nu}$ which have monodromy $\mu$ and $\mu^{-1}$,
respectively, for some $\mu \in G$.  Let $\pt_0 \in
E^{\nu}_{gen}=E_{gen}$ be in the same connected component of
$E^{\nu}$ as $\widetilde{q}_+$ is, and let $\gamma \in G$ be chosen so
that $\widetilde{q}_- \gamma^{-1}$ is in that same component of
$E^{\nu}$.  

\begin{prop}\label{prp:aut-slf-sew}
Any automorphism $\varphi \in \Aut^G_C(E,\pt_1,\dots, \pt_n)$
induces an automorphism $N(\varphi) \in
\Aut^G_{C^{\nu}}(E^{\nu},\pt_1,\dots, \pt_n)$ by pullback to the
normalization.  For any $\pt_0 \in E_{gen}$ the homomorphism
$N$ is injective and is compatible with $\Psi_{\pt_0}$; that
is, the following diagram commutes:

\begin{diagram}
\Aut^G_C(E, \pt_1, \dots, \pt_n) & \rTo^{\Psi_{\pt_0}} & G\\
\dInto^{N} & & \dEq \\ \Aut^G_{C^{\nu}}(E^{\nu},\pt_1,\dots, \pt_n)&
\rTo^{\Psi_{\pt_0}} & G
\end{diagram}

Moreover, we have $$\im(\Psi_{\pt_0} \circ N)=\{g \in
\Psi_{\pt_0}(\Aut^G_{C^{\nu}}E^{\nu})|g \in C(\gamma)\}.$$

\end{prop}

\begin{proof}
Commutativity of the diagram is straightforward to check and
injectivity of $N$ follows from the fact that $\Psi_{\pt_0}$ is
injective.  The fact that the image commutes with $\gamma$ follows
from an argument similar to that for holonomy in
Proposition~\ref{prop:commute}.
\end{proof}

\subsection{Distinguished components of $\MM_{g,n}$}

\ 

Several distinguished components of $\MM_{g,n}$ will be useful for
our construction of $G$-CohFTs.  We describe them and their basic
properties in this subsection.

\subsubsection{The substack $\xi(\bm)$ of
$\MM_{0,3}(\bm)$}\label{substack}

\ 

\begin{df}\label{df:xi}
For any $\bm := (m_1, m_2, m_3) \in G_A^3$, if the product
$\prod_{i=1}^3 m_i$ is not $1$, then we define the stack $\xi(\bm)$ to
be the empty stack.  Otherwise, let $C$ denote the sphere $\CP^1$ with
special points $p_0:=0$ and $p_j = \exp(2\pi j \sqrt{-1}/3) $ for
$j\in \{1,2,3\}$.  For each $j\in \{1,2,3\}$ let $d_i$ be the path in
$C$ determined by following a straight line from $p_0$ to $p_i$.
These paths induce elements $\sigma_i:=\sigma_{d_i} \in
\pi_1(\CP^1-\{p_1,p_2,p_3\}, p_0)$ (see
Definition~\ref{df:path-loop}), which generate the fundamental group
and have trivial product: $\sigma_1 \sigma_2 \sigma_3 =1$.  Thus the
triple $\bd = (d_1, d_2,d_3)$ is in $P_C$.

We define $\xi(\bm)$
to be the connected component of $\MM_{0,3}(\bm)$ containing the
geometric point $\zeta(\bd;\bm)$, as defined in Definition~\ref{df:zeta}.
\end{df}

\begin{rem}\label{unique}
For any $m\in G$, it is clear that the component $\xi(m,m^{-1},1)$ is
the unique component of $\MM_{0,3}(m,m^{-1},1)$ such that all the
points $\pt_i$ lie in the same connected component of the admissible
$G$-cover $E$.
\end{rem}

\begin{lm}\label{lemma}
 For any $\bm \in G_A^3$ we have the following identities for the $\xi(\bm)$.
\begin{enumerate}
\item \label{one} $\rho(\gamma,\gamma,\gamma)(\xi(\bm))= \xi(\gamma
  m_1 \gamma^{-1},\gamma m_2 \gamma^{-1}, \gamma m_3 \gamma^{-1})$ for
  any $\gamma\in G$.

\item $\rho(m_1,1,1)(\xi(\bm))=\rho(1,m_2,1)(\xi(\bm)) =
\rho(1,1,m_3)(\xi(\bm))=\xi(\bm)$.

\item \label{it:brxi}
For the generators $b_1, b_2$ of the braid group $B_3$
\begin{align*}
b_1\xi(\bm) & = \xi(m_1 m_2 m_1^{-1}, m_1, m_3)\\ b_2\xi(\bm) & =
\xi(m_1, m_2 m_3 m_2^{-1}, m_2). \end{align*} Thus for any element
 $b\in B_3$, we have $$b\xi(\bm) = \xi(b\bm),$$ where $b$
acts on the triple $\bm$ via the Hurwitz action (i.e.,
the obvious action where, for example, $b_1(m_1,m_2,m_3) := (m_1 m_2 m_1^{-1}, m_1, m_3)$).

\item\label{sym} Let $s$ be an isomorphism induced from a cyclic
  permutation (also denoted $s$) in $S_3$, then $$s\xi(\bm) =
  \xi(s\bm).$$

\end{enumerate}
\end{lm}

\begin{proof}
The first identity follows from the fact that the global right
action translates all points in the admissible $G$-cover in $\xi$ by
$\gamma$.  Under this action, the $i$-th monodromy $m_i$ changes to
$\gamma^{-1}m_i \gamma$ for all $i=1,\ldots,n$.

The second statement follows from the fact that the action of
$\rho_i$ on the $i$-th point $\pt_i$ is the same (via the map
$\bar{\nu}_{\pt_0}$) as right multiplication acting on the right
$G$-coset $\langle m_i \rangle$; that is, the action
$\varphi_i(m_i)$ is trivial.

The third statement follows from studying the results of  sliding
points $p_j$ around $p_i$, which we now describe in the case of
$b_2$. The case of $b_1$ is essentially the same.

The transformation $T:z\mapsto 1/z$ takes $p_0 = 0$ to $\infty$,
fixes $p_1$, and interchanges $p_2$ and $p_3$.  Let $E'=T_*E :=
(T^{-1})^*E$, $\pt'_2:=T_*(\pt_3)$, $\pt'_3:=T_*(\pt_2)$, and
$\pt'_1:=T_*(\pt_1)$.  The pointed cover $(E', \pt'_1,\pt'_2,
\pt'_3)$ corresponds to the geometric point representing the image
of $\xi(\bm)$ under the action of the transposition $s_{(2,3)}$.
Let $\gamma$ be a straight path from $p_0$ to $\infty$ that passes
between $p_3$ and $p_1$, e.g., the path $\gamma(t) = -i/(1-t)$.
Note that via $\gamma$ we have an isomorphism of (un-pointed)
admissible $G$-covers $E' \cong E_{\chi,p_0}$, where $\chi$ is the
homomorphism $\pi_1(C_{gen},p_0) \rTo G$, given by taking $\gamma
T_*\sigma_i \gamma^{-1}$ to $m_i$ and with the induced
$\pt_{0,\chi}$ being the ``parallel transport" of $T^*(\pt_0)$
along $\gamma$.  The loop $\gamma T_*\sigma_3 \gamma^{-1}$ (around
$T(p_3) = p_2)$ and the loop $\gamma T_*\sigma_1 \gamma^{-1}$
(around $T(p_1)=p_1$) are homotopic to the loops $\sigma_2$ and
$\sigma_1$, respectively.  But the loop $\gamma T_*\sigma_2
\gamma^{-1}$ is homotopic to $\sigma_2 \sigma_3 \sigma_2^{-1}$.
Thus the (un-pointed) $G$-cover $E'$ is isomorphic to
$E'':=E_{\chi,p_0}$, where $\chi$ is the homomorphism taking
$\sigma_1$ to $m_1$, $\sigma_2$ to $m_2 m_3 m_2^{-1}$, and
$\sigma_3$ to $m_2$, that is, to the $G$-cover $E''$ associated to
$\xi(b_2\bm)$. And the points $\pt'_1$ and $\pt'_3$ are the same
as those that are induced on $\xi(b_2\bm)$.  However, the point
$\pt'_2$ is not the same as the point $\pt''_2$ induced on
$\xi(b_2\bm)$; indeed, $\pt''_2$ is induced by parallel transport
from $\pt_0$ along the path $d_2$, whereas $\pt'_2$ is induced by
parallel transport from $\pt_0$ along $\gamma T_*d_3 \gamma^{-1} =
(\sigma_3)^{-1} d_2$.  Thus they differ by the holonomy $m_2$ of
the loop $\sigma_3$; i.e., the claim of the third statement holds.

Finally, the last statement of the lemma follows from the fact that a
rotation (multiplication by $\pm\exp(2\pi i /3)$) of $\CP^1$ will
induce the permutation $s$ on $\xi$.
\end{proof}

The construction of $\xi(\bm)$ depends \emph{a priori} on the choices of
$p_i$ and $d_i$, but we will see in Proposition~\ref{prp:xiindep}  that  it
is independent of these choices.

Before we give that proposition, we need to understand better how the different
 braid actions interact. The fact that the braid action on paths agrees with the usual
  braid action on loops shows that for any $b\in B_3$  we may write $$b\bd =
  (\omega_1 d_{\psi(b(1))},\omega_2 d_{\psi(b(2))} , \omega_3 d_{\psi(b(3))})$$ for some choice of
  $\omega_i \in \pi_1(C_{gen},p_0)$ (and $\psi(b(i))$ is the 
action on $i$ of the permutation
  induced by the standard surjection $B_3 \rTo^{\psi} S_3$).    Let $b\chi$ denote the
  induced holonomy $\sigma_{b d_i}\mapsto m_i$, and let $\gamma_i$ be the
  image of $\omega_i$ in $G$ via $b\chi$.     It is clear that $b  \chi$  is the same
  homomorphism as that induced by taking $d_i \mapsto b m_i$, but the point
  $\pt(bd_i)$ differs from that defined by $\pt(d_{\psi(b(i))})$ by
  $\rho(\gamma_i)$.  That is, we have
\begin{align}
\beta(b) \xi(\bm) &= \beta(b)\zeta(\bd;\bm)\nonumber\\
            & = \zeta(b\bd;\bm)\nonumber\\
            &= \rho(\bgamma) \zeta(\bd;b\bm)\nonumber\\
            &= \rho(\bgamma) \xi(b \bm)\nonumber\\
            &= \rho(\bgamma) b \xi(\bm) .\label{eq:braidcompat}
    \end{align}

\begin{lm}\label{lm:betatriv}
The braid action $\beta$ on $\zeta_C \subset \MM_{0,3}$ factors
through the standard symmetric group action on $\MM_{0,3}$ via the
usual surjection $\psi:B_3 \rTo S_3$ to the symmetric group.  That
is, for any $b\in B_3$, $\bd \in P_C$, and $\bm$, such that $\prod
m_i =1$, we have $$\beta(b)\zeta(\bd;\bm)= \psi(b) \zeta(\bd;\bm).
$$
\end{lm}
\begin{proof}
By transitivity of the $B_n$ action on $P_C$,  for every
$\zeta(\bd';\bm)$ there exists a $b'\in B_n$ such that $\bd' = b'
\bd$, where $\bd$ is the set of paths used to define $\xi$.  So it
suffices to check this only in the case of $\xi(\bm)$; i.e., where
the paths are the standard $\bd$.   Checking the generators of
$B_3$ is now quite easy.  For example, in the case of  $b=b_1$ the
shift  $\bgamma$ is simply $(m_1,1, 1)$  and so
equation~(\ref{eq:braidcompat}) and Lemma~\ref{lemma}
item~(\ref{it:brxi})  gives
\begin{align*}
\beta(b_1)\xi(\bm) &= \rho(m_1,1,1)\xi(b\bm) \\
&= \rho(m_1,1,1)b \xi(\bm)\\
&= s_{1,2} \xi,
\end{align*}
as desired.
\end{proof}

 \begin{prop}\label{prp:xiindep}
For any $\bm\in G^3$ with $\prod_{i=1}^3 m_i =1$, any choice of points $p'_0,
\dots, p'_3 \in \CP^1$, and any choice of paths $d'_i$ from $p'_0$ to $p'_i$ for each $i\in \{1,2,3\}$
with trivial product (i.e., $\bd'=(d'_1,d'_2,d'_3) \in P_C$),
the geometric point of  $\MM_{0,3}(\bm)$ defined by  $\zeta(\bd',\bm)$ lies
in the component $\xi(\bm)$.
\end{prop}

\begin{proof}
Using the action of $\operatorname{PGL}(2,\nc)$ we may assume that
$p'_1=p_1, p'_2=p_2,$ and $p'_3=p_3$.

Moreover, given any path $\delta$ from $p_0$ to $p'_0$, we may replace
the paths $d'_i$ by $d'_i \delta$.  This gives an isomorphism between
the $n$-pointed admissible $G$ cover defined by the $d'_i$ and that
defined by the $d'_i \delta$.  Thus we may assume that $p'_0=p_0$.

Since both sets of paths $\bd =(d_1,d_2,d_3)$ (from the definition
of $\xi(\bm)$) and $\bd'=(d'_1,d'_2,d'_3)$ lie in the set $P_C$,
and since the braid action on $P_C$ is transitive
(Definition-Proposition~\ref{df-prp:trnsBact}), there is an
element $b'\in B_3$ such that $$\zeta(\bd';\bm) = \zeta(b'\bd;\bm)
=\beta(b')\xi(\bm).$$ Moreover, since the endpoint of each $d'_i$
is $p_i$, we must have $$b' \in \ker(\psi:B_3\rTo S_3),$$ that is,
$b'$ lies in the \emph{pure braid group}.

The proposition now follows from Lemma~\ref{lm:betatriv}.

\end{proof}

\subsubsection{Distinguished components of $\MM_{0,4}$}

\ 

\begin{df}\label{def:xifour}
Let $\bm=(m_1,\dots,m_4)$ be chosen so that $\prod_{i=1}^4 m_i = 1$,
and let $$m_{+} :=(m_1m_2)^{-1}, \quad m_{-} =m_{+}^{-1} .$$ We let
$\xi_{0,4}(\bm)$ denote the component of $\MM_{0,4}(\bm)$ which
contains the image of $\xi(m_1,m_2,m_{+}) \times \xi(m_{-}, m_3,m_4)$
under the gluing map
$$\varrho:\MM_{0,3}(m_1,m_2,m_{+}) \times
\MM_{0,3}(m_{-},m_3,m_4) \rTo  \linebreak[0] \MM_{0,4}(\bm).$$
\end{df}

\begin{df}\label{closed}
For any closed substack $Q \subseteq \MM_{g,n}$, consider the
homology class $[Q]$ in $H_\bullet(\MM_{g,n})$.
We define
\[ \db{Q} := 0 \]
when $Q$ is empty, otherwise,
\[ \db{Q} := \frac{1}{\deg(\st_Q)} [Q], \]
where $\deg(\st_Q)$ is the degree of the forgetful morphism
restricted to $Q$: $$\st_Q : Q \rTo \M_{g,n}.$$
\end{df}

\begin{lm}\label{ass}
Using the notation of  Definition~\ref{def:xifour},  let $$m'_{+}
:=(m_4m_1)^{-1}, \quad m'_{-}: =(m'_{+})^{-1} .$$ We further let
$\varrho'$ denote the gluing map composed with the cyclic
permutation $s=(4,3,2,1) \in S_4$, that is,  $\varrho'=s \circ
\varrho_{tree}$: $$\MM_{0,3} (m_4, m_1, m'_{+}) \times \MM_{0,3}(
m'_{-}, m_2, m_3 ) \rTo^{\varrho_{tree}} \linebreak[0]
\MM_{0,4}(m_4,  m_1,  m_2, m_3) \linebreak[0] \rTo^{s}
\MM_{0,4}(\bm),$$
 and we let $\varrho''$ denote the gluing map
$$\varrho'':\MM_{0,3} (m_1,m'_{+}, m_4) \times \MM_{0,3}(m_2,
m'_{-}, m_3) \rTo \linebreak[0] \MM_{0,4}(\bm).$$
\begin{enumerate}
\item The component $\xi_{0,4}(\bm)$ contains the image of $\xi(m_4,
  m_1, m'_{+}) \times \xi( m'_{-}, m_2, m_3 )$ under the map
  $\varrho'$ and the image of $\xi(m_1,m'_{+}, m_4) \times \xi( m_2,
  m'_{-}, m_3)$ under the map $\varrho''$.

\item We have the following equalities in $H_2(\xi_{0,4}(\bm))$:
\[
\db{\varrho(\xi(m_1,m_2,m_+)\times\xi(m_-,m_3,m_4))} =
\db{\varrho'(\xi(m_4,m_1,m'_+)\times\xi(m'_-,m_2,m_3))},
\]
and
\item
\[
\varrho_*(\db{\xi(m_1,m_2,m_+)}\otimes\db{\xi(m_-,m_3,m_4)}) =
\linebreak[0]
\varrho'_*(\db{\xi(m_4,m_1,m'_+)}\otimes\db{\xi(m'_-,m_2,m_3)}).
\]
\end{enumerate}
\end{lm}

\begin{proof}
For any choice $\bm \in G_A^4$ with $\prod_{i=1}^4 m_i = 1$, a
construction similar to that of $\xi(m_1,m_2,m_3)$ on $\CP^1 -
\{p_1,p_2,p_3,p_4\}$, with, say, $p_i := (\sqrt{-1})^i$ and $p_0:=0$,
and with straight-line paths $d_i$ to each $p_i$, gives a pointed
admissible $G$-cover of $\CP^1 - \{p_1,p_2,p_3,p_4\}$, which has two
obvious degenerations. The first degeneration is given by contracting
the great circle defined by $\{z=t(1+\sqrt{-1})|t\in \nr \cup
\infty \}$.  This can easily be seen to be the image of
$\xi(m_1,m_2,m_+) \times \xi(m_-, m_3,m_4)$ under the gluing map
$\varrho_{tree}:\MM_{0,3}(m_1,m_2,m_+) \times \MM_{0,3}(m_-,m_3,m_4)
\rTo \linebreak[0] \MM_{0,4}(\bm)$.  Similarly, the second
degeneration, given by contracting the great circle
$\{z=t(1-\sqrt{-1})|t\in \nr \cup \infty \},$ is the image of
$\xi(m_1,m'_+, m_4) \times \xi(m_2,m'_-, m_3)$ under the gluing map
$\MM_{0,3}(m_1, m'_+, m_4) \times \MM_{0,3}( m_2, m'_-, m_3) \rTo
\linebreak[0] \MM_{0,4}(\bm)$.  The
first claim follows from Lemma ~\ref{lemma} item ~\ref{sym} and the fact that
all these gluing morphisms are well-behaved under cyclic permutations.

To see the second claim, consider the forgetful morphism
\[
\st:\xi_{0,4}(\bm)\rTo\M_{0,4}.
\]
By pulling back the corresponding boundary divisors on $\M_{0,4}$, one
obtains the equality
\[
[\varrho(\xi(m_1,m_2,m_+)\times\xi(m_-,m_3,m_4))] \frac{A}{B}
=
[\varrho'(\xi(m_4,m_1,m'_+)\times\xi(m'_-,m_2,m_3))] \frac{A'}{B},
\]
where $A$ is the order of the automorphism group of
$\varrho(\xi(m_1,m_2,m_+)\times\xi(m_-,m_3,m_4))$, $A'$ is the
order of the automorphism group of
$\varrho'(\xi(m_4,m_1,m'_+)\times\xi(m'_-,m_2,m_3))$, and $B$ is
the order of the automorphism group of a generic point in
$\MM_{0,4}(\bm)$.

Finally, we observe that
\[
\varrho_*([\xi(m_1,m_2,m_+)]\otimes[\xi(m_-,m_3,m_4)]) =
[\varrho(\xi(m_1,m_2,m_+)\times\xi(m_-,m_3,m_4))] \frac{C}{D_+
D_-},
\]
where $C$ is the order of the automorphism group of a generic
point in $\varrho(\xi(m_1,m_2,m_+)\times\xi(m_-,m_3,m_4))$, $D_+$
is the order of the automorphism group of $\xi(m_1,m_2,m_+)$, and
$D_-$ is the order of the automorphism group of
$\xi(m_-,m_3,m_4)$. This equation, together with its counterpart
from $\varrho_*([\xi(m_1,m_2,m_+)]\otimes[\xi(m_-,m_3,m_4)])$ and
the previously derived equation, yields the desired result.
\end{proof}

\subsubsection{Distinguished components of $\MM_{1,1}$}

\ 

\begin{df}
Choose elements $a,b,m_1 \in G$ such that $m_1= [a,b]$. Let
$\varrho_b$ be the composition of the morphisms
\begin{equation}
\xi(m_1,b, a b^{-1} a^{-1})\rTo^{\rho_3(a)} \MM_{0,3}(m_1,b, b^{-1} )
\rTo^{\varrho'_{b}} \MM_{1,1}(m_1),
\end{equation}
where the first morphism is right action by $a$ in the third factor,
and the second morphism is the gluing morphism identifying the 2nd and
3rd marked points.

Similarly, let $\varrho_a$ be the composition of the morphisms
\begin{equation}
\xi(m_1,b a b^{-1}, a^{-1}) \rTo^{\rho_2(b)}\MM_{0,3}(m_1,a,a^{-1}) 
\rTo^{\varrho'_a}\MM_{1,1}(m_1),
\end{equation}
where the first morphism is right action by $b$ in the second factor,
and the second is again the gluing morphism identifying the 2nd and
3rd marked points.

We define $\xi_{1,1}(m_1,a,b)$ to be the component of $\MM_{1,1}(m_1)$
containing the image of $\varrho_b$.
\end{df}

\begin{lm}\label{lm:trace}

The images of $\varrho_a$ and $\varrho_b$  lie in
the same connected component $\xi_{1,1}(m_1,a,b)$ of  $\MM_{1,1}(m_1)$.
Moreover,  the following equation holds in $H_2(\MM_{1,1}(m_1))$:
\begin{equation}\label{eq:trace}
\varrho'_{b*}(\db{\rho_3(a)\xi(m_1,b,a b^{-1} a^{-1})}) =
\varrho'_{a*}(\db{\rho_2(b)\xi(m_1,b a b^{-1}, a^{-1})}).
\end{equation}
\end{lm}
\begin{proof}
The images of $\varrho_a$ and $\varrho_b$ are degenerations of the
same smooth admissible $G$-cover over a smooth torus.  In particular,
consider a smooth, one-pointed torus $(T, p_1)$ with generators
$\alpha$, $\beta$, and $\gamma$ of $\pi_1(T,p_0)$ for some point
$p_0$, with $\gamma$ corresponding to the loop $\sigma_{d}$ induced by
a path $d$ from $p_0$ to $p_1$ (as in Definition~\ref{df:path-loop}),
and $[\alpha,\beta]=\gamma$.  The homomorphism $\chi:\pi_1(T,p_0)\rTo
G$ that takes $\alpha$, $\beta$ , and $\gamma$ to $a$, $b$, and $m_1$,
respectively, defines an admissible $G$-cover $E_{\chi}$, and a point
$\pt_{0,\chi}$.  Parallel transport along $d$ induces a point $\pt_1$
with monodromy $m_1$, giving us a pointed admissible $G$-cover
$E_{\chi,\pt_1}$.

It is straightforward to see that the image of $\varrho_a$
corresponds to the $\alpha$-cycle shrinking to become a node,
while the image of $\varrho_b$ corresponds to the $\beta$-cycle
shrinking to become a node.  Thus both images lie in the same
connected component $\xi_{1,1}(m_1)$ of $\MM_{1,1}(m_1)$.

Equation~(\ref{eq:trace})  follows from the identity
\[
\varrho'_{b*}([\rho_3(a)\xi(m_1,b,a b^{-1} a^{-1})]) \frac{A}{B}
=
\varrho'_{a*}([\rho_2(b)\xi(m_1,b a b^{-1}, a^{-1})])
\frac{A'}{B'},
\]
where $A$ is the order of the automorphism group of
$\rho_3(a)\xi(m_1,b,a b^{-1} a^{-1})$, $A'$ is the order of the
automorphism group of $[\rho_2(b)\xi(m_1,b a b^{-1}, a^{-1})]$,
$B$ is the order of the automorphism group of
$\varrho'_b(\rho_3(a)\xi(m_1,b,a b^{-1} a^{-1}))$, and $B'$ is the
order of the automorphism group of $\varrho'_a(\rho_2(b)\xi(m_1,b
a b^{-1}, a^{-1}))$.

However, $B = B'$, as their corresponding automorphism groups are both isomorphic to $C(a,b)\subseteq G$ (see Proposition~\ref{prp:aut-slf-sew}).
\end{proof}

\section{The category of $G$-graded $G$-modules}\label{3}

\ 

In this section, we briefly review some well-known facts from the category of
$G$-graded $G$-modules (see \cite{Kas,BaKi}) which will be useful in the
sequel.

\subsection{$G$-graded $G$-modules and their $G$-coinvariants}

\ 

\begin{df} \label{df:module} Let $\ch := \bigoplus_{m\in G} \ch_{m}$ be a
finite-dimensional $G_A$-graded vector space which is endowed with the
structure of a right $G$-module  $\rho(\gamma):\ch\irightarrow \ch$ for all
$\gamma$ in $G$, with $\rho(\gamma)$ taking $\ch_{m}$ to $\ch_{\gamma^{-1}
m\gamma}$ for all $m$ in $G$. \emph{$(\ch,\rho)$ is said to be a
$G$-graded $G$-module.}

A \emph{$G$-invariant metric $\eta$ on a $G$-graded $G$-module $\ch$} is a
symmetric, nondegenerate, bilinear form $\eta$ on $\ch$ which is
$G$-invariant (under the diagonal $G$ action) and which respects the
grading, i.e., for all $v_{m_+}$ in $\ch_{m_+}$ and $v_{m_-}$ in $\ch_{m_-}$ we have $\eta(v_{m_+},v_{m_-}) =
0$ unless $m_+ m_- = 1$.
\end{df}

$G$-graded $G$-modules form a category whose objects are $G$-graded
$G$-modules and whose morphisms are homomorphisms of $G$-modules which respect
the $G$-grading. Furthermore, the dual of a $G$-graded $G$-module inherits
the structure of a $G$-graded $G$-module.

\begin{ex}
Any finite-dimensional $G$-module $V$ is a $G$-graded $G$-module
where $\ch_1 := V$ and $\ch_m := 0$ for all $m$ not equal to 1 in
$G$.
\end{ex}

\begin{ex}
The simplest example of a nontrivial $G$-graded $G$-module is $\nc[G]$, the
free vector space generated by $G$, with its natural $G$-grading, endowed
with the $G$-action  $\rho(\gamma)m := \gamma^{-1} m \gamma$ for all $\gamma,
m$ in $G$.
\end{ex}

\begin{df}
Recall that  $\Gb$ is the set of conjugacy classes of $G$,
the conjugacy class of $m$ in $G$ is denoted by $\mb$, and the conjugacy
class of $m^{-1}$ is denoted by $\mb^{-1}$.

A section $s$ of the natural map $G\rTo\Gb$ is said to be
\emph{involutive} if $s(\mb^{-1}) = s(\mb)^{-1}$ for all $\mb$.
\end{df}

\begin{df} Let $(\ch,\rho)$ be a $G$-graded $G$-module. Let
$\pi_G : \ch\to\ch$ be the averaging map
\[
\pi_G(v) :=
\frac{1}{|G|}\sum_{\gamma\in G} \rho(\gamma) v
\]
for all $v$ in $\ch$. Let $\chb$ be the image of $\pi_G$. The vector space $\chb$ is called the \emph{space of
$G$-coinvariants of $\ch$,} and it inherits a grading by $\Gb$,
denoted by
\[
\chb = \bigoplus_{\gb\in\Gb} \chb_{\gb}.
\]
If $\eta$ is a metric on $\ch$, then let $\etab$ be the
restriction of the metric $\frac{1}{|G|} \eta$ to $\chb$.
\end{df}

\begin{rem}
The reason for the factor of $\frac{1}{|G|}$ in the definition of $\etab$ will
become evident when we discuss the geometry of \Gcfts.
\end{rem}

Let us describe $\chb$ in terms of $\ch$.

\begin{prop}\label{prop:gmodule}
Let $(\ch,\rho)$ be a $G$-graded $G$-module with a $G$-invariant metric $\eta$.
\begin{enumerate}
\item \label{prop:gmodule_a}  Consider $v_{\mb}$ in $\chb_\mb$, where $v_\mb =
\sum_{m'\in\mb} v_{m'}$.  For all $m'$ in $\mb$, $v_{m'}$ belongs
to $\ch_{m'}^{C(m')}$, the $C(m')$-invariant subspace of
$\ch_{m'}$. In particular, for all $v_m$ in $\ch_m$,
\[
\pi_G(v_m) = \pi_G(\pi_{C(m)}(v_m)),
\]
where $\pi_{C(m)} : \ch_m\rTo\ch_m^{C(m)}$ is the averaging map
\[
\pi_{C(m)}(v_m) := \frac{1}{|C(m)|}\sum_{\gamma\in C(m)} \rho(\gamma) v_m.
\]
\item \label{prop:gmodule_b} For all $m$ in $G$, the map $\pi_m :
\ch_m^{C(m)}\rTo\chb_\mb$, defined as
\[
\pi_m(v_m) := \pi_G(v_m),
\]
is an isomorphism of vector spaces.
\item \label{prop:gmodule_c} For all $m_\pm$ in $G$ and $v_{m_\pm}$ in
$\ch_{m_\pm}^{C(m_\pm)}$, where $m_+ m_- = 1$, we have
\[
\eta(\pi_{m_+}(v_{m_+}),\pi_{m_-}(v_{m_-})) =
\eta(\pi_G(v_{m_+}),\pi_G(v_{m_-})) = \frac{|C(m_+)|}{|G|}
\eta(v_{m_+},v_{m_-}).
\]
\item  \label{prop:gmodule_d}If $s$ is an involutive section of the natural map $G\rTo
\Gb$,
then
\[
\bigoplus_{\mb\in\Gb} \ch_{s(\mb)}^{C(s(\mb))}\rTo\chb,
\]
taking $v_{s(\mb)}\mapsto \pi_G(v_{s(\mb)})$, is an isomorphism of
vector spaces which is not an isometry.
\item  \label{prop:gmodule_e} $\etab$ is nondegenerate, i.e., $\chb $
is a $\Gb$-graded vector space with metric $\etab$.
\end{enumerate}
\end{prop}
\begin{proof}
To prove part~(\ref{prop:gmodule_a}), consider $w_m$ in $\ch_m$.
We have
\begin{eqnarray*}
\pi_G(w_m) &=& \frac{1}{|G|}\sum_{\gamma'\in G}\rho(\gamma')w_m \\
           &=& \frac{1}{|G|}\sum_{[\gamma]\in (C(m)\backslash G)} \sum_{c\in
C(m)}
\rho(c\gamma) w_m \\
       &=& \frac{|C(m)|}{|G|}\sum_{[\gamma]\in C(m)\backslash G}
\rho(\gamma) \frac{1}{|C(m)|}\sum_{c\in C(m)} \rho(c)w_m \\
       &=& \frac{|C(m)|}{|G|}\sum_{[\gamma]\in C(m)\backslash G}
\rho(\gamma) \pi_{C(m)}(w_m).
\end{eqnarray*}
We conclude that
\begin{equation}
\pi_G(w_m) = \frac{|C(m)|}{|G|}\sum_{[\gamma]\in C(m)\backslash G}
\pi_{C(\gamma^{-1}m\gamma)}(\rho(\gamma) w_m),
\end{equation}
which finishes the proof.

We prove part~(\ref{prop:gmodule_b}) by showing that the map $f_m
: \chb_{\mb}\rTo\ch_m^{C(m)}$, defined by
\begin{equation}
f_m(\sum_{m'\in\mb} v_{m'}) := \frac{|G|}{|C(m)|} v_m,
\end{equation}
is the inverse of $\pi_m$. Notice that the right hand side is
$C(m)$-invariant by part~(\ref{prop:gmodule_a}). Consider $w_m$ in
$\ch_m^{C(m)}$.  We have
\[
f_m(\pi_G(w_m)) = f_m\left(\frac{|C(m)|}{|G|} \sum_{[\gamma]\in
C(m)\backslash G} \rho(\gamma) w_m\right)\]
\[ =\frac{|C(m)|}{|G|} f_m(w_m) = \frac{|C(m)|}{|G|}
\frac{|G|}{|C(m)|} w_m = w_m.
\]
Therefore, $\pi_m$ is an isomorphism.

To prove part~(\ref{prop:gmodule_c}), observe that
\begin{eqnarray*}
\eta(\pi_G(v_{m_+}),\pi_G(v_{m_-}))
&=&\frac{1}{|G|^2}\sum_{\gamma_\pm\in G} \eta(\rho(\gamma_+)
v_{m_+},\rho(\gamma_-) v_{m_-}) \\ &=&
\frac{1}{|G|^2}\sum_{\gamma_\pm\in G} \eta(\rho(\gamma_-^{-1})
\rho(\gamma_+) v_{m_+},v_{m_-}) \\
&=&\frac{1}{|G|^2}\sum_{\gamma_\pm\in G}
\eta(\rho(\gamma_+\gamma_-^{-1}) v_{m_+},v_{m_-})\\
&=&\frac{1}{|G|^2}\sum_{\gamma\in G}\sum_{\gamma_+\in G}
\eta(\rho(\gamma) v_{m_+},v_{m_-})\\
 &=&\frac{1}{|G|}\sum_{\gamma\in
G}\eta(\rho(\gamma) v_{m_+},v_{m_-})\\
&=&\frac{1}{|G|}\sum_{\gamma\in C(m_+)}\eta(\rho(\gamma)
v_{m_+},v_{m_-})\\ &=&\frac{|C(m_+)|}{|G|}\eta(v_{m_+},v_{m_-}),
\end{eqnarray*}
where we used the $C(m_+)$-invariance of $v_{m_+}$ in the last equality.

Part~(\ref{prop:gmodule_d}) follows immediately from
(\ref{prop:gmodule_b}) and (\ref{prop:gmodule_c}). Involutivity of
$s$ is needed to insure that $\bigoplus_{\mb\in\Gb}
\ch_{s(\mb)}^{C(s(\mb))}$ inherits a metric from $\ch$ compatible
with its grading.

To prove part~(\ref{prop:gmodule_e}), let $m_+ m_- = 1$, so that
$\eta$ restricted to $\ch_{m_+}\oplus\ch_{m_-}$ is nondegenerate.
$\ch_{m_\pm}$ is a $C(m_\pm)$-module, so one can write
\begin{equation}\label{eq:representations}
\ch_{m_\pm} = \ch_{m_\pm}^{C(m_\pm)} \oplus \ch_{m_\pm}'
\end{equation}
as $C(m_\pm)$-modules, where $\ch_{m_\pm}'$ is the direct sum of
all nontrivial irreducible representations of $C(m_\pm)$ appearing
in $\ch_{m_\pm}$. Notice that since $m_+ m_- = 1$, we have $C(m_+)
= C(m_-)$.

Since $\eta$ is $C(m_+)$-invariant, $\eta$ restricted to
$\ch_{m_+}\oplus \ch_{m_-}$ is the direct sum of $\eta$ restricted
to $\ch_{m_+}^{C(m_+)}\oplus \ch_{m_-}^{C(m_+)}$ and  $\eta$
restricted to $\ch_{m_+}'\oplus \ch_{m_-}'$. Therefore, $\eta$
restricted to $\ch_{m_+}^{C(m_+)}\oplus \ch_{m_-}^{C(m_-)}$ is
nondegenerate.

Let $v_{m_+}$ be in $\ch_{m_+}^{C(m_+)}$. Suppose that
$\eta(\pi_G(v_{m_+}),\pi_G(v_{m_-})) = 0$ for all $v_{m_-}$ in
$\ch_{m_-}^{C(m_-)}$.  By part~(\ref{prop:gmodule_c}), this is
equivalent to the condition $\eta(v_{m_+},v_{m_-}) = 0$ for all
$v_{m_-}$ in $\ch_{m_-}^{C(m_-)}$. However, $\eta$ restricted to
$\ch_{m_+}^{C(m_+)}\oplus \ch_{m_-}^{C(m_-)}$ is nondegenerate,
therefore, $v_{m_+} = 0$. Thus $\eta$ restricted to $\chb$
is also non-degenerate.
\end{proof}

\subsection{Tensor products and the braid group}

\ 

 As is usual in the representation theory of groups, there are two
kinds of tensor products associated to $G$-graded $G$-modules,

\begin{df}
Let $\ch'$ be a $G'$-graded $G'$-module and $\ch''$ be a $G''$-graded
$G''$-module. Their vector space tensor product $\ch'\otimes \ch''$ is
naturally a $G'\times G''$-graded $G'\times G''$-module called the
\emph{external tensor product of $\ch'$ and $\ch''$}.
\end{df}

On the other hand, the category of $G$-graded $G$-modules has a natural
tensor product which differs from the tensor product of their underlying
vector spaces.

\begin{df} Let $\ch' := \bigoplus_{m\in G} \ch'_m$ and
$\ch'' := \bigoplus_{m\in G} \ch''_m$ be two $G$-graded
$G$-modules. Let
\[
\ch'\btimes\ch'' := \bigoplus_{m\in G} \ch'_m\otimes\ch''_m,
\]
with the induced $G$-module structure, where $G$ acts diagonally.
We call $\ch'\btimes\ch''$ the \emph{tensor product of $\ch'$ and
$\ch''$}.
\end{df}

\begin{rem}
The $G$-graded $G$-module $\nc[G]$ has the important property that
\[
\ch\btimes\nc[G]\cong\nc[G]\btimes\ch\cong\ch
\]
for any $G$-graded $G$-module $\ch$.
\end{rem}

Finally, we note that objects in this category have a natural action of the
braid group, which we now describe.

\begin{df}
Let $\ch$ be a $G$-graded $G$-module. Its $n$-fold tensor product
$\ch^{\otimes n}$ inherits the structure of a right
$G^n \rtimes S_n$-module where the symmetric group $S_n$ acts on
$\ch^{\otimes n}$ by permuting its factors.

For all $i=1,\ldots,n-1$, let $b_i : \ch^{\otimes n} \rTo \ch^{\otimes n}$ be
defined by
\[
b_i(v_{m_1}\otimes\cdots \otimes v_{m_i}\otimes
v_{m_{i+1}}\cdots\otimes v_{m_n}) :=
v_{m_1}\otimes\cdots \otimes (\rho(m_i^{-1})v_{m_{i+1}})\otimes
v_{m_i}\otimes\cdots\otimes v_{m_n}
\]
for all $v_{m_j}$ in $\ch_{m_j}$, $m_j$ in $G$, and $j=1,\ldots,n-1$.
\end{df}

The following proposition is immediate.

\begin{prop}
The elements $\{b_1,\ldots,b_n\}$ on $\ch^{\otimes n}$ give an
action of $B_n$ on $\ch^{\otimes n}$.
\end{prop}

\section{$G$-equivariant cohomological field theories}\label{4}

\ 

In this section, we introduce the notion of a \Gcft\ , defined in
terms of $\MM_{g,n}$, and prove some of its basic properties.

\subsection{\Gcfts\ and $G$-Frobenius algebras}\label{sec:gcohft}

\ 

\begin{df}\label{df:GCFT} A tuple
  $((\ch,\rho),\eta,\{\Lambda_{g,n}\},\vac)$ is said to be a
  \emph{$G$-equivariant Cohomological Field Theory (\Gcft)} if the
  following axioms hold:

\begin{enumerate}
\item ($G$-Graded $G$-module) $(\ch,\rho)$ is a $G$-graded $G$-module. The subspace
$\ch_1$ is called the \emph{untwisted sector of the \Gcft}, and
$\ch_m$, where $m\not=1$, is called a \emph{twisted sector of the
\Gcft.}
\item ($G^n\rtimes S_n$ Invariance) For all $\bm := (m_1,\ldots,m_n)$ in $G^n$
and all stable pairs $(g,n)$, if we denote $\ch_{\bm} := \bigotimes_{i=1}^n
\ch_{m_i}$, then 
$\Lambda_{g,n}$ is  an element of $\bigoplus_\bm
H^\bullet(\MM_{g,n}(\bm))\otimes \ch_{\bm}^*$ 
which is invariant under the diagonal action of $G^n\rtimes S_n$.
\item(Identity) The element $\vac$ in $\ch_1$ is non-zero, and is
  called the \emph{flat identity} or \emph{vacuum vector}.
\begin{enumerate}
\item ($G$-Invariance of the Identity) The vacuum vector $\vac$ is $G$
  invariant, i.e., $\rho(\gamma)\vac = \vac$ for all $\gamma$ in $G$.
\item (Flat Identity) Under the forgetting tails morphism
$\tst:\MM_{g,n+1}(\bm,1)\linebreak[1] \rTo \MM_{g,n}(\bm)$, we
have
\[
\Lambda_{g,n+1}(v_{m_1},\ldots,v_{m_n},\vac) =
\tst^*\Lambda_{g,n}(v_{m_1},\ldots,v_{m_n})
\]
for all $\bm$ in $G^n$, and $v_{m_i}$ in $\ch_{m_i}$ for all
$i=1,\ldots, n$.
\end{enumerate}
\item (Metric) $\eta$ is a symmetric, nondegenerate, bilinear form on $\ch$
such that
\[
\eta(v_{m_1},v_{m_2}) := \int_{\db{\xi(m_1,m_2,1)}}
\Lambda_{0,3}(v_{m_1},v_{m_2},\vac).
\]
 It follows that $\eta(v_{m_1},v_{m_2})=0$ unless $m_1 m_2 = 1$.
 Recall that $\xi$ is defined in Subsection ~\ref{substack} and
 the scaled class $\db{Q}$  in Definition
 ~\ref{closed}.
\item (Factorization) Fix any $m_+ \in G$ and $m_-:=(m_+)^{-1}$.  Let the set
$\{e_{\alpha}\}$ be a basis for $\ch_{m_+}$, the set
$\{\hat{e}_{\beta}\}$ be a basis for $\ch_{m_-}$, and
$\eta^{\alpha \beta}$ be the inverse of the metric
$$\eta:\ch_{m_+} \otimes \ch_{m_-} \rTo \nc$$ relative to these
bases.
\begin{enumerate}
\item\label{5a} For all stable pairs $(g_1,n_1+1)$ and $(g_2,n_2+1)$
  let $g=g_1+g_2$ and $n=n_1+n_2$.  For all $\bm$ in $G^n$ and all
  $(v_{m_1} \dots, v_{m_n}) \in \ch_\bm$ we require
\begin{equation*}
(\varrho^*_{tree}\Lambda_{g,n})(v_{m_1}, \dots,
v_{m_n}) = \sum_{\alpha,\beta}\Lambda_{g_1,n_1 +1}(v_{m_{i_1}},
\dots, v_{m_{i_{n_1}}},e_{\alpha})\eta^{\alpha
\beta}\Lambda_{g_2,n_2+1}(\hat{e}_{\beta},v_{m_{j_1}}, \dots,
v_{m_{j_{n_2}}})
\end{equation*}
for all partitions $\{ i_1,\ldots,i_{n_1} \}\sqcup
\{ j_1,\ldots,j_{n_2} \}$ of the set $\{ 1,\ldots,n\}$.
\item For all stable pairs $(g-1,n+2)$, all $\
\bm \in G^n$, and all $(v_{m_1},\dots, v_{m_n})  \in \ch_{\bm}$, the classes
$\Lambda$ must satisfy
\begin{equation*}
(\varrho^*_{loop}\Lambda_{g,n})(v_{m_1}, \dots,
v_{m_n}) = \sum_{\alpha,\beta} \Lambda_{g-1, n+2}(v_{m_1}, \dots,
v_{m_n}, e_{\alpha}, \hat{e}_{\beta})\eta^{\alpha \beta}.
\end{equation*}

\end{enumerate}
\end{enumerate}
\end{df}

\begin{rem}
If $G$ is the trivial group, then a \Gcft\ coincides with a \cft\ in the
sense of Kontsevich-Manin \cite{KoMa}.
\end{rem}

\begin{ex}\label{ex:grpring}
The simplest example of a \Gcft\ has as its state space
$\ch=\bigoplus_{m \in G}    \ch_{m} :=H^{\bullet}(G)=H^{0}(G) \cong
\nc [G]$ as $G$-graded $G$-modules, i.e., if $\{ e_m \}_{m\in G}$
denotes the obvious basis in $\ch$, then the $G$-action
$\rho(\gamma):\ch_m\to\ch_{\gamma^{-1}m\gamma}$ is
$\rho(\gamma)(e_m) := e_{\gamma^{-1} m \gamma}$ for all $\gamma,
m$ in $G$.

For all $\bm = (m_1,\ldots,m_n)$ in $G^n$, let
\[
\Lambda_{g,n}(e_{m_1},\ldots,e_{m_n}) := \mathbf{e}^*
\mathbf{1}_\bm,
\]
where $\mathbf{e}:\MM_{g,n}\rTo G^n$, and $\vac_\bm$ in $H^0(G^n)$
denotes the fundamental class of the point $\bm$ in $G^n$.

It follows that
\[
\eta(e_{m_1},e_{m_2}) := \int_{\db{\xi(m_1,m_2,1)}}
\Lambda_{0,3}(e_{m_1},e_{m_2},\vac) = \delta_{m_1,m_2^{-1}}.
\]
\end{ex}

\begin{df}
We will call the \Gcft\ of the previous example the \emph{group ring
\Gcft}, and we will denote it simply by $\nc[G]$ whenever it is clear from
context that we mean the group ring \Gcft\ and not just the ring
itself.
\end{df}

\begin{rem}
We will see in the next section that this \Gcft\ induces the
$G$-Frobenius algebra $\nc[G]$, and a standard argument (along the
lines of \cite{Tu}) shows that the two constructions are actually
equivalent, thus we are justified in the terminology and notation
of the previous definition.
\end{rem}

\subsection{Tensor products of  equivariant \cfts}

\ 

Given two  equivariant \cfts, one can construct a new one by
taking their tensor product. As in the case of $G$-graded
$G$-modules, there are two tensor products associated to \Gcfts.
The first, the external tensor product, associates to a \Gcft\ and
a \Gpcft\ a $G\times G'$-\cft. The second is a tensor product in
the category of \Gcfts.

\begin{prop}\label{prop:ExternalFP}
For all $\bm'$ in ${G'}^n$ and $\bm''$ in ${G''}^n$, let $\bm'\times\bm''$
denote the element $((m'_1,m''_2),\ldots,\linebreak[0](m'_n,m''_n))$ in $(G'\times G'')^n$.
Consider the commuting diagram
\[
\begin{diagram}
\M^{G'\times G''}_{g,n}(\bm'\times\bm'') &\rTo^{\Upsilon} &
\M^{G'}_{g,n}(\bm')\times_{\M_{g,n}} \M^{G''}_{g,n}(\bm'')
 & \rTo^{\pr''} & \MMpp_{g,n}(\bm'') \\
& &\dTo^{\pr'} & & \dTo^{\st''}
\\
& & \MMp_{g,n}(\bm')  & \rTo^{\st'} & \M_{g,n} \\
\end{diagram}
\]
where $\st'$ and $\st''$ forget the pointed admissible covers and
$\M^{G'}_{g,n}(\bm')\times_{\M_{g,n}} \M^{G''}_{g,n}(\bm'')$ is
the fibered product with projections $\pr'$ and $\pr''$.
The map $\Upsilon$ takes an object $(E\rTo C;
\pt_1,\ldots,\pt_n)$ to $((E'\rTo C;
\pt'_1,\ldots,\pt'_n),(E''\rTo C; \pt''_1,\ldots,\pt''_n))$, where
$E'$ is the variety $E/G''$ and $\pt'_i$ is the marked point on
$E'$ induced by $\pt_i$, 
$E''$ is the variety $E/G'$ and $\pt''_i$ is the marked
point on $E''$ induced by $\pt_i$.

\begin{enumerate}
\item The morphism $\Upsilon$ preserves the $(G'\times G'')^n$ and $S_n$
actions. \label{prop:ExternalFP_a}
\item The morphism $\pr'$ is $G'^n$-equivariant and $\pr''$ is
$G''^n$-equivariant.\label{prop:ExternalFP_b}
\item The morphisms $\pr', \pr'',\st',\st''$ are $S_n$-equivariant. \label{prop:ExternalFP_c}
\item The morphisms $\Upsilon, \pr', \pr'',\st',\st''$ commute with the gluing
morphisms.\label{prop:ExternalFP_d}
\end{enumerate}
\end{prop}

\begin{proof}
For Part~(\ref{prop:ExternalFP_a}) note that  $\Upsilon$ is a morphism because both $E'\rTo C$ and $E''\rTo C$ are admissible $G'$-, respectively $G''$-covers with the
proper monodromies. 
The equivariance under the actions of $(G'\times G'')^n$ and $S_n$ is manifest.

Similarly, Parts~(\ref{prop:ExternalFP_b}) and~(\ref{prop:ExternalFP_c}) are
manifest.

We now treat part~(\ref{prop:ExternalFP_d}) in the case of the loop for the
morphism $\pr'$. For all $m'_\pm$ in $G'$ and $m''_\pm$ in $G''$ such that
$m_+' m_-' = m_+'' m_-'' = 1$, consider the diagram
\begin{equation}
\begin{diagram}
\M^{G'\times G''}_{g-1,n+2}(\bm'\times \bm'',(m_+',m_+''),(m'_-,m_-''))
 & \rTo^{\varrhot} & \M^{G'\times G''}_{g,n}(\bm' \times \bm'') \\
\dTo^{\prc'} & & \dTo^{\pr'}
\\
\MMp_{g-1,n+2}(\bm',m_+',m'_-)  & \rTo^{\varrhot'} & \MMp_{g,n}(\bm') \\
\end{diagram}
\end{equation}
where $\varrhot$ and $\varrhot'$ are the gluing morphisms and $\prc'$ and
$\pr'$ are the canonical projections. Part~(\ref{prop:ExternalFP_d}) states
that this diagram commutes, which follows immediately from the definition of
the morphisms involved. Similarly, the analogous diagrams for $\Upsilon$,
$\pr''$, $\st'$, and $\st''$ also commute.  The proof in the case of the tree
is identical and will be omitted.
\end{proof}

\begin{crl}
Let $(\ch',\eta',\{\Lambda'_{g,n}\},\vac')$ be a \Gpcft\ and
$(\ch'',\eta'',\{\Lambda''_{g,n}\},\vac'')$ be a \Gppcft. If we define
\begin{equation}\label{eq:externaltensor}
\Lambda_{g,n}(v'_{m'_1}\otimes v''_{m''_1},\ldots, v'_{m'_n}\otimes
v''_{m''_n}) := \Upsilon^*(
(\pr'^*\Lambda'_{g,n}(v'_{m'_1},\ldots, v'_{m'_n}))\cup
(\pr''^*\Lambda''_{g,n}(v''_{m''_1},\ldots, v''_{m''_n})))
\end{equation}
for all $v'_{m'_i}$ in $\ch'_{m'_i}$ and $v''_{m''_i}$ in
$\ch''_{m''_i}$, where the morphisms $\pr'$ and $\pr''$ are
defined as in Proposition \ref{prop:ExternalFP}, then
$(\ch'\otimes\ch'',\eta'\otimes\eta'',\{\Lambda_{g,n}\},\vac'\otimes\vac'')$
is a $G'\times G''$-\cft.
\end{crl}
\begin{proof}

Let $G := G'\times G''$. Using the tensor product of a $G'$-graded
$G'$-module and $G''$-graded $G''$-module, $\ch'\otimes\ch''$
inherits the structure of a $G$-graded $G$-module. The
$G$-invariance of $\vac'\otimes \vac''$ follows.

The               $G^n$- and $S_n$-invariance follow from Proposition
\ref{prop:ExternalFP}(\ref{prop:ExternalFP_b}) and
(\ref{prop:ExternalFP_c}), respectively.

The flatness of the identity follows immediately from the definition of
$\Lambda_{g,n}$.

The metric axiom follows from observation that since $\M_{0,3}$ is a point,
the fibered product $\M^{G'}_{0,3}(\bm')\times_{\M_{0,3}}
\M^{G''}_{0,3}(\bm'')$ is equal to $\M^{G'}_{0,3}(\bm') \times
\M^{G''}_{0,3}(\bm'')$.

We prove the factorization axiom in the case of the loop---the case of
the tree is similar.  Let us adopt the notation from Proposition
\ref{prop:ExternalFP} and define $v'_{\bm'}\times v''_{\bm''}$ to be
$(v'_{m'_1}\otimes v''_{m''_1},\ldots, v'_{m'_n}\otimes v''_{m''_n})$
for all $v'_{\bm'}$ in $\ch'_{\bm'}$ and $v''_{\bm''}$ in
$\ch''_{\bm''}$.

From the definition of $\Lambda$ we have
\begin{equation*}
\varrhot^*\Lambda_{g,n}(v'_{\bm'}\times v''_{\bm''}) =
((\pr'\times\pr'')\circ\Deltat\circ\Upsilon\circ\varrhot)^*
(\Lambda'_{g,n}(v'_{\bm'}) \otimes \Lambda''_{g,n}(v''_{\bm''})),
\end{equation*}
where $\Deltat$ is the diagonal morphism 
$$\Deltat: \M^{G'}_{g,n}(\bm')\times_{\M_{g,n}}\M^{G''}_{g,n}(\bm'')
\rTo \M^{G'}_{g,n}(\bm')\times_{\M_{g,n}}\M^{G''}_{g,n}(\bm'') \times
\M^{G'}_{g,n}(\bm')\times_{\M_{g,n}}\M^{G''}_{g,n}(\bm'').$$ Let
$\Delta$ denote the diagonal morphism
$$\Delta:\MM_{g,n}(\bm\times\bm') \rTo 
\MM_{g,n}(\bm\times\bm')\times \MM_{g,n}(\bm\times\bm')$$ 
and $\Deltac$ denote the diagonal morphism 
associated to 
$\MM_{g-1,n+2}(\bm'\times\bm'',(m'_+,m''_+),(m'_-,m''_-))$
for   any 
$m'_\pm$ in $G'$ and $m''_\pm$ in $G''$, such that $m'_+ m'_- =
m''_+ m''_- = 1$. We have
\begin{eqnarray*}
(\pr'\times\pr'')\circ\Deltat\circ\Upsilon\circ\varrhot &=&
(\pr'\times\pr'')\circ(\Upsilon\times\Upsilon)\circ\Deltat\circ \varrhot \\ 
&=&
(\pr'\times\pr'')\circ(\Upsilon\times\Upsilon)\circ(\varrhot
\times\varrhot)\circ \Deltac \\ 
&=& (\varrhot'
\times\varrhot'')\circ
(\pr'\times\pr'')\circ(\Upsilon\times\Upsilon)\circ\Deltac,
\end{eqnarray*}
where the first equality follows from the identity
$(\Upsilon\times\Upsilon)\circ\Delta = \Deltat\circ\Upsilon$, the
second from the identity $(\varrhot\times\varrhot)\circ\Deltac =
\Delta\circ\varrhot$, and the third from Proposition
\ref{prop:ExternalFP}(\ref{prop:ExternalFP_d}). Putting these
together, we obtain
\begin{eqnarray*}
\varrhot^*\Lambda_{g,n}(v'_{\bm'}\times v''_{\bm''}) &= &
((\pr'\circ\Upsilon)^*\varrho'^*\Lambda'_{g,n}(v'_{\bm'}))\cup
((\pr''\circ\Upsilon)^*\varrho''^*\Lambda''_{g,n}(v''_{\bm''})) \\
&=&
((\pr'\times\pr'')\circ(\Upsilon\times\Upsilon)\circ\Delta)^*
(\varrho'^*\Lambda'_{g,n}(v'_{\bm'}))\otimes
\varrho''^*\Lambda''_{g,n}(v''_{\bm''}))) \\
&=&
((\pr'\times\pr'')\circ\Deltat\circ\Upsilon)^*
(\varrho'^*\Lambda'_{g,n}(v'_{\bm'}))\otimes
\varrho''^*\Lambda''_{g,n}(v''_{\bm''})))\\
&=&
\Upsilon^*((\pr'^*\varrho'^*(\Lambda'_{g,n}(v'_{\bm'})))\cup
(\pr''^* \varrho''^*(\Lambda''_{g,n}(v''_{\bm''}))))\\
&=&\Upsilon^*((\pr'^*\Lambda'_{g-1,n+2}(v'_{\bm'},e'_{\alpha[m_+']},
e'_{\alpha[m_-']}))
\cup \\
& & (\pr''^*\Lambda''_{g-1,n+2}(v''_{\bm''},e''_{\beta[m_+'']},
e''_{\beta[m_-'']})) \eta'^{\alpha[m_+']\alpha[m_-']}
\eta''^{\beta[m_+'']\beta[m_-'']}
\end{eqnarray*}
as desired, where $\{ e'_{\alpha[m'_\pm]} \}$ is a basis for
$\ch'_{m'_\pm}$ and  $\{ e''_{\beta[m''_\pm]} \}$ is a basis for
$\ch''_{m''_\pm}$. 

 This completes the case of the loop. The case
of the tree is identical and will be omitted.
\end{proof}

\begin{df}
Let $\mathfrak{G}'=(\ch',\eta',\{\Lambda'_{g,n}\},\vac')$ be a
\Gpcft\ and
$\mathfrak{G}'':=(\ch'',\eta'',\{\Lambda''_{g,n}\},\vac'')$ be a
\Gppcft. Their \emph{external tensor product} $\mathfrak{G}'
\etimes \mathfrak{G}''$ is the $G'\times G''$-\cft\ $(\ch' \otimes
\ch'', \eta' \otimes \eta'', \{ \Lambda_{g,n} \}, \vac' \otimes
\vac'')$, where $\Lambda_{g,n}$ is defined by Equation
(\ref{eq:externaltensor}).
\end{df}

The category of \Gcfts\ also has a tensor product induced from the diagonal
morphism on $\MM_{g,n}$.

\begin{df}
Let $\mathfrak{G}'=(\ch',\eta',\{\Lambda'_{g,n}\},\vac')$ and
$\mathfrak{G}''=(\ch'',\eta'',\{\Lambda''_{g,n}\},\vac'')$ be \Gcfts,  then
consider the tuple  $(\ch,\eta,\{\Lambda_{g,n}\},\vac)$ given by
\begin{enumerate}
\item $\ch = \ch'\btimes\ch''$ as $G$-graded $G$-modules,
\item For all $v'_{m_1}\otimes v''_{m_1}$ in $\ch_{m_1}$ and $v'_{m_2}\otimes
v''_{m_2}$ in $\ch_{m_2}$,
\[\eta(v'_{m_1}\otimes v''_{m_1},v'_{m_2}\otimes v''_{m_2}) :=
       \eta'(v'_{m_1},v'_{m_2}) \eta''(v''_{m_1},v''_{m_2}),
\]
\item $\vac := \vac' \otimes \vac''$, and
\item \[\Lambda_{g,n}(v'_{m_1}\otimes v''_{m_1}, \ldots, v'_{m_n}\otimes
v''_{m_n}) \linebreak[1] := \Lambda'_{g,n}(v'_{m_1}, \ldots,
v'_{m_n}) \cup \Lambda''_{g,n}(v''_{m_1}, \ldots, v''_{m_n}). \]
\end{enumerate}
$(\ch,\eta,\{\Lambda_{g,n}\},\vac)$ is said to be the \emph{tensor
product of the \Gcfts\ } $(\ch',\eta',\{\Lambda'_{g,n}\},\vac')$
and $(\ch'',\eta'',\{\Lambda''_{g,n}\},\vac'')$ and is denoted
$\mathfrak{G}' \itimes \mathfrak{G}''.$
\end{df}

\begin{prop}
The tensor product of two \Gcfts\ is a \Gcft.
\end{prop}
\begin{proof}
The proof follows, first, from the fact that the diagonal morphism
\[
\Delta:\MM_{g,n}(\bm)\to\MM_{g,n}(\bm)\times \MM_{g,n}(\bm)
\]
induces a morphism
\[
H_\bullet(\MM_{g,n}(\bm))\to H_\bullet(\MM_{g,n}(\bm))\otimes
H_\bullet(\MM_{g,n}(\bm)),
\]
which respects the gluing, the $S_n$ actions, and the $G^n$ action, and
second, from the fact that the cup product is induced via pullback of the
diagonal morphism. The definitions of the flat identity and the metric are
easily verified.
\end{proof}

\begin{rem}\label{rem:FactoringTensors}

Let $\ch'$ and $\ch''$ be two $G$-graded $G$-modules. The $G$-module
structure on $\ch'\itimes\ch''$ is induced from the $G\times G$-module
structure on the external tensor product $\ch'\etimes\ch''$ via the diagonal
homomorphism $G\rInto G\times G$. An analogous phenomenon occurs in the
category of $\Gcfts$, where the role of the homomorphism $G\rInto G\times G$
is replaced by a natural inclusion $\MM_{g,n}(\bm)\rInto\M_{g,n}^{G\times
G}(\bm\times\bm)$ for all stable pairs $(g,n)$ and $\bm$ in $G^n$. This
inclusion respects  the actions of $G^n$ and $S_n$ as well as the gluing
morphisms. Consequently, the tensor product in the category of \Gcfts\
``factors through'' the external tensor product. 

This natural inclusion is obtained as follows. The diagonal morphism
$\Delta:\MM_{g,n}(\bm)\rInto\MM_{g,n}(\bm)\times \MM_{g,n}(\bm)$  can be 
written as the composition 
\[ \MM_{g,n}(\bm)\rInto^{\widehat{\Delta}}
\MM_{g,n}(\bm)\times_{\M_{g,n}} \MM_{g,n}(\bm)\rInto^{\hat{j}}
\MM_{g,n}(\bm)\times \MM_{g,n}(\bm),
\]
where $\widehat{\Delta}$ is the diagonal morphism into the fibered
product, and $\hat{j}$ is the obvious inclusion.  However,
$\MM_{g,n}(\bm)\times_{\M_{g,n}} \MM_{g,n}(\bm)$ is isomorphic to
$\M^{G\times G}_{g,n}(\bm\times\bm)$ via $\Upsilon$. Observe that
$\widehat{\Delta}$ and $\widehat{j}$ both preserve the actions of
$S_n$ and $G^n$ and the gluing operations.
\end{rem}

The \Gcft\ $\nc[G]$ is initial among all \Gcfts, in the following
sense.
\begin{prop}\label{prop:groupring}
Let $\mathfrak{G}:=(\ch,\eta,\{\Lambda_{g,n}\},\vac)$ be any
\Gcft.   The tensor product of $\nc[G]$ with $\mathfrak{G}$
satisfies
\[
\nc[G]\itimes\mathfrak{G} \cong \mathfrak{G} \itimes \nc[G]
\cong \mathfrak{G}.
\]
\end{prop}

The proof is immediate from the definition of tensor product.

\subsection{$G$-Frobenius algebras}\label{sec:GFA}

\ 

Recall that a Frobenius algebra is a special \cft. This statement
admits a generalization to \Gcfts\ and $G$-Frobenius algebras, as we will see in 
Theorem~\ref{thm:GFA}.

\begin{df}\label{df:trace}
Let us adopt the notation that $v_m$ is a vector in $\ch_m$ for
any $m\in G$. A tuple $((\ch,\rho),\cdot,\vac,\eta)$ is said to be
a \emph{(non-projective) $G$-Frobenius algebra \cite{Ka,Ka2,Tu}}
provided that the following hold:
\begin{enumerate}
\item ($G$-graded $G$-module) $(\ch,\rho)$ is a
$G$-graded $G$-module.
\item (Self-invariance) For all $\gamma$ in $G$,
$\rho(\gamma):\ch_{\gamma}\to\ch_{\gamma}$ is the identity map.
\item (Metric) $\eta$ is a symmetric, nondegenerate, bilinear form on $\ch$
such that $\eta(v_{m_1},v_{m_2})$ is nonzero only if $m_1 m_2 =
1$.
\item ($G$-graded Multiplication) The binary product
 $(v_{1}, v_{2})\mapsto v_{1}\cdot v_{2}$, called the \emph{multiplication} on $\ch$,
 preserves the
$G$-grading (i.e., the multiplication takes $\ch_{m_1}\otimes
\ch_{m_2}$ to $\ch_{m_1 m_2}$) and is distributive over addition.
\item (Associativity) The multiplication is associative; i.e.,
\[
(v_{1}\cdot v_{2})\cdot v_{3} = v_{1}\cdot (v_{2}\cdot v_{3})
\]
for all $v_{1}$, $v_{2}$, and $v_{3}$ in $\ch$.
\item (Braided Commutativity) The multiplication is invariant with respect
to the braiding:
\[
v_{m_1}\cdot v_{m_2} = (\rho(m_1^{-1}) v_{m_2})\cdot v_{m_1},
\]
for all $m_i \in G$ and all $v_{m_i}\in \ch_{m_i}$ with $i=1,2$.
\item ($G$-equivariance of the Multiplication)
\[ (\rho(\gamma) v_{1})\cdot (\rho(\gamma) v_{2}) =
\rho(\gamma)(v_{1}\cdot v_{2})
\] for all $\gamma$ in $G$,
and all $v_1, v_2 \in \ch$.
\item ($G$-invariance of the Metric)
\[  \eta(\rho(\gamma) v_{1},\rho(\gamma) v_{2}) = \eta(v_{1}, v_{2})
\]
for all $\gamma$ in $G$, and all $v_1, v_2 \in \ch$.
\item (Invariance of the Metric)
\[
\eta(v_{1}\cdot v_{2},v_{3}) = \eta(v_{1}, v_{2}\cdot
v_{3})
\]
for all $v_1,v_2,v_3 \in \ch$.
\item \label{identity} ($G$-invariant Identity) The element $\vac$ in
$\ch_1$ is the identity element under the multiplication, and which satisfies
\[
\rho(\gamma)\vac = \vac
\]
 for all $\gamma$ in $G$.
\item (Trace Axiom)\label{trace} For all $a,b$ in $G$ and $v$ in
$\ch_{[a,b]}$, let $L_v$ denote left multiplication by $v$:
\[
\mathrm{Tr}_{\ch_{a}}(L_v \rho(b^{-1})) =
\mathrm{Tr}_{\ch_{b}}(\rho(a)L_v).
\]

\end{enumerate}

\end{df}

\begin{rem}
When $G$ is the trivial group, a $G$-Frobenius algebra is a
\emph{Frobenius algebra}, a unital, commutative, associative
algebra with an invariant metric. Given a general $G$-Frobenius
algebra $\ch$, there are two ways that one can construct a
Frobenius algebra from it. The first Frobenius algebra is obtained
by considering the subalgebra $\ch_1$. The second approach is to
consider $\chb$, the algebra of $G$-coinvariants of $\ch$, with
its induced multiplication and identity. The metric on $\ch$
induces a metric on $\chb$ which makes $\chb$ into a Frobenius
algebra.
\end{rem}

\begin{rem}
If $\ch$ is a $G$-Frobenius algebra, then it follows from the axioms of a $G$-Frobenius algebra that the action of the braid group
on the multiplication factors through the symmetric group. More
precisely, let $\mu:\ch^{\otimes 3}\to \nc$ be given by
$\mu(v_{m_1},v_{m_2},v_{m_3}) := \eta(v_{m_1}\cdot v_{m_2}, v_{m_3})$
and let $b_1, b_2$ denote the generators of the braid group $B_3$, then
$\mu\circ b_i\circ b_i = \mu$ for all $i=1,2$.
\end{rem}

\begin{thm}\label{thm:GFA}
Let $((\ch, \rho),\eta,\{\Lambda_{g,n}\},\vac)$ be a \Gcft.
Define a multiplication $\cdot$ on $\ch$ as follows:  For any $m_1,m_2 \in
G$, let $m_3 = (m_1 m_2)^{-1}$. For all $v_{m_1}$ in $\ch_{m_1}$ and
$v_{m_2}$ in $\ch_{m_2}$, define
\[
v_{m_1} \cdot v_{m_2} := \int_{\db{\xi(m_1,m_2,m_3)}}
\Lambda_{0,3}(v_{m_1},v_{m_2},e_\alpha) \eta^{\alpha \beta}
f_\beta,
\]
where $\{ e_\alpha \}$ is a basis for 
  $\ch_{m_3 }$,
$\{ f_\beta \}$ is a basis for 
    $\ch_{m_3^{-1}}$, and
$\eta^{\alpha \beta}$ is the inverse of the metric in those bases.

The tuple $((\ch,\rho),\cdot, \vac,\eta)$ is a \emph{$G$-Frobenius
algebra}.
\end{thm}
\begin{proof}
The $G$-module $(\ch,\rho)$, the metric $\eta$, and
the identity element $\vac$ in the \Gcft\ are the same for
the $G$-Frobenius algebra.

The invariance of the metric follows from the fact that
\[
s \xi(m_1,m_2,m_3) = \xi(m_2,m_3,m_1),
\]
where $s$ is the isomorphism induced from the cyclic permutation
in $S_3$ (this is proved in Lemma ~\ref{lemma}).

Notice also that since $\xi(m_1,m_2,m_3)$ is empty unless
$m_1m_2m_3=1$, the product is naturally graded.

The product is not commutative, in general, because
$\xi(m_1,m_2,m_3)\not=\xi(m_2,m_1,m_3)$. However, it is braided
commutative, because
\begin{equation}\label{eq:br-com}
\xi(m_1,m_2,m_3) = b^{-1}_1\xi(m_1 m_2 m_1^{-1},m_1,m_3) = \sigma
\rho_1(m_1) \xi(m_1 m_2 m_1^{-1},m_1,m_3),
\end{equation}
with $\sigma$ the transposition $(1,2)\in S_3$, as shown in Lemma
~\ref{lemma}. The relation (\ref{eq:br-com}) on $\xi$ implies the
braided commutativity via the equation

\begin{eqnarray*}
v_{m_1} \cdot v_{m_2} &=& \int_{\db{\xi(m_1,m_2,m_3)}}
\Lambda_{0,3}(v_{m_1},v_{m_2},e_\alpha) \eta^{\alpha \beta}
f_\beta\\ &=&\int_{\db{(\sigma\rho_1(m_1))\xi(m_1 m_2
m_1^{-1},m_1,m_3)}} \Lambda_{0,3}(v_{m_1},v_{m_2},e_\alpha)
\eta^{\alpha \beta}f_\beta\\ &=&\int_{\db{\xi(m_1 m_2
m_1^{-1},m_1,m_3)}}
(\sigma\rho_1(m_1))^*(\Lambda_{0,3}(v_{m_1},v_{m_2},e_\alpha))
\eta^{\alpha \beta} f_\beta\\ &=&\int_{\db{\xi(m_1 m_2
m_1^{-1},m_1,m_3)}} \rho_1(m_1)^*\sigma^*
(\Lambda_{0,3}(v_{m_1},v_{m_2},e_\alpha)) \eta^{\alpha \beta}
f_\beta\\ &=&\int_{\db{\xi(m_1 m_2 m_1^{-1},m_1,m_3)}}
\rho_1(m_1)^*(\Lambda_{0,3}(v_{m_2},v_{m_1},e_\alpha))
\eta^{\alpha \beta} f_\beta\\ &=&\int_{\db{\xi(m_1 m_2
m_1^{-1},m_1,m_3)}}
\Lambda_{0,3}(\rho_1(m_1^{-1})v_{m_2},v_{m_1},e_\alpha)
\eta^{\alpha \beta} f_\beta\\ &=&(\rho(m_1^{-1})v_{m_2}) \cdot
v_{m_1}.
\end{eqnarray*}

Again, using the braided commutativity for the classes $\xi$ and
the invariance of $\bf{1}$,  we can show that $\bf{1}$ is indeed a
unit for the multiplication, since

\begin{eqnarray*}
v_{m_1} \cdot \vac  &=&  \int_{\db{\xi(m_1,1,m_1^{-1})}}
\Lambda_{0,3}(v_{m_1},{\vac },e_\alpha) \eta^{\alpha \beta}
f_\beta\\ &=& \int_{\db{\xi(m_1,m_1^{-1},1)}}
\Lambda_{0,3}(v_{m_1},e_\alpha,{\vac }) \eta^{\alpha \beta}
f_\beta\\ &=& \eta(v_{m_1},\vac) \eta^{\alpha \beta} f_\beta=
v_{m_1},
\end{eqnarray*}
where we introduced a basis $(e_\alpha)$ of   $\ch_{m_1^{-1}}$ and
a basis $(f_\beta)$ of $\ch_{m_1}$.

The property that ${\vac }$ is a unit implies that the invariance
of the metric follows from
\begin{equation}\label{eq:metr}
\eta(v_{m_1},v_{m_2})= \eta(v_{m_1}\cdot v_{m_2},{\vac } ).
\end{equation}

Equation (\ref{eq:metr}) in turn follows from
\begin{eqnarray*}
\eta(v_{m_1} \cdot v_{m_2},{\vac })&=&
 \int_{\db{\xi(m_1,m_2,m_3)}}
\Lambda_{0,3}(v_{m_1},v_{m_2},e_\alpha) \eta^{\alpha \beta}
 \int_{\db{\xi(m_3,1,1)}}
\Lambda_{0,3}(f_{\beta},{\vac } ,{\vac }) \\
&=&\int_{\db{\xi(m_1,m_2,m_3)}}
\Lambda_{0,3}(v_{m_1},v_{m_2},e_\alpha) \eta^{\alpha \beta}
 \eta(f_{\beta} ,{\vac })\\
&=&\int_{\db{\xi(m_1,m_2,m_3)}} \Lambda_{0,3}(v_{m_1},v_{m_2},{\vac})\\
 &=&\eta(v_{m_1}, v_{m_2} ),
\end{eqnarray*}
where we use the notation $m_3:=(m_1m_2)^{-1}$, and we let
$\{e_\alpha\}$ be a basis of $\ch_{m_3}$ and $\{f_\beta\}$ be a
basis of $\ch_{m_3^{-1}}$.

The $\rho(\gamma)$-invariance of $\ch_{\gamma}$ follows from the
second part of Lemma \ref{lemma}:
\begin{eqnarray*}
\rho(\gamma)v_{\gamma} &=&  \rho(\gamma)v_{\gamma}\cdot {\vac } =
\int_{\db{\xi(\gamma, 1,\gamma^{-1})}}
\Lambda_{0,3}(\rho(\gamma)v_{\gamma},\vac,e_{\alpha})\eta^{\alpha
\beta} f_{\beta}\\ &=& \int_{\db{\rho(\gamma,1,1)\xi(\gamma,
1,\gamma^{-1})}}
\Lambda_{0,3}(\rho(\gamma)v_{\gamma},\vac ,e_{\alpha}) \eta^{\alpha
\beta}  f_{\beta}\\ &=&\int_{\db{\xi(\gamma, 1,\gamma^{-1})}}
\Lambda_{0,3}(v_{\gamma},\vac ,e_{\alpha})\eta^{\alpha \beta}
f_{\beta}\\ &=& v_{\gamma}.
\end{eqnarray*}
Again we use bases $\{e_\alpha\}$ of $\ch_{\gamma^{-1}}$ and
$\{f_{\beta}\}$ of $\ch_{\gamma}$.

The self invariance, together with the invariance of the metric,
imply the symmetry of the metric:
\[\eta(v_{m},v_{m^{-1}}) = \eta(v_{m}v_{m^{-1}},\vac)
=\eta(\rho(m^{-1})(v_{m^{-1}})v_m,\vac)= \eta(v_{m^{-1}},v_m).
\]

The $G$-invariance of the metric follows from the $G^n$-invariance
of $\Lambda$ and the $\rho$-invariance of the unit ${\vac }$
via
\begin{eqnarray*}
\eta(\rho(\gamma)v_{m_1},\rho(\gamma) v_{m_2} )&=&
\int_{\db{\xi(\gamma^{-1} m_1\gamma,\gamma^{-1} m_2\gamma, 1)}}
\Lambda_{0,3}(\rho(\gamma)v_{m_1},\rho(\gamma)v_{m_2},{\vac })\\
&=&\int_{\db{\rho(\gamma,\gamma,\gamma)\xi (m_1, m_2, 1)}}
\Lambda_{0,3}(\rho(\gamma)v_{m_1},\rho(\gamma)v_{m_2},\rho(\gamma){\vac})\\ 
&=&\int_{\db{\xi( m_1, m_2,
1)}}\Lambda_{0,3}(v_{m_1},v_{m_2},{\vac })\\ &=& \eta(v_{m_1},
v_{m_2}),
\end{eqnarray*}
where we used the first property of $\xi$ of Lemma \ref{lemma}.

The above in turn gives the $G$-equivariance of the multiplication
\begin{eqnarray*}
\rho(\gamma)v_{m_1}\cdot \rho(\gamma) v_{m_2}&=&
\int_{\db{\xi(\gamma^{-1} m_1\gamma,\gamma^{-1} m_2\gamma,
(\gamma^{-1} m_3 \gamma)^{-1})}}
\Lambda_{0,3}(\rho(\gamma)v_{m_1},\rho(\gamma)v_{m_2},e_{\alpha})
\eta^{\alpha \beta} f_{\beta}\\
&=&\int_{\db{\rho(\gamma,\gamma,\gamma)(\xi( m_1, m_2, m_3))}}
\Lambda_{0,3}(\rho(\gamma)v_{m_1},\rho(\gamma)v_{m_2},\rho(\gamma)e'_{\alpha})
\eta^{\alpha \beta} \rho(\gamma)f'_{\beta}\\ &=&\int_{\db{\xi(
m_1, m_2,m_3}}
\Lambda_{0,3}(v_{m_1},v_{m_2},e'_{\alpha})\eta^{\prime \alpha
\beta} \rho(\gamma)f'_{\beta}\\ &=&\rho(\gamma)(v_{m_1}\cdot
v_{m_2}),
\end{eqnarray*}
where we used $m_3:=(m_1m_2)^{-1}$, a basis $\{e_\alpha\}$ of
$\ch_{\gamma^{-1}m_3\gamma}$, $\{f_\beta\}$ of
$\ch_{\gamma^{-1}m_3^{-1}\gamma}$, and the transformed bases
$\{e'_{\alpha}:= \rho(\gamma^{-1})e_{\alpha}\}$ of $\ch_{m_3}$ and
$\{f'_{\beta}:= \rho(\gamma^{-1})f_{\beta}\}$ of $\ch_{m_3^{-1}}$.
Also, we used the notation $\eta'^{ \alpha \beta}$ for the inverse
metric of $\eta_{\alpha \beta}=\eta(e'_\alpha, f'_\beta)$, the
$G$-invariance of the metric $\eta'_{\alpha
  \beta}=\eta_{\alpha\beta}$, and the first property of Lemma
\ref{lemma}.

Associativity follows from Lemma \ref{ass} in the following
way:
\begin{eqnarray*}
(v_{m_1}\cdot v_{m_2})\cdot v_{m_3} &=&\int_{\db{\xi( m_1,
m_2,m_+}} \Lambda_{0,3}(v_{m_1},v_{m_2},e_{\alpha})\eta^{\alpha
\beta} \int_{\db{\xi( m_-,m_3,m_4)}}
\Lambda_{0,3}(f_{\beta},v_{m_3},k_{\gamma})\eta^{\gamma \delta}
l_{\delta}\\ &=& \int_{\db{\xi( m_1, m_2,m_+)\times \xi(m_-, m_3,
m_4)}} \varrho^*\Lambda_{0,4}(v_{m_1},v_{m_2},v_{m_3},k_{\gamma})
\eta^{\gamma \delta}l_{\delta}\\ &=&\int_{\db{\varrho(\xi( m_1,
m_2,m_+)\times \xi(m_-, m_3, m_4)}}
\Lambda_{0,4}(v_{m_1},v_{m_2},v_{m_3},k_{\gamma}) \eta^{\gamma
\delta}l_{\delta}\\ &=&\int_{\db{\varrho'(\xi(m_4 , m_1,m'_+)
\times \xi(m'_-,m_2,m_3))}}
\Lambda_{0,4}(v_{m_1},v_{m_2},v_{m_3},k_{\gamma})\eta^{\gamma
\delta} l_{\delta}\\ &=&\int_{\db{\xi(m_4 , m_1,m'_+)\times
\xi(m'_-,m_2,m_3)}} \varrho^{\prime
*}\Lambda_{0,4}(v_{m_1},v_{m_2},v_{m_3},k_{\gamma}) \eta^{\gamma
\delta} l_{\delta}\\ &=&\int_{\db{\xi(m_4 , m_1,m'_+)}}
\Lambda_{0,3}(k_{\gamma},v_{m_1},e_{\alpha})\eta^{\alpha \beta}
\int_{\db{ \xi(m'_-,m_2,m_3)}}
\Lambda_{0,3}(f_{\beta},v_{m_2},v_{m_3})\eta^{\gamma \delta}
l_{\delta}\\ &=&\int_{\db{\xi( m_2, m_3,m'_-)}}
\Lambda_{0,3}(v_{m_2},v_{m_3},f_{\beta})\eta^{\beta \alpha}
\int_{\db{\xi(m_1,m'_+,m_4}}
\Lambda_{0,3}(v_{m_1},e_{\alpha},k_{\gamma})\eta^{\gamma \delta}
l_{\delta}\\ &=&v_{m_1}\cdot (v_{m_2}\cdot v_{m_3}),
\end{eqnarray*}
where we used the $S_3$-invariance of $\Lambda$, the fourth
property of \ref{lemma}, and the symmetry of the metric. Also, we
introduced the notation $m_4=(m_1m_2m_3)^{-1}$ and used the
notations of  Lemma \ref{ass} for
$m_{\pm},m'_{\pm},\varrho,\varrho'$ ; i.e., $m_+:=(m_1m_2)^{-1},
m_-:= m_1 m_2, m'_+=m_2 m_3$, and $m'_-:=(m_2m_3)^{-1}$.
Furthermore, $\{e_{\alpha}\}$ is a basis of $\ch_{m_+}$,
$\{f_{\beta}\}$ is a basis of $\ch_{m_-}$, $\{k_{\gamma}\}$ is a
basis of $\ch_{m_4}$, and $\{l_{\delta}\}$ is a basis of
$\ch_{m_4^{-1}}$.

Lastly, the proof of the trace axiom follows using Lemma
\ref{lm:trace}:
\begin{eqnarray*}
\mathrm{Tr}_{\ch_{a}}(L_v \rho(b^{-1}))&= &
\eta(\eta^{\alpha\beta} f_{\beta},
v_{aba^{-1}b^{-1}}\cdot\rho(b^{-1}) e_{\alpha})\\
&=&\eta(\eta^{\alpha\beta}f_{\beta},
\int_{\db{\xi(m_1,bab^{-1},a^{-1})}}
\Lambda_{0,3}(v_{m_1},\rho(b^{-1})
e_{\alpha},f_{\gamma})\eta^{\gamma\delta} e_{\delta})\\
&=&\int_{\db{\xi(a^{-1},a,1)}}
\Lambda_{0,3}(f_{\beta},e_{\delta},{\vac }) \eta^{\alpha\beta}
\int_{\db{\xi(m_1,bab^{-1},a^{-1})}}
\Lambda_{0,3}(v_{m_1},\rho(b^{-1})
e_{\alpha},f_{\gamma})\eta^{\gamma\delta}\\
&=&\int_{\db{\xi(m_1,bab^{-1},a^{-1})}}
\Lambda_{0,3}(k_{\lambda},\rho(b^{-1})e_{\alpha},f_{\beta})
\eta^{\alpha\beta} \int_{\db{\xi(1,m_1,m_1^{-1})}}
\Lambda_{0,3}({\vac },v_{m_1},l_{\mu})\eta^{\lambda\mu}\\
&=&\int_{\db{\varrho_{b*}\xi(m_1,bab^{-1},a^{-1})}}
\Lambda_{1,1}(k_{\lambda}) \int_{\db{\xi(1,m_1,m_1^{-1})}}
\Lambda_{0,3}({\vac },v_{m_1},l_{\mu})\eta^{\lambda\mu}\\
&=&\int_{\db{\varrho_{a*}\xi(m_1,b,ab^{-1}a^{-1})}}
\Lambda_{1,1}(k_{\lambda}) \int_{\db{\xi(1,m_1,m_1^{-1})}}
\Lambda_{0,3}({\vac },v_{m_1},l_{\mu})\eta^{\lambda\mu}\\ &=&
\int_{\db{\xi(m_1,b,ab^{-1}a^{-1})}}
\Lambda_{0,3}(k_{\lambda},g_{\gamma},\rho(a)^{-1}h_{\delta})
\eta^{\gamma\delta} \int_{\db{\xi(1,aba^{-1}b^{-1},m_1^{-1})}}
\Lambda_{0,3}({\vac },v_{m_1},l_{\mu})\eta^{\lambda\mu}\\
&=&\int_{\db{\xi(b^{-1},b,1)}}
\Lambda_{0,3}(h_{\tau},g_{\gamma},{\vac }) \eta^{\sigma\tau}
\int_{\db{\xi(m_1,b,ab^{-1}a^{-1})}}
\Lambda_{0,3}(v_{m_1},g_{\sigma},\rho(a^{-1})h_{\delta})\eta^{\gamma\delta}\\
&=&\mathrm{Tr}_{\ch_{b}}(\rho(a)L_v ),
\end{eqnarray*}
where we used Lemma \ref{ass}, as well as Lemma \ref{lm:trace} with
its notation for the maps $\varrho_a,\varrho_b$ and $m_1 = [a,b]$,
and introduced the bases $\{e_\alpha\}$ of $\ch_{a}$,
$\{f_\beta\}$ of $\ch_{a^{-1}}$, $\{g_\gamma\}$ of $\ch_b$,
$\{h_\delta\}$ of $\ch_{b^{-1}}$, $\{k_\lambda\}$ of $\ch_{m_1}$,
and $\{l_{\mu}\}$ of $\ch_{m_1^{-1}}$.

\end{proof}

We can now justify naming the \Gcft\ $\nc[G]$ of Example~\ref{ex:grpring} the  \emph{group ring \Gcft}.

\begin{prop}
In the group ring \Gcft, the metric $\eta$ on $\ch=\bigoplus_{g\in G} \nc$
satisfies
\[
\eta(e_{m_1}, e_{m_2}) = \delta_{m_1,m_2^{-1}}
\]
for all $m_1,m_2$ in $G$.

The multiplication is given by
\[
e_{m_1} \cdot e_{m_2} = e_{m_1 m_2}
\]
for all $m_1, m_2$ in $G$. The identity element is $\vac := e_1$.
The resulting $G$-Frobenius algebra is isomorphic to the group ring $\nc[G]$.
\end{prop}
\begin{proof}
The multiplication operation is
\[
e_{m_1} \cdot e_{m_2} =
\int_{\db{\xi(m_1,m_2,(m_1 m_2)^{-1})}} \Lambda_{0,3}(e_{m_1},e_{m_2},e_{(m_1
m_2)^{-1}}) e_{m_1 m_2}   = e_{m_1 m_2}.
\]
The metric and identity element follow by a similar calculation.
\end{proof}

\section{\cfts\ and Quotients of \Gcfts}\label{5}

\ 

In this subsection, we explain how to obtain a \cft\ from a \Gcft\
by taking the appropriate quotient with respect to $G$.
Geometrically, going from a \Gcft\ to a \cft\ corresponds to going
from $\MM_{g,n}$ to $\M_{g,n}$, where the $\Lambda_{g,n}$ are
allowed to only act upon elements of $\chb$. We perform this
procedure in two steps. The first step is to go from $\MM_{g,n}$
to $\M_{g,n} (\BG)$. The second step is to go from $\M_{g,n}(\BG)$
to $\M_{g,n}$.

\subsection{From $\MM_{g,n}$ to $\M_{g,n}(\BG)$}

\ 

We begin with a useful lemma.

\begin{lm}\label{lm:useful}
For all $\bmb$ in $\Gb^n$, the forgetful morphism
$\stt_{\bmb}:\MM_{g,n}(\bmb)\rTo\M_{g,n}(\BG;\bmb)$ induces a ring isomorphism
$\stt_{\bmb}^*: H^\bullet(\M_{g,n}(\BG;\bmb)) \rTo
H^\bullet(\MM_{g,n}(\bmb))^{G^n}$.
\end{lm}
\begin{proof}

Let ${\cc}$ be the constant sheaf of complex numbers on
$\MM_{g,n}(\bmb)$, and let ${\cc}'$ be the constant sheaf on
$\M_{g,n}(\BG;\bmb)$.

Since $\stt$ is finite, the Leray spectral sequence degenerates,
giving  $$H^p(\MM_{g,n}(\bmb),{\cc}) =
H^p(\M_{g,n}(\BG;\bmb),\stt_*({\cc})).$$     Since these sheaves are
all sheaves of vector spaces over $\nc$, they are all
divisible, hence the coinvariant map $\pi_{G^n}$ is well defined
and preserves invariants; i.e.,  if $i: (\stt_*{\cc})^{G^n} \rTo
\stt_*{\cc}$ is the natural inclusion, then $\pi_{G^n} \circ i = \id$.
Thus taking $G^n$-invariants is the same as applying the map
$\pi_{G^n}$, and is exact.  So a general homological argument
gives that $$(H^p(\M_{g,n}(\BG;\bmb),\stt_*{\cc}))^{G^n} =
H^p(\M_{g,n}(\BG;\bmb), (st_*{\cc})^{G^n}),$$  and we have
$$(H^p(\MM_{g,n}(\bmb),{\cc}))^{G^n} = H^p(\M_{g,n}(\BG;\bmb),
(\stt_*{\cc})^{G^n}).$$

On the other hand, we have $\stt^*({\cc}') = {\cc}$, so by adjointness we
have a map $j: {\cc}' \rTo \stt_* {\cc}$.  Composing with $\pi_{G^n}$, we
get a map of sheaves $\pi_{G^n}\circ j: {\cc}' \rTo (\stt_*{\cc})^{G^n}$.
On each fiber this map is an isomorphism, since for a fixed
admissible cover $E\rTo C$ the fiber $F:=\stt^{-1}([E\rTo C])$ is
a disjoint union of points with transitive ${G^n}$-action inducing
the ${G^n}$-action on $\stt_*{\cc} = \bigoplus_{f \in F} {\cc}'$, and $j$
is just given by $q \mapsto (q,q,...,q)$.  Some straightforward
work shows that for any vector space $V$ and any set $F$ with
transitive ${G^n}$-action, the vector space $V\times F$ has as its
${G^n}$-invariants exactly the image of the map $j:V \rTo V\times
F$, taking $v$ to $(v,v,...,v)$.  In particular, this holds for
$V={\cc}'$.  Since the fiber $F$ and the ${G^n}$-action on $F$ are
unchanged under small deformation, this shows that the morphism of
sheaves $\pi_{G^n}\circ j$ induces an isomorphism on stalks, and
thus is an isomorphism of sheaves. So we have
\begin{eqnarray*}
H^p(\M_{g,n}(\BG;\bmb),{\cc}') &=& H^p(\M_{g,n}(\BG;\bmb),(\stt_*({\cc}))^{G^n}) \\
&=& H^p(\M_{g,n}(\BG;\bmb),\stt_*{\cc})^{G^n} = H^p(\MM_{g,n},{\cc})^{G^n}
\end{eqnarray*}
as desired.

\end{proof}

\begin{prop} \label{prop:Lambdah}
Let $(\ch,\eta,\Lambda_{g,n},\vac)$ be a \Gcft.  
There exist uniquely-determined 
classes
$\Lambdah_{g,n}$ in $\bigoplus_{\bmb\in\Gb^n}
H^\bullet(\M_{g,n}(\BG;\bmb))\otimes \chb_\bmb^*$ such that
\[
\stt^*_\bmb\Lambdah_{g,n}(v_\bmb) = \Lambda_{g,n}(v_\bmb)
\]
for all $v_\bmb$ in $\chb_\bmb$.
\end{prop}

\begin{proof}
Consider $v_\bmb$ in $\chb_\bmb$ for $\bmb$ in $\Gb^n$.  For all $\bg$
in $G^n$ we have
\[
\rho(\bg)^*(\Lambda_{g,n}(v_\bmb)) =
\Lambda_{g,n}(\rho(\bg^{-1})_*v_\bmb)) = \Lambda_{g,n}(v_\bmb)),
\]
where the first equality is by the (diagonal) $G^n$-invariance of
$\Lambda_{g,n}$  and the second is by the definition of $\chb$.
Therefore, $\Lambda_{g,n}(v_\bmb)$ belongs to
$H^\bullet(\M_{g,n}(\bmb))^{G^n}$, and we are done by the previous
lemma.
\end{proof}

Fix an element $m_+$ in $G$ and let $m_- := m_+^{-1}$. To each such choice,
we have the following associated commutative diagram, which we will use extensively  hereafter, and whose morphisms and other terms we explain below:

\begin{equation}
 \begin{diagram} \label{eq:MMtoM(BG)}
 & &
F^G_{\Gammah} & \rTo^{\mut} &  \MM_{\Gammah} &
 \rInto^{\itt} & \MM_{g,n}(\bmb) \\
 & \ruTo_{\rt} \\
 \MM_{\Gammatc} & & \dTo^\prt & & \dTo^{\stt'} & &
 \dTo^\stt \\
 & \\
\dTo_{\stt''}  & & F_{\Gammah}(\BG) & \rTo^{\muh} &  \M_{\Gammah}(\BG) &
 \rInto^{\ih} & \M_{g,n}(\BG;\bmb)\\
& \ldTo_{\rh}
  & & & & & \\
\M_{\Gammahc} & & \dTo^\prh & & \dTo^{\sth'} & &  \dTo^\sth \\
& \rdTo^{\sth''} \\
  & & \M_{\Gammac} & \rTo^{\mu} &  \M_{\Gamma} &
 \rInto^{i} & \M_{g,n}
\end{diagram}
\end{equation}

The above diagram has two cases. The first case corresponds to the situation
where all graphs are decorated stable graphs of genus $g$ with $n$ tails
which are trees of the form
\begin{equation}\label{eq:treea}
\Gammat   = \treemb \quad \mathrm{and} \quad \Gammatc =
\treembcut
\end{equation}

\begin{equation}\label{eq:treeb}
\Gammah = \treembb \quad \mathrm{and} \quad \Gammahc =
\treembbcut
\end{equation}

and
\begin{equation}\label{eq:treec}
\Gamma = \tree \quad \mathrm{and} \quad \Gammac = \treecut,
\end{equation}
where $N_+ := \{ i_1,\ldots,i_{n_+} \}$, 
is the index set of the labels for the tails
on the left half of each graph above,
$N_- := \{ j_1,\ldots,j_{n_-} \}$,
is the index set of the labels for the tails
on the right half of each graph above,
$N_+\sqcup N_- = \{1,\ldots,n\}$, and $g_+ + g_- = g$.

The second case corresponds to the situation where all graphs are decorated
stable graphs of genus $g$ with $n$ tails which are loops of the following
form:
\begin{equation}\label{eq:loopa}
\Gammat = \loopmb \quad \mathrm{and} \quad \Gammatc =
\loopmbcut
\end{equation}
\begin{equation}\label{eq:loopb}
\Gammah = \loopmbb \quad \mathrm{and} \quad \Gammahc =
\loopmbbcut
\end{equation}
and
\begin{equation}\label{eq:loopc}
\Gamma = \loopp \quad \mathrm{and} \quad \Gammac = \loopcut.
\end{equation}

Note that in both cases, the graph $\Gammat$ has tails labeled by
conjugacy classes $\mb_i$, but its one edge is labeled by a specific
choice of $m_+$ and $m_-$.

Now we explain the various terms and morphisms.  $\M_\Gamma$ is the
closure of the locus in $\M_{g,n}$ whose dual graph is $\Gamma$, $i$
is the inclusion morphism, and $\mu$ is the normalization morphism
associated to cutting the internal edge of $\Gamma$. Similarly,
$\M_\Gammah(\BG)$ denotes the closure of the locus in
$\M_{g,n}(\BG;\bmb)$ whose associated dual graph has tails decorated
by $\bmb$ and whose monodromies around one side of the node lie in
$\mb_+$ and whose monodromies around the other side of the node lie in
$\mb_-$. The morphism $\ih$ is the inclusion, and $F_\Gammah(\BG)$ is
the fibered product $\M_{\Gammac}\times_{\M_{\Gamma}}
\M_{\Gammah}(\BG)$.  The morphisms $\muh$ and $\prh$ are the canonical
projections of the fibered product.  To explain $\rh$, we first note
that $F_{\Gammah}(\BG)$ is the the stack of triples consisting of a
cut curve $C'$ in $\M_{\Gammac}$, an admissible $G$-cover $E\to C$ in
$\M_{\Gammah}(\BG)$, and an isomorphism $\alpha$ in $\M_{\Gamma}$ from
the glued curve $\mu(C')$ to $C$.  The morphism $\rh$ takes such a
triple to the pullback of $E$ along the composition $\alpha \circ
\mu$.  The morphisms $\sth$, $\sth'$, and $\sth''$ simply forget their
respective twisted curve structures.

Similarly, $\MM_\Gammah$ is the closure of the locus of pointed
$G$-covers with dual graph $\Gammah$, so all tails are labeled by
conjugacy classes $\mb_i$; and $\MM_{\Gammatc}$ is the closure of
the locus of pointed $G$-covers with dual graph $\Gammatc$, so
their tails are labeled with conjugacy classes $\mb_i$, but on the
two sides of the node their monodromies are the specific group
elements $m_+$ and $m_-$.  The morphisms $\stt'$ and $\stt''$
simply forget the marked points in the $G$-cover.

The stack $F^G_{\Gammah}$ is the fibered product
$F_\Gammah(\BG)\times_{\M_{\Gammah}(\BG)} \MM_{\Gamma} =
\M_\Gammac\times_{\M_{\Gamma}} \MM_{\Gamma}$, and the morphisms $\mut$
and $\prt$ are the canonical projections.  The morphism $\rt$ is
induced by the pair of the gluing map $\varrho:\MM_{\Gammatc} \rTo
\MM_{\Gammah}$ and the map $\sth'' \circ \stt'': \MM_{\Gammatc} \rTo
\M_{\Gammac}$ (actually the gluing map has as its target $\MM_{g,n}$,
but it factors through the substack $\MM_{\Gammah}$).  In particular,
we can write the gluing morphism on $\MM$ as
\begin{equation}
\varrhot_\Gammat = \itt\circ\mut\circ\rt,
\end{equation}
while the corresponding gluing morphism on $\M$ can be written as
\begin{equation}
\varrho_\Gamma = i\circ\mu.
\end{equation}

\begin{rem}
The morphisms $i$, $\ih$, $\itt$ are regular embeddings. The remaining
morphisms in the diagram are both flat and proper.
\end{rem}

\begin{nota}
For all $\mb$ in $\Gb$, let $|C(\mb)|$ denote the order of the subgroup
$C(m')$ of $G$ for any $m'$ in $\mb$, as it is independent of
the choice of $m'$.
\end{nota}

\begin{thm}\label{thm:cutting_edges}
Let $\{\Lambdah_{g,n}\}$ be a collection of classes associated to
a \Gcft\ $\{\Lambda_{g,n}\}$, as in Proposition~\ref{prop:Lambdah}.
Fix any conjugacy class $\mb_+$, and let $\mb_- := \mb_+^{-1}$. Let
$\Gammah$ be a decorated stable graph of genus $g$ with $n$-tails
which is either a tree, as in Equation (\ref{eq:treeb}), or a
loop, as in Equation (\ref{eq:loopb}). Let $v_\bmb$ belong to
$\chb_\bmb$.

When $\Gammah$ is a tree then
\begin{equation} \label{eq:cutting_tree}
\rh_*\muh^*\ih^*\Lambdah_{g,n}(v_\bmb) = \frac{\deg(\sth')}{\deg(\sth'')}
            \sum_{\beta[\mb_+],\beta[\mb_-]}
\Lambdah_{g_+,n_++1}(v_{\bmb_{N_+}},e_{\beta[\mb_+]})
\etah^{\beta[\mb_+],\beta[\mb_-]}
\Lambdah_{g_-,n_-+1}(e_{\beta[\mb_-]},v_{\bmb_{N_-}}),
\end{equation}
where $N_+\sqcup N_- = \{ 1,\ldots,n\}$ is the partition
corresponding to the tree, $n_\pm = |N_\pm|$ and
$v_{\bmb_{N_\pm}}$ denotes the $n_\pm$-tuple $\prod_{i\in N_\pm}
v_{m_i}$, the collection $\{ e_{\beta[\mb_\pm]} \}$ is a basis of
$\chb_{\bm_\pm}$, and $g_++g_- = g$.   And $\etah^{\beta[\mb_+]
\beta[\mb_-]}$ is the inverse of the metric $\etah$ on $\chb$ in
the basis $\{ e_{\beta[\mb_\pm]} \}$, where
\begin{equation}\label{eq:etah}
\etah(v_{\mb_+},v_{\mb_-}) := |C(\mb_+)| \etab(v_{\mb_+},v_{\mb_-})
\end{equation}
for all $v_{\mb_\pm}$ in $\chb_{\mb_\pm}$.

When $\Gammah$ is a loop then
\begin{equation} \label{eq:cutting_loop}
\rh_*\muh^*\ih^*\Lambdah_{g,n}(v_\bmb) = \frac{\deg(\sth')}{\deg(\sth'')}
 \sum_{\beta[\mb_+],\beta[\mb_-]}
\Lambdah_{g-1,n+2}(v_{\bmb},e_{\beta[\mb_+]}, e_{\beta[\mb_-]})
\etah^{\beta[\mb_+],\beta[\mb_-]}.
\end{equation}
In either case, denote the right hand side of equations
(\ref{eq:cutting_tree}) and (\ref{eq:cutting_loop}) by
$\frac{\deg(\sth')}{\deg(\sth'')} \Lambdah_{\Gammahc}$.
\end{thm}

\begin{rem}
This theorem suggests that $\Lambdah_{g,n}$ should be regarded as an
analog of the virtual class $c^{1/r}_{g,n}$ on $\M_{g,n}^{1/r}$, the moduli
stack of $r$-spin curves \cite{JKV}. Equations (\ref{eq:cutting_tree}) and
(\ref{eq:cutting_loop}) should be regarded as an analog of the Cutting-Edges
axiom.
\end{rem}

\begin{proof}[Proof: (of Theorem~\ref{thm:cutting_edges})]
Let $m_+$ be any representative of the conjugacy
class $\mb_+$ and let $m_- := m_+^{-1}$.  Consider the associated
commuting diagram (\ref{eq:MMtoM(BG)}) and graphs (\ref{eq:treea}) to
(\ref{eq:loopc}).

For $v_\bmb$ in $\chb_\bmb$, let $$I:=\stt''_* \varrhot_\Gammat^* \Lambda_{g,n}(v_\bmb).$$
We have
\begin{eqnarray*}
\varrhot_\Gammat^*\Lambda_{g,n}(v_\bmb) &=& \varrhot_\Gammat^*\stt^*\Lambdah_{g,n}(v_\bmb) \\
&=& (\stt\circ\itt\circ\mut\circ\rt)^* \Lambdah_{g,n}(v_\bmb) \\
&=& (\ih\circ\muh\circ\prt\circ\rt)^* \Lambdah_{g,n}(v_\bmb).
\end{eqnarray*}
Therefore,
\begin{eqnarray}
I &=& \stt''_* (\ih\circ\muh\circ\prt\circ\rt)^* \Lambdah_{g,n}(v_\bmb)\nonumber \\
&=& (\rh\circ\prt\circ\rt)_* (\ih\circ\muh\circ\prt\circ\rt)^*
\Lambdah_{g,n}(v_\bmb)\nonumber\\
&=& \deg(\prt\circ\rt) \rh_*\muh^*\ih^*\Lambdah_{g,n}(v_\bmb)\label{eq:I_one}
\end{eqnarray}
because $\prt\circ\rt$ is finite and surjective.

For all $m$ in $G$, let $\{ e_{\alpha[m]} \}$ be a basis for
$\ch_m$ such that $\{ e_{\alpha[m]} \}$ is the disjoint union of a
basis $\{ e_{\mu[m]} \}$ for $\ch_m^{C(m)}$  and a basis $\{
e_{\nu[m]} \}$ for $\ch_m'$ as in
Equation~(\ref{eq:representations}), such that for all $\gamma$ in
$G$,
\begin{equation}\label{eq:ginvbasis}
\rho(\gamma)e_{\mu[m]} = e_{\mu[\gamma^{-1} m\gamma]}.
\end{equation}

Assume that $\Gammah$ is a tree, then let $\MM_+(m'_+) :=
\MM_{g_+,n_++1}(\bmb_{N_+},m'_+)$ and $\MM_-(m'_-) \linebreak[0] :=
\MM_{g_-,n_-+1}(m'_-,\bmb_{N_-})$ for all $m'_\pm$ in $\mb_\pm$.
Let $\MM_\pm(\mb_\pm) :=
\coprod_{m'_\pm\in\mb_\pm}\MM_\pm(m'_\pm)$. We can write
$\MM_{\Gammatc} = \MM_+(m_+)\times \MM_-(m_-)$. Similarly, let
\[
\Lambda_+(v_{m_+}) := \Lambda_{g_+,n_++1}(v_{\bmb_{N_+}},v_{m_+})
\]
and
\[
\Lambda_-(v_{m_-}) := \Lambda_{g_-,n_-+1}(v_{m_-},v_{\bmb_{N_-}})
\]
for all $v_{m_\pm}$ in $\ch_{m_\pm}$.  Furthermore, let
$\Lambdah_\pm(v_{\mb_\pm})$ be defined by
\[
\Lambda_\pm(v_{\mb_\pm}) = \stt^*\Lambdah_\pm(v_{\mb_\pm})
\]
for all $v_{\mb_\pm}$ in $\ch_{\mb_\pm}$.

The \Gcft\ axioms imply that
\begin{eqnarray*}
I &=& \sum_{\alpha[m_\pm]} \stt''_* (\Lambda_+(e_{\alpha[m_+]})
\eta^{\alpha[m_+] \alpha[m_-]} \Lambda_-(e_{\alpha[m_-]})) \\ &=&
\sum_{\alpha[m_\pm]} \stt_{\mb_+ *}\Lambda_+(e_{\alpha[m_+]})
\eta^{\alpha[m_+] \alpha[m_-]} \stt_{\mb_-
*}\Lambda_-(e_{\alpha[m_-]}),
\end{eqnarray*}
where we have the natural forgetful morphisms $\stt_{\mb_+}:\MM_+(\mb_+) \rTo
\M_{g_+,n_++1}(\BG; \bmb_{N_+},\mb_+)$ and
$\stt_{\mb_-}:\MM_-(\mb_-) \rTo
\M_{g_-,n_-+1}(\BG; \mb_-,\bmb_{N_-})$.

Therefore, since $\Lambda$ is $G$-equivariant and the fibers of
$\stt_{\mb_\pm}$ are $G$-orbits, we have
\begin{eqnarray*}
I &=&\sum_{\alpha[m_\pm]} \stt_{\mb_+
*}\Lambda_+(\pi_G(e_{\alpha[m_+]})) \eta^{\alpha[m_+] \alpha[m_-]}
\stt_{\mb_- *}\Lambda_-(\pi_G (e_{\alpha[m_-]})) \\ &=&
\sum_{\mu[m_\pm]} \stt_{\mb_+ *}\Lambda_+(\pi_G(e_{\mu[m_+]}))
\eta^{\mu[m_+] \mu[m_-]} \stt_{\mb_- *}\Lambda_-(\pi_G
(e_{\mu[m_-]})) \\ &=& \sum_{\mu[m_\pm]} \stt_{\mb_+
*}\stt_{\mb_+}^* \Lambdah_+(\pi_G(e_{\mu[m_+]})) \eta^{\mu[m_+]
\mu[m_-]} \stt_{\mb_- *}\stt_{\mb_-}^* \Lambdah_-(\pi_G
(e_{\mu[m_-]})) \\ &=& \sum_{\mu[m_\pm]} \deg(\stt_{\mb_+})
\deg(\stt_{\mb_- }) \Lambdah_+(\pi_G(e_{\mu[m_+]})) \eta^{\mu[m_+]
\mu[m_-]} \Lambdah_-(\pi_G (e_{\mu[m_-]})) \\ &=&\sum_{\mu[m_\pm]}
\deg(\stt_{\mb_+}\times \stt_{\mb_- })
\Lambdah_+(\pi_G(e_{\mu[m_+]})) \eta^{\mu[m_+] \mu[m_-]}
\Lambdah_-(\pi_G (e_{\mu[m_-]})),
\end{eqnarray*}
where the first equality holds because $\stt_*\Lambda_{g,n}$ belongs to
$H^\bullet(\M_{g,n}(\BG))\otimes\chb^{*\otimes n}$, and the second follows
from the choice of basis and Proposition \ref{prop:gmodule}(\ref{prop:gmodule_a}).

Furthermore, let $\stt''_{(m'_+,m'_-)}$ denote the forgetful morphism
\[
\MM_+(m_+')\times\MM_-(m_-')\rTo \M_{g_+,n_++1} (\BG;\bmb_{N_+},
\mb_+)\times \M_{g_-,n_-+1}(\BG; \mb_-,\bmb_{N_-})
\]
 for all
$m'_\pm$ in $\mb_\pm$, then
\[
\deg(\stt_{\mb_+}\times\stt_{\mb_-}) = \sum_{m'_\pm\in \mb_\pm}
\deg(\stt''_{(m'_+,m'_-)}) = \frac{|G|^2}{|C(\mb_+)|^2}
\deg(\stt''_{(m_+,m_-)}),
\]
where in the second equality, we have used that $\deg(\stt''_{(m'_+,m'_-)})$
is independent of the choice $m'_\pm$ in $\mb_\pm$, the fact that $\mb'_\pm$
contains $\frac{|G|}{|C(\mb_\pm)|}$ elements, and that $|C(\mb)| =
|C(\mb^{-1})| $ for all conjugacy classes $\mb$ in $G$.
Thus,
\begin{equation}
I = \sum_{\mu[m_\pm]} \deg(\stt'') \frac{|G|^2}{|C(\mb_+)|^2}
\Lambdah_+(\pi_G(e_{\mu[m_+]})) \eta^{\mu[m_+] \mu[m_-]}
\Lambdah_-(\pi_G (e_{\mu[m_-]})),
\end{equation}
but $\stt'' = \rh\circ\prt\circ\rt$, hence,
\begin{equation}\label{eq:I_two}
I = \sum_{\mu[m_\pm]} \deg(\rh) \deg(\prt\circ\rt) \frac{|G|^2}{|C(\mb_+)|^2}
\Lambdah_+(\pi_G(e_{\mu[m_+]})) \eta^{\mu[m_+] \mu[m_-]}
\Lambdah_-(\pi_G (e_{\mu[m_-]})).
\end{equation}

Equating Equations~(\ref{eq:I_one}) and (\ref{eq:I_two}) and canceling
factors of $\deg(\prt\circ\rt)$, we obtain
\begin{equation}
\rh_*\muh^*\ih^*\Lambdah_{g,n}(v_\bmb) = \frac{|G|^2}{|C(\mb_+)|^2}
\deg(\rh) \sum_{\mu[\mb_\pm]}
\Lambdah_+(\pi_G(e_{\mu[m_+]}))
\eta^{\mu[\mb_+] \mu[\mb_-]}
\Lambdah_-(\pi_G(e_{\mu[m_-]})).
\end{equation}

Let $\epsilon_{\mu[\mb_\pm]} := \pi_G(e_{\mu[m_\pm]})$. Notice
that the left hand side only depends upon $\mb_\pm$ because of
Equation (\ref{eq:ginvbasis}). Since $\{ e_{\mu[m_\pm]} \}$ is a
basis for $\ch_{m_\pm}^{C(m_\pm)}$, then by Proposition
\ref{prop:gmodule}(\ref{prop:gmodule_b}), $\{
\epsilon_{\mu[\mb_\pm]} \}$ is a basis for $\chb_{\mb_\pm}$. Let
\[
\etab_{\mu[\mb_+]\mu[\mb_-]} :=
\etab(\epsilon_{\mu[\mb_+]},\epsilon_{\mu[\mb_-]})
= \frac{1}{|G|} \eta(\epsilon_{\mu[\mb_+]},\epsilon_{\mu[\mb_-]})
= \frac{|C(m_+)|}{|G|^2} \eta_{\mu[m_+] \mu[m_-]}
\]
where $\eta_{\mu[m_+] \mu[m_-]} =
\eta(e_{\mu[m_+]},e_{\mu[m_-]})$. Therefore, taking inverses,
\begin{equation}\label{eq:etab_basis}
\etab^{\mu[\mb_+]\mu[\mb_-]} = \frac{|G|^2}{|C(\mb_+)|}
\eta^{\mu[m_+]\mu[m_-]}
\end{equation}
and
\[
\etah^{\mu[\mb_+]\mu[\mb_-]} = \frac{|G|^2}{|C(\mb_+)|^2}\eta^{\mu[m_+]\mu[m_-]}
\]
by Equation (\ref{eq:etah}), so
\begin{equation} \label{eq:cutting_tree_almost}
\rh_*\muh^*\ih^*\Lambdah_{g,n}(v_\bmb) = \deg(\rh)
\sum_{\mu[\mb_\pm]}
\Lambdah_{g_+,n_++1}(v_{\bmb_{N_+}},\epsilon_{\mu[\mb_+]}))
\etah^{\mu[\mb_+] \mu[\mb_-]}
\Lambdah_{g_-,n_-+1}(\epsilon_{\mu[\mb_-]},v_{\bmb_{N_-}}).
\end{equation}
To conclude, note that $\prh = \sth''\circ\rh$ and $\deg(\prh) =
\deg(\sth')$, so
\[
\deg(\rh) = \frac{\deg(\prh)}{\deg(\sth'')} =
\frac{\deg(\sth')}{\deg(\sth'')}.
\]
This finishes the tree case.

Suppose now that $\Gammah$ is a loop and that $\MM_{\Gammatc} =
\MM_{g-1,n+2}(\bmb,m_+,m_-)$. Following the analogous steps to the
case of the tree, we obtain the counterpart of Equation~(\ref{eq:I_two}):
\begin{equation}
I = \sum_{\mu[m_\pm]} \deg(\rh) \deg(\prt\circ\rt) \frac{|G|^2}{|C(\mb_+)|^2}
\Lambdah_{g-1,n+2}(v_{\bmb},\pi_G(e_{\mu[m_+]}),\pi_G(e_{\mu[m_-]}))
\eta^{\mu[m_+] \mu[m_-]}.
\end{equation}
Proceeding further, the counterpart of Equation~(\ref{eq:cutting_tree_almost}) is
\begin{equation}
\rh_*\muh^*\ih^*\Lambdah_{g,n}(v_\bmb) = \deg(\rh)
\sum_{\mu[\mb_\pm]}
\Lambdah_{g-1,n+2}(v_{\bmb},\epsilon_{\mu[\mb_+]}),\epsilon_{\mu[\mb_-]})
\etah^{\mu[\mb_+] \mu[\mb_-]}.
\end{equation}
The rest of the proof is essentially the same as in the case of the tree.
\end{proof}

\subsection{From $\M_{g,n}(\BG)$ to $\M_{g,n}$}

\ 

\begin{df}  Let $(\ch,\eta,\Lambda_{g,n},\vac)$ be a \Gcft. Define
$\Lambdab_{g,n} := \sth_*\Lambdah_{g,n}$ in
$H^\bullet(\M_{g,n})\otimes\ch^{*\otimes n}$.
\end{df}

\begin{thm}\label{thm:Coinvariants}
If $(\ch,\eta,\Lambda_{g,n},\vac)$ is a \Gcft,\ then
$(\chb,\etab,\Lambdab_{g,n},\vac)$ forms a \cft.
\end{thm}
\begin{proof}
We begin by observing that
\begin{equation}\label{eq:rightsquare}
\frac{1}{\deg{\sth}}i^*\sth_* = \frac{1}{\deg{\stt'}}\sth'_*
\ih^*,
\end{equation}
since the lower right square is not Cartesian, due to ramification over
$\M_\Gamma$.

Next, we observe that
\begin{equation}\label{eq:leftsquare}
\mu^*\sth'_* = \prh_*\muh^*,
\end{equation}
since the lower left square is Cartesian by definition.

Therefore, for all $v_\bmb$ in $\chb_\bmb$,
\begin{eqnarray*}
\varrho_\Gamma^*\Lambdab_{g,n}(v_\bmb) & = &
\varrho_\Gamma^* \sth_* \Lambdah_{g,n}(v_\bmb) \\
&=& \mu^* i^*\sth_*\Lambdah_{g,n}(v_\bmb) \\
&=& \mu^* \left(\frac{\deg(\sth)}{\deg(\sth')} \sth'_*\ih^*\right)\Lambdah_{g,n}(v_\bmb)
\\
&=& \frac{\deg(\sth)}{\deg(\sth')} \prh_*\muh^*\ih^*\Lambdah_{g,n}(v_\bmb)
\\
&=& \frac{\deg(\sth)}{\deg(\sth')} \sth''_*(\rh_* \muh^*\ih^*)
\Lambdah_{g,n}(v_\bmb)
\\
&=& \frac{\deg(\sth)}{\deg(\sth')} \sth''_*
\left(\frac{\deg(\sth')}{\deg(\sth'')}
\Lambdah_{\Gammahc}\right)\\ &=& \frac{\deg(\sth)}{\deg(\sth'')}
\sth''_*\Lambdah_{\Gammahc},
\end{eqnarray*}
where Equations (\ref{eq:cutting_tree}) and (\ref{eq:cutting_loop}) have
been used in the sixth equality.

Assume that $\Gammah$ is a tree.  Adopting the notation from the
proof of Theorem~\ref{thm:cutting_edges}, we obtain
\begin{eqnarray*}
\varrho_\Gamma^*\Lambdab_{g,n}(v_\bmb) &=& \frac{\deg(\sth)}{\deg(\sth'')}
            \sum_{\beta[\mb_+],\beta[\mb_-]}
\sth''_*(\Lambdah_+(e_{\beta[\mb_+]})
\etah^{\beta[\mb_+] \beta[\mb_-]}
\Lambdah_-(e_{\beta[\mb_-]})) \\
&=& \frac{\deg(\sth)}{\deg(\sth'')}
            \sum_{\beta[\mb_+],\beta[\mb_-]}
\Lambdab_+(e_{\beta[\mb_+]}) \etah^{\beta[\mb_+] \beta[\mb_-]}
\Lambdab_-(e_{\beta[\mb_-]}),
\end{eqnarray*}
where $\Lambdab_\pm := \st_{\mb_\pm *} \Lambdah_\pm$.  This can be rewritten as
\begin{equation}\label{eq:cross_multiplied}
\deg(\sth'') \varrho_\Gamma^*\Lambdab_{g,n}(v_\bmb) = \deg(\sth)
\sum_{\beta[\mb_+],\beta[\mb_-]} \Lambdab_+(e_{\beta[\mb_+]})
\etah^{\beta[\mb_+] \beta[\mb_-]} \Lambdab_-(e_{\beta[\mb_-]}).
\end{equation}
Following \cite{JK}, let
\begin{equation}\label{eq:Omega}
\Omega_{g,n}(\bmb) := \deg(\sth)
\end{equation}
for all $\bmb$ in $\Gb^n$.  We have
\[ \deg(\sth'') = \Omega_{g_+,n_++1}(\mb_{N_+},\mb_+)
\Omega_{g_-,n_-+1}(\mb_-,\mb_{N_-}),
\]
and Equation~(\ref{eq:cross_multiplied}) becomes, after multiplying both
sides by $|C(\mb_+)|$, using the definition of $\etah$, and summing over all
conjugacy classes $\mb_\pm$ such that $\mb_- = \mb_+^{-1}$,
\begin{eqnarray*}
\sum_{\mb_\pm:\mb_- = \mb_+^{-1}} &|C(\mb_+)|&
\Omega_{g_+,n_++1}(\mb_{N_+},\mb_+)
\Omega_{g_-,n_-+1}(\mb_-,\mb_{N_-})
\varrho_\Gamma^*\Lambdab_{g,n}(v_\bmb) = \\ & \Omega_{g,n}(\bmb)&
\sum_{\mb_\pm:\mb_- = \mb_+^{-1}} \sum_{\beta[\mb_+],\beta[\mb_-]}
\Lambdab_+(e_{\beta[\mb_+]}) \etab^{\beta[\mb_+] \beta[\mb_-]}
\Lambdab_-(e_{\beta[\mb_-]}).
\end{eqnarray*}
But Lemma 3.5(1) from \cite{JK} states that
\[
\sum_{\mb_\pm:\mb_- = \mb_+^{-1}} |C(\mb_+)| \Omega_{g_+,n_++1}(\mb_{N_+},\mb_+)
\Omega_{g_-,n_-+1}(\mb_-,\mb_{N_-}) = \Omega_{g,n}(\bmb).
\]
Therefore, by canceling $\Omega_{g,n}(\bmb)$ from both sides, we obtain
the desired result.

In the case of the loop, we have
\begin{eqnarray*}
\varrho_\Gamma^*\Lambdab_{g,n}(v_\bmb) &=& \frac{\deg(\sth)}{\deg(\sth'')}
            \sum_{\beta[\mb_+],\beta[\mb_-]}
\sth''_*(\Lambdah_{g-1,n+2}(v_\bmb, e_{\beta[\mb_+]},e_{\beta[\mb_-]}))
\etah^{\beta[\mb_+] \beta[\mb_-]} \\
&=& \frac{\deg(\sth)}{\deg(\sth'')}
            \sum_{\beta[\mb_+],\beta[\mb_-]}
\Lambdab_{g-1,n+2}(v_\bmb, e_{\beta[\mb_+]},e_{\beta[\mb_-]})
\etah^{\beta[\mb_+] \beta[\mb_-]}.
\end{eqnarray*}
Multiplying both sides by $\deg(\sth'') |C(\mb_+)|$, plugging in
$\deg(\sth'') = \Omega_{g-1,n+2}(\bmb,\mb_+,\mb_-)$, $\deg(\sth) =
\Omega_{g,n}(\bmb)$, and then summing over all conjugacy classes $\mb_\pm$
such that $\mb_- = \mb_+^{-1}$, we obtain
\begin{eqnarray*}
\sum_{\mb_\pm:\mb_- = \mb_+^{-1}} &|C(\mb_+)|& \Omega_{g-1,n+2}(\bmb,\mb_+,\mb_-)
\varrho_\Gamma^*\Lambdab_{g,n}(v_\bmb) = \\
&\Omega_{g,n}(\bmb)& \sum_{\mb_\pm:\mb_- = \mb_+^{-1}}
\sum_{\beta[\mb_+],\beta[\mb_-]}
\Lambdab_{g-1,n+2}(v_\bmb,e_{\beta[\mb_+]},e_{\beta[\mb_-]})
\etab^{\beta[\mb_+] \beta[\mb_-]}.
\end{eqnarray*}
Since Lemma 3.5(2) from \cite{JK} states that
\[
\sum_{\mb_\pm:\mb_- = \mb_+^{-1}} |C(\mb_+)|
\Omega_{g-1,n+2}(\bmb,\mb_+,\mb_-) = \Omega_{g,n}(\bmb),
\]
we may cancel $\Omega_{g,n}(\bmb)$ from both sides to obtain the desired
result.

This completes the proof of the factorization axiom of the \cft.

The invariance under the symmetric group is manifest.

The flat identity axiom follows from considering the following commuting
diagram:
\begin{diagram}
\MM_{g,n+1}(\bmb,1) & \rTo^{\tst} & \MM_{g,n}(\bmb) \\
\dTo_{\stt_1} & & \dTo_\stt\\
\M_{g,n+1}(\BG;\bmb,1)& \rTo^\tsh & \M_{g,n}(\BG;\bmb) \\
\dTo_{\sth_1} & & \dTo_\sth\\
\M_{g,n+1}& \rTo^\ts & \M_{g,n}\\
\end{diagram}
The horizontal morphisms are forgetting-tails morphisms and are both
flat and proper.  The vertical morphisms are forgetful morphisms and
are all quasi-finite, flat, and proper.

By Lemma~\ref{lm:useful} We have 
\begin{eqnarray*}
\stt_1^*\Lambdah_{g,n+1}(v_{\bmb},\vac) 
& = &\Lambda_{g,n+1}(v_\bmb,\vac) \\
& = &\tst^*\Lambda_{g,n}(v_\bmb)\\
& = &\tst^*\stt^*\Lambdah_{g,n}(v_\bmb)\\
& = &\stt_1^*\tsh^*\Lambdah_{g,n}(v_\bmb)
\end{eqnarray*}
By the uniqueness of the classes $\Lambdah$ (again, see Lemma~\ref{lm:useful})
we conclude that
\begin{equation}
\tsh^*\Lambdah_{g,n}(v_\bmb) = \Lambdah_{g,n+1}(v_\bmb,\vac).
\end{equation}

On the other hand, while the bottom square of this diagram is not
Cartesian, it is almost so---the stack $\M_{g,n+1}(\BG;\bmb,1)$
is the universal \emph{orbicurve} over $\M_{g,n}(\BG;\bmb)$, and
it is birational to its coarse moduli space, the universal curve
over $\M_{g,n}(\BG;\bmb)$.  Thus, we have
\[
\ts^*\Lambdab_{g,n}(v_\bmb) = \ts^*\sth_*\Lambdah_{g,n}(v_\bmb) =
\st_{1 *}\tsh^*\Lambdah_{g,n}(v_\bmb) =
\sth_{1 *}\Lambdah_{g,n+1}(v_\bmb,\vac) = \Lambdab_{g,n+1}(v_\bmb,\vac).
\]

The last property that must be verified is
\begin{equation}\label{eq:etab}
\etab(v_{\mb_+},v_{\mb_-}) = \Lambdab_{0,3}(v_{\mb_+},v_{\mb_-},\vac)
\end{equation}
for all $v_{\mb_\pm}$ in $\chb_\pm$, where we have identified
$H^\bullet(\M_{0,3})$ with the ground ring $\nc$. Since this
identity holds trivially if $\mb_- \not= \mb_+^{-1}$, let us
assume that $\mb_-=\mb_+^{-1}$.

We have the morphisms
\[
\coprod_{m'_+\in\bmb_+} \xi_{0,3}(m'_+,m_+^{'-1},1)\rTo^{\stt_\xi}
\M_{0,3}(\BG; \mb_+, \mb_+^{-1}, 1) \rTo^\sth \M_{0,3},
\]
and we let $\st_\xi := \sth\circ \stt_\xi$. Since $\eta$ is defined by
\[
\eta(v_{m_+},v_{m_-}) = \int_{\db{\xi(m_+,m_-,1)}} \Lambda_{0,3}(v_{m_+},v_{m_-},\vac)
\st_{\xi *} \Lambda_{0,3}(v_{m_+},v_{m_-},\vac),
\]
we have
\begin{eqnarray*}
\Lambdab_{0,3}(v_{\mb_+},v_{\mb_-},\vac) &= &
\sth_*\Lambdah_{0,3}(v_{\mb_+},v_{\mb_-},\vac) \\ &=& \sth_*
(\stt_{\xi *} \stt_\xi^* \frac{1}{\deg{\stt_\xi}})
\Lambdah_{0,3}(v_{\mb_+},v_{\mb_-},\vac) \\ &=&
\frac{1}{\deg{\stt_\xi}} \sth_* \stt_{\xi *}
\Lambda_{0,3}(v_{\mb_+},v_{\mb_-},\vac) \\ &=&
\frac{1}{\deg{\stt_\xi}} \st_{\xi
*}\Lambda_{0,3}(v_{\mb_+},v_{\mb_-},\vac) \\ &=&
\frac{1}{\deg{\stt_\xi}} \eta(v_{\mb_+},v_{\mb_-}),
\end{eqnarray*}
but
\[
\deg(\stt_\xi) = |\mb_+| |C(\mb_+)| = |G|,
\]
since a generic point of $\M_{0,3}(\BG;\mb_+,\mb_-,1)$ has automorphism
group isomorphic to $C(m_+)$. Therefore, Equation~(\ref{eq:etab}) is satisfied.
\end{proof}

\begin{rem}
The \cft\ $(\chb,\etab,\{ \Lambdab_{g,n}\},\vac)$ constructed
above has more structure than a generic \cft, as it is
$\Gb$-graded; that is, $(\chb,\etab)$ is a $\Gb$-graded vector
space with metric, and for all $v_\bmb$ in $\chb_\bmb$, the class
$\Lambdab_{g,n}(v_\bmb)$ vanishes unless there exist
representatives $m'_i$ in $\mb_i$ for all $i=1,\ldots,n$ such that
$\prod_{i=1}^n m'_i$ belongs to the subgroup $[G,G]^g$. This is
follows from the fact that $\M_{g,n}(\BG;\bmb)$ is empty unless
this holonomy condition holds.
\end{rem}

\begin{prop}
Let $(\ch,\eta,\{ \Lambda_{g,n}\},\vac)$ be a \Gcft. For all nonzero $\lambda$
in $\nc$, $(\ch,\lambda^{-2} \eta ,\{ \lambda^{2g-2}\Lambda_{g,n}\},\vac)$ is
a \Gcft.
\end{prop}
The proof is immediate from the definition.

\begin{rem}\label{rem:StringCC}
One can eliminate the annoying factor of $\frac{1}{|G|}$ in the
definition of $\etab$ by choosing $\lambda$ such that $\lambda^2 =
\frac{1}{|G|}$.  In this case, the associated ``quotient'' by $G$
of the \Gcft\ $(\ch,|G| \eta ,\{ |G|^{1-g} \Lambda_{g,n}\},\vac)$
is the \cft\ $(\chb,|G| \etab,\{ |G|^{1-g}\Lambda_{g,n}\},\vac)$,
but $|G| \etab$ is equal to the restriction of $\eta$ to $\chb$.
\end{rem}

\subsection{Coinvariants of $G$-Frobenius algebras}

\ 

Let $((\ch,\rho),\cdot,\vac,\eta)$ be a $G$-Frobenius algebra. We
now have two ways to endow its space of coinvariants $\chb$ with
the structure of a Frobenius algebra. The first is purely
algebraic. The tuple $(\chb,\cdot,\vac,\eta)$ is a Frobenius
algebra where the multiplication on $\chb$ is inherited by
restriction from the multiplication on $\ch$ and the metric $\eta$
is the restriction of the metric on $\ch$.

The second is to apply the geometric procedure described in the
previous subsection to $\ch$, regarded as a \Gcft, to induce the
structure of a Frobenius algebra on $\chb$. It turns out that
these two Frobenius structures are identical after a rescaling.

In order to simplify the proof, we note that the structure of the
$G$-Frobenius algebra $((\ch,\rho),\cdot,\vac,\eta)$ can also be
described as the tuple $((\ch,\rho),\mu,\vac),$ where $\mu$
belongs to $\ch^{* \otimes 3}$ and is defined by
\begin{equation}\label{eq:mu}
\mu(v_{m_1},v_{m_2},v_{m_3}) := \eta(v_{m_1}, v_{m_2} \cdot v_{m_3}),
\end{equation}
since it follows that $\eta(v_{m_1},v_{m_2}) = \eta(v_{m_1}, v_{m_2} \cdot
\vac)$. If $\mut$ denotes the restriction of $\mu$ to $\chb$, then the data
$(\chb,\mut,\vac)$ is an equivalent description of the Frobenius algebra
structure on $\chb$ induced by restriction.

\begin{prop}\label{prop:agreement}
Let $((\ch,\rho),\mu,\vac)$ be a $G$-Frobenius algebra arising from a
\Gcft\ $(\ch,\eta,\{ \Lambda_{g,n} \}, \vac)$. The Frobenius algebra
structure on $\chb$ arising from the \cft\ $(\ch,\etab,\{
\Lambdab_{g,n} \}, \vac)$ is $(\chb,\mub,\vac)$, where
\[
\mub = \frac{1}{|G|} \mut,
\]
and $\mut$ is the restriction of $\mu$ to $\chb$.
\end{prop}

\begin{rem}
The Frobenius algebra $(\chb,\mub,\vac)$ can also be described as the
tuple $(\chb,\cdot,\etab,\vac),$ where the multiplication $\cdot$ on
$\chb$ is inherited from the multiplication on $\ch$, but where
$\etab$ is the restriction of $\frac{1}{|G|} \eta$ to $\chb$.
\end{rem}

\begin{proof} (of Proposition \ref{prop:agreement})

Since $\mub(v_\bmb) =
\Lambdab_{0,3}(v_\bmb)$, after identifying $H^\bullet(\M_{0,3})$ with $\nc$,
we need only prove that
\[
\Lambdab_{0,3}(v_\bmb) = \frac{1}{|G|} \mu(v_\bmb)
\]
for all $v_\bmb$ in $\chb_\bmb$. In order to proceed, let us introduce some
notation.

For all $\bm := (m_1,m_2,m_3)$ belonging to $G^3$ such that $m_1 m_2
m_3 = 1$, we have the following forgetful morphisms
\[
\MM_{0,3}(\bmb)\rTo^\stt \M_{0,3}(\BG;\bmb)\rTo^\sth \M_{0,3},
\]
and we let $\st :=\sth\circ\stt$.

Furthermore, if $Q$ is a substack of $\MM_{0,3}(\bmb)$ then we let
$\stt_Q$ denote the restriction of $\stt$ to $Q$. Let
\[
\xi := \coprod_{\bm'\in\chi(\bmb)} \xi(\bm'),
\]
where
\[
\chi(\bmb) := \{ (m'_1,m'_2,m'_3)\in \bmb |  m'_1 m'_2 m'_3 = 1 \}.
\]
Henceforth, fix an element $\bm$ in $\chi(\bmb)$ once and for all.

Let us adopt the notation that for any $v_\bmb$ in $\chb_\bmb$,
and for any $\bm'\in\bmb$, the vector $v_{\bm'}$ 
denotes the $\bm'$-graded component of $v_{\bmb}$; that is,
\[
v_\bmb =: \sum_{\bm'\in\bmb} v_{\bm'}.
\]
Note that $v_{\bm'}$ belongs to the subspace of $C(\bm')$-invariant
vectors in $\ch_{\bm'}$, where $C(\bm') := C(m'_1)\times C(m'_2)\times
C(m'_3)$.

For all $\bm'$ in $\chi(\bmb)$, we have
\[
\mu(v_\bm') = \frac{1}{\deg(\st_{\xi(\bm')})} \st_{\xi(\bm')
*}\Lambda_{0,3}(v_\bm').
\]
Otherwise, $\mu(v_\bm') = 0$. Since $\mut$ is the restriction of $\mu$ to
$\chb$,
\begin{eqnarray*}
\mut(v_\bmb) &=& \sum_{\bm'\in\chi(\bmb)}
\frac{1}{\deg(\st_{\xi(\bm')})} \st_{\xi(\bm')
*}\Lambda_{0,3}(v_{\bm'}) \\ &=& \frac{1}{\deg(\st_{\xi(\bm)})}
\sum_{\bm'\in\chi(\bmb)} \st_{\xi(\bm') *}\Lambda_{0,3}(v_{\bm'})
\\ &=& \frac{1}{\deg(\st_{\xi(\bm)})} \st_{\xi *}
\Lambda_{0,3}(v_\bmb),
\end{eqnarray*}
where the second equality comes from the fact that the degree of $\st$
restricted to any connected component of $\MM_{0,3}(\bmb)$ is independent of
the choice of connected component. This statement follows from the fact that
every connected component of $\MM_{0,3}(\bmb)$ is $\rho(\bg)\xi(\bm)$ for
some $\bg$ in $G^3$, but $\rho(\bg)$ is an isomorphism.

However, we have
\begin{eqnarray*}
\st_{\xi *} \Lambda_{0,3}(v_\bmb) &=& \sum_{\bm'\in\chi(\bmb)}
\st_{\xi(\bm') *} \Lambda_{0,3}(v_{\bm'}) \\ &=& |\chi(\bmb)|
\st_{\xi(\bm) *} \Lambda_{0,3}(v_\bm) \\ &=& |G|
\Omega_{0,3}(\bmb) \st_{\xi(\bm) *} \Lambda_{0,3}(v_\bm),
\end{eqnarray*}
where the second equality follows from the observation that every
connected component of $\xi$ can be obtained by the action of some
element of $G^3$, and the fact that $v_\bmb$ and $\Lambda_{0,3}$
are $G^3$-invariant. The third equality is from Proposition 3.4 of
\cite{JK}, where $\Omega_{0,3}$ is defined in Equation
(\ref{eq:Omega}). Therefore, we obtain
\begin{equation}
\mut(v_\bmb) = \frac{|G| \Omega_{0,3}(\bmb)}{\deg(\st_{\xi(\bm)})}
\st_{\xi(\bm) *} \Lambda_{0,3}(v_\bm).
\end{equation}
On the other hand, the definition of $\Lambdah_{0,3}$ implies that
\[
\Lambdah_{0,3}(v_\bmb) = \frac{1}{\deg(\stt)}\stt_*\Lambda_{0,3}(v_\bmb),
\]
hence
\begin{eqnarray*}
\Lambdab_{0,3}(v_\bmb) &=& \sth_*\Lambdah_{0,3}(v_\bmb) \\ &=&
\frac{1}{\deg(\stt)} \sth_* \stt_* \Lambda_{0,3}(v_\bmb),
\end{eqnarray*}
and we obtain
\begin{equation}\label{eq:total_b}
\Lambdab_{0,3}(v_\bmb) = \frac{1}{\deg(\stt)} \st_*
\Lambda_{0,3}(v_\bmb).
\end{equation}
Using the fact that $\deg(\st_Q)$ is independent of the choice of
connected component $Q$ of $\MM_{0,3}(\bmb)$, we can write
\begin{equation} \label{eq:denominator_b}
\deg(\stt) = A(\bmb) \deg(\stt_{\xi(\bm)}),
\end{equation}
where $A(\bmb)$ is the number of connected components of
$\MM_{0,3}(\bmb)$. Similarly, let $I(\bmb)$ consist of all
elements $\bg$ in $G^3$ such that the collection $\{
\rho(\bg)\xi(\bm) \}$ is in one-to-one correspondence with the
connected components of $\MM_{0,3}(\bmb)$, then
\begin{eqnarray}
\st_*\Lambda_{0,3}(v_\bmb) &=& \sum_{\bg\in I(\bmb)}
\st_{\rho(\bgamma)\xi(\bm) *} \Lambda_{0,3}(v_{\rho(\bg)\bm})
\nonumber\\ &=& \sum_{\bg\in I(\bmb)} \st_{\xi(\bm)*}
\rho(\bgamma^{-1})_*\Lambda_{0,3}(v_{\rho(\bg)\bm})\nonumber \\
&=&\sum_{\bg\in I(\bmb)} \st_{\xi(\bm)*}
\rho(\bgamma)^*\Lambda_{0,3}(v_{\rho(\bg)\bm}) \nonumber\\
&=&\sum_{\bg\in I(\bmb)} \st_{\xi(\bm)*}
\Lambda_{0,3}(\rho(\bgamma^{-1})v_{\rho(\bg)\bm}) \nonumber\\
&=&\sum_{\bg\in I(\bmb)} \st_{\xi(\bm)*} \Lambda_{0,3}(v_{\bm})
\nonumber\\ &=& A(\bmb) \st_{\xi(\bm) *}\Lambda_{0,3}(v_\bm)
\label{eq:temp},
\end{eqnarray}
where the first equality is the sum over contributions from each
connected components of $\MM_{0,3}(\bmb)$, and the second is from
the fact that, for all $\bg$ in $G^3$, we have
\begin{equation}\label{eq:numerator_b}
\st_{\xi(\bm)} = \st_{\rho(\bg)\xi(\bm)}\circ\rho(\bg).
\end{equation}
The fourth equality is from the $G^3$-invariance of $\Lambda_{0,3}$ and the fifth is
from the $G^3$-invariance of $v_\bmb$.
Putting together Equations (\ref{eq:total_b}),  (\ref{eq:denominator_b}),
and (\ref{eq:temp}), we obtain
\begin{eqnarray*}
\Lambdab_{0,3}(v_\bmb) &=& \frac{1}{\deg(\stt_{\xi(\bm)})}
\st_{\xi(\bm) *}\Lambda_{0,3}(v_\bm) \\ &=&
\frac{\Omega_{0,3}(\bmb)}{\deg(\st_{\xi(\bm)})} \st_{\xi(\bm)
*}\Lambda_{0,3}(v_\bm),
\end{eqnarray*}
since
\[
\deg(\st_{\xi(\bm)}) = \deg(\stt_{\xi(\bm)}) \deg(\sth) =
\deg(\stt_{\xi(\bm)})
\Omega_{0,3}(\bmb).
\]
However,
\begin{eqnarray*}
\st_*\Lambda_{0,3}(v_\bmb) &=& \sum_{\bm'\in\chi(\bmb)}
\st_{\xi(\bm') *}\Lambda_{0,3}(v_{\bm'}) \\ &=& |\chi(\bmb)|
\st_{\xi(\bm) *}\Lambda_{0,3}(v_{\bm}) \\ &=&
|G|\Omega_{0,3}(\bmb) \st_{\xi(\bm) *}\Lambda_{0,3}(v_{\bm}).
\end{eqnarray*}
Putting this all together, we obtain the desired result
\begin{eqnarray*}
\Lambdab_{0,3}(v_\bmb)
&=& \frac{\Omega_{0,3}(\bmb)}{\deg(\st_{\xi(\bm)})} \frac{1}{|G|
\Omega_{0,3}(\bmb)} \st_{\xi *}\Lambda_{0,3}(v_\bmb) \\
&=&\frac{1}{|G|} \mut(v_\bmb).
\end{eqnarray*}
\end{proof}

The results of this section can be applied to the example of the group
ring \Gcft\ and its associated $G$-Frobenius algebra to yield the
(stringy) orbifold cohomology of a point with trivial $G$-action.

\begin{prop}
The Frobenius algebra $\chb$ induced from the $G$-Frobenius
algebra $\ch = \nc[G]$ is the Frobenius algebra $Z(\nc[G])$, the
center of the group ring, with its induced multiplication,
identity, and the metric $\etab$.

The resulting Frobenius algebra is isomorphic to the orbifold (stringy)
quantum cohomology of $\BG$, the classifying stack of $G$.
\end{prop}
We refer the reader to \cite{JK} where the calculation is worked
through in detail.

\subsection{The quotient stack $[\MM_{g,n}/G^n]$}

\ 

The process of obtaining a \cft\ from a \Gcft\ involved the stack
$\M_{g,n}(\BG)$. However, there is another stack that one could
have used instead, namely, the quotient stack $\MMq_{g,n} :=
[\MM_{g,n}/G^n]$ and its substacks    $\MMq_{g,n}(\bmb) :=
[\MM_{g,n}(\bmb)/G^n]$. We will show that one can construct a
\cft\ by replacing $\M_{g,n}(\BG)$ by $\MMq_{g,n}(\bmb) :=
[\MM_{g,n}(\bmb)/G^n]$, but that the resulting \cft\ is isomorphic
to the original one.

We have the following
sequence of forgetful morphisms
\begin{equation}
\MM_{g,n}(\bmb)\rTo^{\stcp}\MMq_{g,n}(\bmb)\rTo^{\stc}\M_{g,n}(\BG;\bmb)
\rTo^\sth\M_{g,n},
\end{equation}
where $\stt := \stc\circ \stcp$.
The stack $\MMq$ is a smooth, Deligne-Mumford stack, and all of these
morphisms are proper and flat. Observe that while the morphism $\stc$ induces
an isomorphism at the level of the corresponding coarse moduli spaces, they
are not isomorphic as stacks, since an object in $\MMq_{g,n}(\bmb)$ has a larger
automorphism group than the corresponding object in $\M_{g,n}(\BG;\bmb)$.

\begin{df}
Let $((\ch,\rho),\eta,\{\Lambda_{g,n}\},\vac)$ be a \Gcft. Define the
elements $\Lambdac_{g,n}$ in $\bigoplus_{\bmb\in\Gb^n}\linebreak[0]
H^\bullet(\MMq_{g,n}(\bmb))\otimes\chb_\bmb^*$ via
\begin{equation}\label{eq:lambdac}
\Lambdac_{g,n}(v_\bmb):=\stc^* \Lambdah_{g,n}(v_{\bmb})
\end{equation}
for all $v_\bmb$ in $\chb_\bmb$ and $\bmb$ in $\Gb^n$. Define
$\Lambdab'_{g,n}$ in $H^\bullet(\M_{g,n})\otimes\chb_\bmb^*$ via
\[
\Lambdab'_{g,n}(v_\bmb) := (\sth\circ\stc)_* \Lambdac_{g,n}(v_\bmb).
\]
Let
\[
\etab'(v_{\mb_+},v_{\mb_-}) := \Lambdab'_{0,3}(v_{\mb_+},v_{\mb_-},\vac)
\]
for all $v_{\mb_\pm}$ in $\chb_{\mb_\pm}$.
\end{df}

\begin{prop}
Let $((\ch,\rho),\eta,\{\Lambda_{g,n}\},\vac)$ be a \Gcft.
\begin{enumerate}
\item  We have the identity \[
\Lambda_{g,n}(v_\bmb) = {\stc^{'*}} \Lambdac_{g,n}(v_\bmb).
\]
\item We also have $$\Lambdab'_{g,n} (v_\bmb):= \left(\prod_{i=1}^n
\frac{1}{k_{\mb_i}}\right) \Lambda_{g,n}(v_\bmb),$$ where $k_\mb$
is the order of the cyclic subgroup generated by any
representative of $\mb$ in $\Gb$.
\item $(\chb,\etab',\{\Lambdab'_{g,n}\},\vac)$ is a \cft.
\item The linear map $\phi:\chb\rTo\chb$, where
\[
\phi(v_\mb) := k_\mb v_\mb
\]
for all $v_\mb$ in $\chb_\mb$ and $\mb$ in $\Gb$, is an
isomorphism between the \cfts\
$(\chb,\etab,\{\Lambdab_{g,n}\},\vac)$ and
$(\chb,\etab',\{\Lambdab'_{g,n}\},\vac)$.

\end{enumerate}\end{prop}
\begin{proof}
Since $\stc^{' *}\Lambdac_{g,n}(v_\bmb) = \stt^*\Lambdah_{g,n}(v_\bmb) = \stc^{'*} \stc^*
\Lambdab_{g,n}(v_\bmb)$, we obtain
\[
\Lambda_{g,n}(v_{\bmb})=\stc'^* \Lambdac_{g,n}(v_\bmb).
\]
For the second part, apply $\stc_*$ to
both sides of equation (\ref{eq:lambdac}) and use the fact that
\[
\deg(\stc) = \prod_{i=1}^n \frac{1}{k_{\mb_i}}
\]
to get
\[
\stc_*\Lambdac_{g,n}(v_\bmb) = \left(\prod_{i=1}^n \frac{1}{k_{\mb_i}}\right)
\Lambdah_{g,n}(v_\bmb).
\]
Thus,
\[
\Lambdah_{g,n}(v_\bmb) = \left(\prod_{i=1}^n k_{\mb_i}\right) \stc_* \Lambdac_{g,n}(v_\bmb).
\]
However,
\begin{eqnarray*}
\Lambdab_{g,n}(v_\bmb) &=& \sth_*\Lambdah_{g,n}(v_\bmb) \\
&=& \left(\prod_{i=1}^n k_{\mb_i}\right) \sth_*\stc_* \Lambdac_{g,n}(v_\bmb)
\\
&=& \left(\prod_{i=1}^n k_{\mb_i}\right) \Lambdab'_{g,n}(v_\bmb).
\end{eqnarray*}
This establishes the second part of the proposition.

Clearly, $\phi^*\Lambdab'_{g,n} = \Lambdab_{g,n}$, $\phi^*\etab' =
\etab$, and $\phi(\vac) = \vac$.  Since
$(\chb,\etab,\{\Lambdab_{g,n}\},\vac)$ is a \cft, so is
$(\chb,\etab',\{\Lambdab'_{g,n}\},\vac)$, and $\phi$ is an
isomorphism.
\end{proof}

\begin{rem}
A similar rescaling was observed in \cite{AGV}, and the previous
proposition could be regarded as its origin in the framework of
\Gcfts.
\end{rem}

\section{$G$-stable maps}\label{6}

\ 

In this section we briefly describe the main source of examples of \Gcfts;
namely, Gromov-Witten style classes on the moduli space of $G$-stable maps.

\begin{df}
A \emph{genus $g$, $n$-pointed $G$-stable map} over a base $T$ into a global
quotient $\xg $ is a $G$-equivariant morphism $f:E\rTo X$ from an admissible
$G$-cover $\pi:E\rTo C$ of a genus $g$ prestable curve $C/T$ with $n$
sections $\pt_i:T\rTo E$ such that the induced morphism of stacks $\bar{f}:
[E/G]\rTo \xg $ with marked points $p_i :=\pi\circ \pt_i$ is an $n$-pointed
orbifold (a.k.a.~twisted)  stable map of genus $g$ (as defined in
\cite{CR1,AGV}).

We denote the stack of genus $g$, $n$-pointed $G$-stable maps by
$\MM_{g,n} (X)$, and if $\beta\in H_2(X/G,\nz)$, then we denote
the substack of maps whose image lies in the homology class
$\beta$ by $\MM_{g,n} (X,\beta)$.
\end{df}

\begin{thm}
If the quotient $\xg $ admits a projective coarse moduli space $X/G$,
then the stack $\MM_{g,n} (X , \beta)$ is a proper Deligne-Mumford
stack, which itself admits a projective coarse moduli space.
\end{thm}

The proof follows from the results of \cite{AGV} in essentially the
same way that Theorem~\ref{thm:MMDM} follows from the results of
\cite{ACV}.

There is a natural forgetful morphism $\stq:\MM_{g,n} (X ,
\beta) \rTo \MM_{g,n}$ obtained by forgetting the morphism $f$ and
contracting components in a  manner similar to that described in
the definition of the forgetting tails morphism of
Section~\ref{2}.   There are also natural evaluation morphisms
$ev_i$ from $\MM_{g,n} (X , \beta)$ to the \textsl{inertia variety of
  $X$}, $$\X:=\{(x,g)|
x\in X, \ g\in \stab(x)\} = \coprod_{g\in G} X^g \subseteq X\times
G,$$ with $ev_i((f:E\rTo X, \pt_i)) = (f(\pt_i), m_i)$, where
$m_i$ is the monodromy of $E$ around $\pt_i$ and $X^g$ is the
fixed point locus in $X$ of the subgroup $\langle g \rangle\subseteq G$. These
are compatible in the sense that the following diagram commutes
$$\begin{diagram} \MM_{g,n} (X, \beta) & \rTo^{ev_i} & \X
\\ \dTo^{\stq} &
& \dTo^{\pr_2} \\ \MM_{g,n} & \rTo^{ev_i} &                     G\\
\end{diagram},
$$ where the map $\pr_2$ is the projection onto the second factor and
the lower map $ev_i$ is the $i$th component of the map $\be$ of
Definition~\ref{df:e}.

\begin{df}
We denote by $\MM_{g,n}(X,\beta,\bm)$ the component
$\stq^{-1}(\MM_{g,n}(\bm))$ that maps to $\bm \in G^n$ via $\be
\circ \stq$.
\end{df}

\begin{df}
Let $\ch(X) := H^{2\bullet}(\widehat{X};\Theta) = \bigoplus_{m\in
G} \ch(X)_m$, where $\ch(X)_m := H^{2 \bullet}(X^m;\Theta)$, and
$\Theta$ is the usual ring (see \cite{Ma})
associated to $X$ with generators $\{
q^\beta \}$ over $\nc$,  satisfying $q^{\beta+\beta'} = q^{\beta}
q^{\beta'}$.
\end{df}

\begin{rem}
Of course, one could allow odd-dimensional cohomology classes as
well, after inserting the necessary signs for skew-symmetry, but
for simplicity we will work only with even-dimensional classes.
\end{rem}

In a subsequent paper \cite{IP}, we will describe the details of
how the classes 
\[
\Lambda^{G,X}_{g,n}(v_1,\dots,v_n) :=
\sum_\beta \stq_*(\prod_{i=1}^n ev_i^*(v_i) \cap
[\MM_{g,n}(X)]^{vir}) q^\beta
\]
form a \Gcft, and how the \cft\ of coinvariants of
$\{\Lambda^{G,X}_{g,n}\}$ agrees with the orbifold
Gromov-Witten classes of Chen-Ruan \cite{CR2}.

In the remainder of this section we will briefly treat two special
cases.  In Subsection~\ref{sec:betazero} we describe the case of
$\beta=0$, and show that it gives the ring $H^{\bullet}(X,G)$ of
Fantechi and G\"ottsche---and therefore the stringy orbifold
cohomology of Chen and Ruan---as special cases.  In
Subsection~\ref{sec:trivG} we describe the \Gcft\
$\{\Lambda^{G,X}_{g,n}\}$ for all $\beta$ in the case that $G$ acts
trivially on $X$.

\subsection{The degree zero case, the Fantechi-G\"ottsche ring, and Chen-Ruan
orbifold cohomology} \label{sec:betazero}

\ 

We will now study the case of degree-zero $G$-stable maps in more
detail. We will explicitly prove that the degree-zero $G$-stable
maps endow $\ch(X)$ with the structure of a $G$-Frobenius algebra,
the genus-zero part of which agrees with the ring $H^\bullet(X,G)$
in \cite{FG}. 

\emph{Throughout this section, we will use the ground ring $\nc$
  instead of $\Theta$ in the definition of $\ch(X)$, since we are
  restricting to degree-zero maps.  We will also assume that $X$ is a
  smooth variety with projective coarse moduli space, unless otherwise
  stated.}

\begin{df}
Let $\str : \MM_{g,n}(X,0,\bm) \rTo
\MM_{g,n}(\bm)$ denote the morphism ${\mathrm{st}}_{(X,\beta=0)}$.
We define $$\xi(X,0,\bm) := \str^{-1}(\xi(\bm)).$$

Similarly, if $m,a,b \in G$ are chosen such that $m \in [a,b],$ we
let $$\xi_{1,1}(X,0,(m,a,b)):=\str^{-1}(\xi_{1,1}(m,a,b)).$$

We also define $X^{\langle \bm \rangle}$ to be the locus in $X$ of points fixed by
the subgroup $\langle \bm \rangle \le G $ generated by all of the
elements $m_1, \dots, m_n$ in $\bm$.
\end{df}

Since the marked points $\pt_i$ in the universal $G$-cover $\E$ over
$\xi(\bm)$ all lie in the same connected component of $\E$, it is
straightforward to see that any $G$-stable map $f$ into $X$ of degree
$0$ that maps by $\str$ to $\xi(\bm)$ is determined only by the
underlying $G$-cover (the point $\str([E\to C] \in\xi(\bm)$) and by
the point $f(\pt_1) = \cdots = f(\pt_n)$.  Moreover, the point
$f(\pt_i)$ must have a stabilizer that includes the monodromy element
$m_i$, so the following proposition is now easy to see.

\begin{lm}\label{lm:KeyIsomorphism}
The substack $\xi(X,0,\bm)$ of $\MM_{0,3}(X,0,\bm)$ is canonically
isomorphic to the product $$\xi(X,0,\bm) = \xi(\bm) \times X^{\langle
  \bm \rangle},$$ and the substack $\xi_{1,1}(X,0,(m,a,b))$ is
canonically isomorphic to the product
$$\xi_{1,1}(X,0,(m,a,b))=\xi_{1,1}(m,a,b) \times X^{\langle m,a,b \rangle}.$$
\end{lm}

\begin{proof}
For an object in $\MM_{0,3}(X,0,\bm),$ the isomorphism is given by
the morphism $(E\rTo^f X^{\langle \bm \rangle}; \linebreak[0]
\pt_1,\ldots,\pt_n) \mapsto
(E;\pt_1,\ldots,\pt_n)\times f(\pt_1)$, where $E$ is the $G$-cover
(we have suppressed the underlying curve $C$ since it is
determined by $E$), and its inverse is given by
$(E;\pt_1,\ldots,\pt_n)\times q \mapsto (E\rTo^f X^{\langle \bm\rangle};
\pt_1,\ldots,\pt_n)$, where if $E'$ is the connected component of
$E$ containing the marked points $\pt_1,\ldots,\pt_n$, then
$f(E'):=q$ and $f(\rho(\gamma)\pt') := \rho(\gamma)f(\pt')$ for
all $\gamma$ in $G$ and $\pt'$ in $E'$. The maps for
$\xi_{1,1}(X,0,(m,a,b))$ are similar.
\end{proof}

\subsubsection{The minimal cover $\xi'(\bm)$}

\ 

\begin{df}
Let $G$ be a finite group and fix $\bm$ in $G^n$ such that
$\prod_{i=1}^3 m_i = 1$, and let $G' := \langle \bm \rangle$ denote the
subgroup of $G$ generated by the components of $\bm$. Let
$\xi'(\bm)$ denote the connected component of $\MMp_{0,3}(\bm)$
which is defined in the same way as $\xi(\bm)$ but with the group
$G$ replaced by $G'$.
\end{df}

\begin{lm} \label{lm:SecondIsomorphism}
Let $G$ be a finite group $\bm\in G^3$ with $\prod^3_{i=1} m_i=1$,
and $G' = \langle \bm\rangle $. Consider the morphism
$\Ic:\xi(\bm)\rTo\xi'(\bm)$ taking the object $(E\rTo
C;\pt_1,\ldots,\pt_n)$ to the object $(E'\rTo C;\pt_1,\ldots,\pt_n)$,
where $E'$ is the connected component of $E$ which contains
$\pt_i$ for all $i=1,\ldots,n$.  The morphism $\Ic$ is an
isomorphism.
\end{lm}

\begin{proof}
Since $E'$ is a $G'$-cover (see Appendix of \cite{FG}), $\hat{I}$
is a morphism.

The inverse morphism takes $(E'\rTo C;\pt'_1,\ldots,\pt'_n)$ to 
$(E\rTo C; \pt_1,\ldots,\pt_n)$, where $E = E'\times_{G'}
G$ and $G'$ acts on $E$ from the right in the usual way, $G'$ acts
on $G$ by left multiplication, and $\pt_i := [\pt_i',1]$ for all
$i=1,\ldots,n$.
\end{proof}

Consider the following commutative diagram
\begin{equation}
\begin{diagram}
                   &            & X & \\
                  & \ruTo^{f'} & &\luTo^{f} \\
  \ce'\times X^{\langle \bm \rangle} && \rTo^{\It} &  &\ce\times
  X^{\langle \bm \rangle} \\
  & \rdTo & & & &\rdTo &\\
 \dTo^{\pit'}   &     &\cc' \times X^{\bm}     &    &  \dTo^{\pit}
 &    & \cc \times X^{\bm}    &   \\
 & \ldTo && & &\ldTo\\
 \xi'(\bm)\times X^{\langle \bm \rangle} & &\rTo^I   &  &
 \xi(\bm)\times X^{\langle \bm\rangle}
 &         &  \\
  \dTo^{\pr_{\xi'}}      &&     &      & \dTo^{\pr_{\xi}}
 &         &  \\
 \xi'(\bm)     &         &\rTo^\Ic  & & \xi(\bm)
 &         &  \\
\end{diagram}
\end{equation}
where $\It$ and $I$ are the isomorphisms induced by $\Ic$, $\cc'$
and $\cc$ are the universal curves,  $\ce'$ and $\ce$ are the
universal $G'$ and $G$ covers, respectively, and $f$, $f'$ are
the universal stable maps.

\begin{prop}
$I^* R\pit^G_*(f^* TX)$ is canonically isomorphic to $R\pit_*^{'
G'}(f'^*(TX))$ in the $K$-theory of $\xi'(\bm)\times X^{\langle \bm \rangle}$, where
$R\pit^G_*$ denotes the $G$-invariant derived push-forward, and
$R\pit^{'G'}_*$ denotes the $G'$-invariant derived push-forward.

\end{prop}
\begin{proof}
The fiber of $I^* R\pit^G_*(f^* TX)$ over $\xi'(\bm)\times q$ for
all $q$ in $X^{\langle \bm \rangle}$ is $H^\bullet(\ce\times q,\ct)$, where the
sheaf $\ct$ over $\ce\times q$ is $f^*(TX)$. Since $\ce'$ is the
connected component of $\ce$ containing $\pt_i$ for all
$i=1,\ldots, n$ we have  $\ct{\vert}_{\pt\times q} = T_q X$ for
all $\pt$ in $\ce'$. Henceforth, let us regard $\ct$ as a bundle
over $\ce$ to avoid notational clutter.

Observe that $\ct$ is a $G$-equivariant trivial bundle on $\ce$.
Denote the restriction of $\ct$ to $\ce'$ by $\ct'$ and observe
that it is a $G'$-equivariant bundle. We will now construct a
bundle from $\ct'$ on $\ce'$ which is isomorphic as a
$G$-equivariant bundle to $\ct$ on $\ce$ as follows.

Consider the bundle $\ct'\otimes\co_G$ on $ \ce'\times G$. We
observe that $\ce'$ is a right $G'$-space and $G$ is a left
$G$-space by left multiplication. Similarly, there is a right $G'$
action on $\ct'$ and a left $G'$-action on $\co_G$. Therefore,
$\ct'\otimes\co_G$ over $ \ce'\times G$ is a $G'$-equivariant
bundle with respect to the diagonal $G'$ action. Quotienting by
$G'$ and using the identification of $\co_G$ with $\nc[G]$, we
obtain $\ct'\otimes_{\nc[G']}\nc[G]$ over $\ce'\times_{G'} G,$
which is a $G$-equivariant bundle, where an element $\gammat$ in
$G$ acts upon an element of the base as $[e',\gamma]\mapsto
[e',\gamma\gammat]$, and similarly in the bundle. We now have the
isomorphism of $G$-equivariant vector bundles
\[
\begin{diagram}
\ct'\otimes_{\nc[G']}\nc[G] & \rTo^{\widetilde{\lambda}} & \ct \\
\dTo &     & \dTo \\ \ce'\times_{G'} G & \rTo^\lambda &    \ce , \\
\end{diagram}
\]
where $\widetilde{\lambda}([v',\gamma]) := \rho(\gamma) v'$, and
$\lambda([e',\gamma]) := \rho(\gamma) v'$, and where
$\rho(\gamma)$ indicates the right $G$ action.

Therefore,
\begin{eqnarray*}
H^\bullet(\ce,\ct) &=& H^\bullet(\ce'\times_{G'}
G,\ct'\otimes_{\nc[G']}\co_G) \\ &=& H^\bullet(\ce'\times
G,\ct'\otimes\nc[G])^{G'} \\ &=& (H^\bullet(\ce',\ct')\otimes
H^\bullet(G,\co_G))^{G'} \\ &=& (H^\bullet(\ce',\ct')\otimes
\nc[G])^{G'} \\ &=& H^\bullet(\ce',\ct')\otimes_{\nc[G']} \nc[G].
\end{eqnarray*}
Taking $G$-invariants, we have
\[
H^\bullet(\ce,\ct)^G = (H^\bullet(\ce',\ct')\otimes_{\nc[G']} \nc[G])^G
\cong H^\bullet(\ce',\ct')^{G'}.
\]
The latter is precisely the fiber of $R\pit_*^{' G'}(f'^*(TX))$ over
$\xi'(\bm)\times q$.
\end{proof}

\begin{prop}
When $\beta=0$, the sheaf $R^1\pit^G_*(f^* TX)$ is locally free on
$\xi(\bm)\times X^{\langle\bm\rangle }=\xi(X,0,\bm)$ and the virtual fundamental
class of $\xi(X,0,\bm)$ is simply the top Chern class $c_{top}(R^1
\pit^G_*(f^*TX))$.
\end{prop}

\begin{proof}
This follows immediately from the construction of $\MM_{g,n}(X)$
as a fibered product of sections over $\M_{g,n}([X/G])$, the stack of
orbifold stable maps to $[X/G]$, and the fact that the proposition
holds there (see e.g., \cite{AGV}).
\end{proof}

\begin{df}
Let $c(\bm) := c_\mathrm{top}(R^1\pit^G_*(f^* TX))$ and $c'(\bm)
:= c_\mathrm{top}(R^1\pit^{'G'}_*(f'^* TX))$, where
$c_\mathrm{top}$ denotes the top Chern class.
\end{df}

\begin{crl}
For all $\bm$ in $G^3$ such that $\prod_{i=1}^3 m_i = 1$, we have
\begin{equation}\label{eq:VFCPullback}
I^* c(\bm) = c'(\bm).
\end{equation}
\end{crl}

We now prove that the $3$-point correlator responsible for the multiplication
in the $G$-Frobenius algebra can be identified by the isomorphism
in Lemma \ref{lm:KeyIsomorphism}.
\begin{prop}
For all  $\bm$ in $G^3$  such that $\prod_{i=1}^3 m_i = 1$ and
$\alpha_{m_i}$ in $H^\bullet(X^{m_i})$, let
$\Lambda^{\xi'}_{0,3}(\alpha_\bm)$ in $H^\bullet(\xi'(\bm))$ and
$\Lambda_{0,3}^\xi(\alpha_\bm)$ in $H^\bullet(\xi'(\bm))$ be given
by
\[
\Lambda_{0,3}^{\xi'}(\alpha_\bm) := \pr_{\xi' *}(\prod_{i=1}^3
(ev^{'*}_{m_i} \alpha_{m_i})) \cup c'(\bm))
\]
and
\[
\Lambda_{0,3}^{\xi}(\alpha_\bm) := \pr_{\xi *}((\prod_{i=1}^3
(ev_{m_i}^* \alpha_{m_i})) \cup c(\bm)).
\]
We have
\begin{equation*}
\Ic^*\Lambda_{0,3}^\xi(\alpha_\bm) =
\Lambda_{0,3}^{\xi'}(\alpha_\bm),
\end{equation*}
where $ev_{m_i}:\xi(V,0,\bm)\rTo X^{m_i}$ and
$ev'_{m_i}:\xi'(\bm)\times X^{\langle\bm\rangle }\rTo X^{m_i}$ are the evaluation
morphisms, and $\pr_{\xi'} :\xi'(\bm) \times X^{\langle\bm\rangle }\rTo \xi'(\bm)$
and $\pr_{\xi} : \xi(\bm)\times X^{\langle\bm\rangle }\rTo \xi(\bm)$ are the
projections, which can be identified with the morphism forgetting
the $G$-stable maps.
\end{prop}
\begin{proof}
\begin{eqnarray*}
\Ic^*\Lambda_{0,3}^\xi(\alpha_\bm) &=& \Ic^*\pr_{\xi
*}((\prod_{i=1}^3 (ev_{m_i}^* \alpha_{m_i})) \cup c(\bm)) \\ &=&
\Ic^{-1}_*\pr_{\xi *} ((\prod_{i=1}^3 (ev_{m_i}^* \alpha_{m_i}))
\cup c(\bm))\\ &=& ((\Ic^{-1}\circ\pr_{\xi})_* ((\prod_{i=1}^3
(ev_{m_i}^* \alpha_{m_i})) \cup c(\bm)) \\ &=& (\pr_{\xi'}\circ
I^{-1})_* ((\prod_{i=1}^3 (ev_{m_i}^* \alpha_{m_i})) \cup c(\bm))
\\ &=& \pr_{\xi' *} I^* ((\prod_{i=1}^3 (ev_{m_i}^* \alpha_{m_i}))
\cup c(\bm)) \\ &=& \pr_{\xi' *} ((\prod_{i=1}^3 ((ev_{m_i}\circ
I)^* \alpha_{m_i})) \cup I^*c(\bm)) \\ &=& \pr_{\xi' *}
((\prod_{i=1}^3 ((ev_{m_i}\circ I)^* \alpha_{m_i})) \cup
I^*c(\bm)) \\ &=& \pr_{\xi' *} ((\prod_{i=1}^3 (ev^{'*}_{m_i}
\alpha_{m_i})) \cup c'(\bm)) \\ &=&
\Lambda_{0,3}^{\xi'}(\alpha_\bm),
\end{eqnarray*}
where we have used Equation (\ref{eq:VFCPullback}) in the
penultimate equality.
\end{proof}

\begin{crl}
For all $\bm$ in $G^3$ such that $\prod_{i=1}^3 m_i = 1$, and for
$\alpha_{m_i}$ in $\ch(X)_{m_i}=H^\bullet(X^{m_i})$, we have
\[
\mu(\alpha_\bm) = \int_{\db{\xi(\bm)}}
\Lambda_{0,3}^{\xi}(\alpha_\bm) = \int_{\db{\xi'(\bm)}}
\Lambda_{0,3}^{\xi'}(\alpha_\bm),
\]
where $\mu$ is defined as in Equation (\ref{eq:mu}).
\end{crl}

$\mu$ completely determines the multiplication and metric. We will now prove
that it yields a $G$-Frobenius algebra.

\subsubsection{The genus-zero part of the $G$-Frobenius algebra}

\ 

For this subsection, we can assume, without loss of generality,
that $G' = G$ in light of the results of the previous section.

\begin{df}
Since the virtual fundamental class $c(\bm)$ belongs to
$H^\bullet(\xi(\bm)\times X^{\langle\bm\rangle })\cong H^\bullet(\xi(\bm))\otimes
H^\bullet(X^{\langle\bm\rangle })$, define $\ctt(\bm)$ in $H^\bullet(X^{\langle\bm\rangle })$ to be
the unique class such that
\[
c(\bm) = \vac_{\xi(\bm)}\otimes \ctt(\bm),
\]
where $\vac_{\xi(\bm)}$ is the unit in $H^\bullet(\xi(\bm))$.
\end{df}

We will now write $\mu(v_\bm)$ as an integral over $X^{\langle\bm\rangle }$.
\begin{prop}
For all $v_\bm$ in $\ch(X)_\bm$ where $\prod_{i=1}^3 m_i = 1$, we
have
\begin{equation}\label{eq:Finalmu}
\mu(v_\bm) = \int_{[X^{\langle\bm\rangle }]} (\prod_{i=1}^3 j_{m_i}^*
v_{m_i})\cup\ctt(\bm),
\end{equation}
where $j_{m_i} : X^{\langle\bm\rangle }\rInto X^{m_i}$ is the inclusion and
$[X^{\langle\bm\rangle }]$ is the fundamental class of the variety $X^{\langle\bm\rangle }$. In
particular, when $m_i = 1$ for all $i=1,2,3$, then $\ctt(1,1,1) =
\vac$, the unit in $H^\bullet(X)$. The restriction of the
multiplication and metric to the untwisted sector $\ch(X)_1$ agree
with the usual cup product and metric from $H^{\bullet}(X)$.
Furthermore, $(\chb(X),\mub,\vac)$ is isomorphic as a $\Gb$-graded
Frobenius algebra to the Chen-Ruan orbifold cohomology of $\xg $.
\end{prop}
\begin{proof}
The first statement is a straightforward calculation. The second follows from
the observation that the appropriate obstruction bundle vanishes when $m_1 =
m_2 = m_3 = 1$. The third follows from the following remark and Section 2 of
\cite{FG}.
\end{proof}

\begin{rem}
The vector bundle $R^1\pit_*^G(f^* TX)\rTo \xi(\bm)\times X^{\langle\bm\rangle }$
is \emph{not} the pullback of a vector bundle via the projection
$\xi(\bm)\times X^{\langle\bm\rangle }\rTo X^{\langle\bm\rangle }$ because the the automorphism
group of a $G$-cover (which is isomorphic to $H(\bm)$) in
$\xi(\bm)$ acts non-trivially on  $R^1\pit_*^G(f^* TX)\rTo
\xi(\bm)$ as the action of the automorphism group commutes with
the action of $G$. Nevertheless, one can interpret the bundle
$R^1\pit_*^G(f^* TX)\linebreak[1] \rTo \xi(\bm)\times X^{\langle\bm\rangle }$ as an
$H(\bm)$-equivariant vector bundle $R^1\pit_*^G(f^*
TX)\linebreak[1]\rTo X^{\langle\bm\rangle }$. This bundle can be identified with
the bundle $F(m_1,m_2)\rTo X^{\langle\bm\rangle }$ introduced in \cite{FG}.
Therefore, their cohomology class $c(m_1,m_2)$ can be identified
with $\ctt(\bm)$, which is a class on $X^{\langle\bm\rangle }$, so Equation
(\ref{eq:Finalmu}) is consistent with their multiplication.

They also prove that the vector bundle $F(m_1,m_2)$ restricted
to a  connected component $U$ of $X^{\langle\bm\rangle }$ has rank $a(m_1,U)+a(m_2,U)-a(m_1
m_2,U)-codim(U\subseteq X^{m_1 m_2})$. To explain this notation, 
let $X$ have dimension $D$, $q$ belong to $X$, and $m$ belong to the
isotropy subgroup of $G$ at $q$. Denote the set of eigenvalues of the action
of $m$ on $T_q X$ by $\{ \exp(-2\pi i r_1),\dots,\exp(-2\pi i r_D) \}$ for
%%TK Fixed def of age since we have a right action by inserted - signs in
%%last line.
all $j=1,\ldots,D$ where $r_j$ belongs to the interval $[0,1)$.  The
\emph{age of $m$ in $q$}, $a(m,q)$, is defined to be $\sum_{j=1}^D
r_j$. Since $a(m,q)$ depends only upon the connected component containing
$q$, $a(m,U)$ is defined to be $a(m,q)$ for any $q$ in $U$.
\end{rem}

\begin{prop}
The triple $(\ch(X),\mu,\vac)$ satisfies all of the axioms of a
$G$-Frobenius algebra except, perhaps, for the trace axiom. Our
multiplication, metric, and identity agrees with that \cite{FG} on the
ring $H^\bullet(X,G)$.  Furthermore, $\eta$ on $\ch(X)$ has nonzero
homogeneous components $\ch(X)_{m_+} \otimes \ch(X)_{m_-} \rTo \nc$
only if $m_+ m_- = 1$, in which case $\ch(X)_{m_+} = \ch(X)_{m_-} =
H^\bullet(X^{m_+},\nc)$, and $\eta$ agrees with the usual Poincar\'e
pairing.
\end{prop}
\begin{proof}
This result follows from the previous remark and \cite{FG}.
\end{proof}

\begin{rem}
Proposition \ref{prop:agreement} explains the origin of the factor of
$\frac{1}{|G|}$ in the definition of $\etab$ from the viewpoint of
intersection theory. This factor may be removed, if desired, as per Remark
\ref{rem:StringCC}.
\end{rem}

\subsubsection{The trace axiom}

\ 

We will now prove that the trace axiom, which is a genus-one
condition, holds for $(\ch(X),\mu,\vac)=H^{\bullet}(X,G)$.

\begin{prop}
The Trace Axiom (Definition \ref{df:trace} (\ref{trace})) holds
for the triple $(\ch(X),\mu,\vac).$
\end{prop}

%%Typersetter: the correction starts here.
\begin{proof}
The proof of the trace axiom in Theorem \ref{thm:GFA} shows that
it suffices for us to check that the cutting loops property
\ref{df:GCFT}(\ref{5a}) holds in the special cases of
$\varrho_a:\xi(m_1,bab^{-1},a^{-1}) \rTo \xi_{1,1}(m_1,a,b)$ and
$\varrho_b:\xi(m_1,b,ab^{-1}a^{-1})\rTo \xi_{1,1}(m_1,a,b)$ for
the virtual class.  
We may assume that $G = \langle m,a,b\rangle$, and we denote by $\bm'$ the triple $(m,bab^{-1},b^{-1})$.  Let $H$ 
denote the subgroup $\langle \bm' \rangle$, so that  
$\xi(X,0, \bm') = \xi(\bm')\times X^H$ and $\xi_{1,1}(X,0,(m,a,b)) = \xi_{1,1}(m,a,b) \times X^{G}$. 
It suffices to check that $$j^*c_{top}(R^1\varpi^G_* {f'}^*TX) \cup
e_{\alpha} \cup e_{\beta}\eta^{\alpha\beta} = (\varrho_a\times \id)^*
c_{top}(R^1\pi^G_*f^*TX)$$ (and the same 
for $\varrho_b$), where $e_{\alpha}$ runs over a basis of the Chow ring $A^*(X^a)$ and $e_{\beta}$ runs over a basis for $A^*(X^{a^{-1}})$, and where the morphisms are those of the following
diagram.
\begin{equation}
\begin{diagram}
                                                        &                                               &X      \\
                                                        & \ruTo^{f'}                            &       & \luTo(2,4)^{f}\\
E' \times X^{\langle m,bab^{-1},a^{-1}\rangle }         & \lTo^{\tilde{j}}      & E' \times X^{\langle m,a,b\rangle }& & \\
\dTo^{\varpi}                           &                       & \dTo^{\phi}                        &  &  \\
                                        &                               & E \times X^{\langle m,a,b\rangle }& \rTo^{\vrta}& \ce\times X^{\langle m,a,b\rangle}\\
                                                        &                               & \dTo^{\pit}                           &                       &\dTo^{\pi} \\
\xi(m,bab^{-1},a^{-1})\times X^{\langle m,bab^{-1},a^{-1}\rangle }& \lTo^{j} & \xi(m,bab^{-1},a^{-1})\times X^{\langle m,a, b\rangle } & \rTo^{\varrho_a \times \id} &  \xi_{1,1}(m,a,b)\times X^{\langle m,a, b\rangle } \\
\end{diagram}
\end{equation}
Here $f$ and $f'$ are the universal stable maps from the
universal admissible covers $\ce \times X^G$ and $E' \times
X^{H}$, respectively. The map $j$ is the obvious inclusion
$j:\xi(\bm') \times X^{G }\rTo \xi(\bm') \times X^{H }$, 
and the spaces $E\times X^G$ and $E'\times X^G$ are, respectively, the 
restrictions of the universal admissible covers $\ce\times X^G$ and $E'\times X^H$ to 
$\xi(\bm') \times X^{G }$.   Finally, $\phi$ is the composition of
$\varrho_2(b)$ with the normalization taking the ``unglued'' admissible cover
$E' \times X^G$ of the three-pointed sphere to the
(``glued'') admissible cover $E\times X^G$ of a nodal genus-one curve.

Since $\varrho_a$ is the composition of 
a regular embedding and a flat morphism, and $j$ is a regular embedding, we have
\begin{align*}
j^*c_{top}(R^1\varpi^G_* {f'}^*TX) & = c_{top}(R^1(\pit\circ\phi)^G_*\tilde{j}^*{f'}^*TX) \\
&= c_{top}(R^1\pit^G_* \phi_*\tilde{j}^*{f'}^*TX),
\end{align*}
and
\begin{align*}
(\varrho_a\times \id)^*c_{top}(R^1\pi^G_* {f}^*TX) &= c_{top}(R^1\pit^G_*\vrta^*{f}^*TX)\\
        &= c_{top}(R^1\pit^G_* \phi_*\tilde{j}^*{f'}^*TX).
\end{align*}
We have an obvious short exact sequence on $E\times X^{\langle m,a,b\rangle}$:
\begin{equation}
0 \to \vrta^*f^*TX \to \phi_*\tilde{j}^* {f'}^* TX \to (\phi_*\tilde{j}^* {f'}^* TX)/(\vrta^*f^*TX)
 \to  0  
\end{equation}

Since $\phi$ is the normalization of the nodal curve $E$, obtained by translating a point with monodromy $bab^{-1}$ by $b$ and then gluing to a point with monodromy $a^{-1}$, it follows that the quotient  $(\phi_*\tilde{j}^*{f'}^*TX)/(\vrta^*f^*TX)$ is only supported on the nodal locus, and that
the pushforward $$\pit^G_* \left( (\phi_*\tilde{j}^*{f'}^*TX)/(\vrta^*f^*TX)\right)$$  is equal (in K-theory) to 
$$T(X^{bab^{-1}}\times X^{a^{-1}})/TX^a|_{X^{G}} \cong TX^a|_{X^G}.$$

By the long exact cohomology sequence associated to this short exact sequence, we get the K-theoretic equality
\begin{align}
R^1\pit^G_*\vrta^*{f}^*TX &=  R^1\pit^G_* \phi_*\tilde{j}^*{f'}^*TX \oplus \pit^G_* \left( (\phi_*\tilde{j}^* {f'}^* TX)/(\vrta^*f^*TX)\right) \ominus \pit^G_* \phi_*\tilde{j}^*{f'}^*TX \oplus \pit^G_*\vrta^*{f}^*TX\\
&=R^1\pit^G_* \phi_*\tilde{j}^*{f'}^*TX \oplus  TX^a|_{X^G} \ominus \pit^G_* \phi_*\tilde{j}^*{f'}^*TX \oplus \pit^G_*\vrta^*{f}^*TX
\end{align}

Furthermore, since $H^0(E,\co_E)$ is isomorphic to the trivial $G$-module
$\nc$, and $H^0(E', \co_{E'})$ is isomorphic to the $G$-module
$\nc[H\backslash G]$, we 
have 
\begin{equation}
\pit^G_*\vrta^*{f}^*TX \cong TX^G,\end{equation} and
\begin{equation}\pit^G_* \phi_*\tilde{j}^*{f'}^*TX \cong TX^H|_{X^G}.\end{equation}
That is to say, \begin{equation}R^1\pit^G_*\vrta^*{f}^*TX =  R^1\pit^G_* \phi_*\tilde{j}^*{f'}^*TX \oplus \mathfrak{E},\end{equation}
where $\mathfrak{E}$ is the excess intersection bundle of the diagram
\begin{equation}
\begin{diagram}
\xi(\bm')\times X^G                     & \rTo^j        & \xi(\bm')\times X^H\\
\dTo^{q}                                                &                       & \dTo^{\delta}\\
X^a                                                             & \rTo^{\Delta}& X^{bab^{-1}} \times X^{a^{-1}},
\end{diagram}
\end{equation}
where the map $q$ is the composition of the obvious inclusion followed by the
second projection $\xi(\bm')\times X^G \rTo \xi(\bm')\times X^a \rTo  X^a$,
the map $\Delta$ is the composition of the diagonal followed by the action
$\varrho(b)$ in the first factor and inversion in the second: 
$X^a \rTo X^a \times X^a \rTo X^{bab^{-1}} \times X^{a^{-1}}$, and the map $\delta$ is the product of the evaluation maps: $\delta = ev_2 \times ev_3$.  

The excess intersection formula now gives that
\begin{equation}c_{top}(R^1\pit^G_*\vrta^*{f}^*TX) =  c_{top}(R^1\pit^G_*
\phi_*\tilde{j}^*{f'}^*TX) \cup j^*\delta^*\Delta_* \vac,\end{equation} and it is
straightforward to see that this last term is the desired sum $e_{\alpha}
\cup e_{\beta} \eta^{\alpha \beta}$.  
\end{proof}
%%%%%%%%%%end corrections
\begin{rem}
Finally, we note that the $G$-Frobenius algebra $(\ch(X), \mu, \vac)$ enjoys
some functoriality properties, as Fantechi-G\"ottsche have showed that it
pulls back along \'etale maps \cite[pg. 11]{FG}.
\end{rem}

\subsubsection{Tensor products}

\ 

We now work out the tensor products of the equivariant \cfts\
described above and show that they reduce to the obvious notions
of tensor products for $G$-Frobenius algebras.

\begin{prop}
Let $X'$ be a smooth, projective variety with a $G'$-action and let
$ ((\ch(X'),\rho'), \mu', \vac')$ be the $G'$-Frobenius algebra associated
to contributions from maps of degree zero where $\mu'$ is defined by Equation
(\ref{eq:mu}).  Let $X''$ be a smooth, projective variety with a $G''$-action
and let $ ((\ch(X''),\rho''), \mu'', \vac'')$ be its similarly associated
$G''$-Frobenius algebra.

\begin{enumerate}
\item Consider $X'\times X''$ with its $G'\times G''$ action.
The associated $G'\times G''$-Frobenius algebra $((\ch(X'\times
X''),\rho), \mu, \vac)$ is canonically isomorphic to the external
tensor product of $((\ch(X'),\rho'), \mu', \vac')$ and
$((\ch(X''),\rho''), \mu'', \vac'')$.
\item Suppose that $G' = G'' = G$, and consider $X'\times X''$ with its
diagonal $G$ action. Its associated $G$-Frobenius algebra
$((\ch(X'\times X''),\rho),\mu,\vac)$ is canonically isomorphic to
the tensor product of $((\ch(X'),\rho'), \mu', \vac')$ and
$((\ch(X''),\rho''), \mu'', \vac'')$.
\end{enumerate}
\end{prop}
\begin{proof}
To prove the first part, we need to understand the behavior of the
obstruction bundle when $g=0$, $n=3$, and $\beta=0$.

For all $\bm'$ in $G'^n$ and $\bm''$ in $G''^n$, consider the
substack $\xi^{G'\times G''}(\bm'\times\bm'')$ of $\M^{G'\times
G''}_{0,3}(\bm'\times\bm'')$, where we adopt the notation
$\bm'\times\bm''$ from Proposition \ref{prop:ExternalFP}. Since
$\M_{0,3}$ is a point, $\xi^{G'\times G''}(\bm'\times\bm'')$ is
canonically isomorphic to $\xi^{G'}(\bm')\times \xi^{G''}(\bm'')$.
Similarly, $\xi^{G'\times G''}(X'\times X'',0,\bm'\times\bm'')$ is
canonically isomorphic to $\xi^{G'}(\bm')\times X'^{\langle\bm'\rangle }\times
\xi^{G''}(\bm'')\times X''^{\langle\bm''\rangle }$. Thus, we have
\[
X'\times X'' \lTo^{f'\times f''} \ce'\times X'^{\langle \bm' \rangle }\times
\ce''\times X''^{\langle \bm'' \rangle }\rTo^{\pit'\times\pit''}
\xi^{G'}(\bm')\times X'^{\langle \bm' \rangle }\times \xi^{G''}(\bm'')\times
X''^{\langle \bm'' \rangle },
\]
where $\ce'$ and $\ce''$ are the universal curves, and $f'$ and $f''$ are
universal evaluation morphisms. It is easy to see that there is a canonical
isomorphism
\[
R^\bullet(\pit'\times\pit'')_*^{G'\times G''}(f'\times f'')^*(T(X'\times
X'')) \cong R^\bullet\pit'^{'G'}_*( f'^*(TX'))\oplus R^\bullet\pit''^{G''}_*
(f''^*(TX'')).
\]
By multiplicativity of the top Chern class, we obtain
\begin{equation}\label{eq:VFCTensor}
c(\bm'\times\bm'') = c'(\bm')\otimes c''(\bm''),
\end{equation}
where $c(\bm'\times\bm'')$ is the virtual fundamental class of
$\xi^{G'\times G''}(X'\times X'',0,\bm'\times\bm'')\cong
\xi^{G'}(X',0,\bm')\times \xi^{G''}(X'',0,\bm'')$, $c'(\bm')$ is the virtual
fundamental class of $\xi^{G'}(X',0,\bm')$, and $c''(\bm'')$ is the virtual
fundamental class of $\xi^{G''}(X'',0,\bm'')$.

It is a straightforward exercise using Equation
(\ref{eq:VFCTensor}) and the identification $\ch(X'\times
X'')\cong \ch(X')\otimes \ch(X'')$ to show that $$\mu =
\mu'\otimes \mu''.$$ The trace axiom then follows immediately from
the fact that trace of a tensor product is the product of the
corresponding traces over each tensor factor.  This finishes the
proof of the first part of the proposition.

The second part of the proposition follows from Remark
\ref{rem:FactoringTensors} and the first part of this proposition.

\end{proof}

\subsection{Trivial $G$-actions}\label{sec:trivG}

\ 

In the special case that the action of $G$ on $X$ is trivial, the
data of a $G$-stable map to $X$ is the same as a stable map from
the underlying curve $C$ to $X$ and the data of a pointed 
admissible $G$-cover, that is $$\MM_{g,n}(X)=\MM_{g,n}
\times_{\M_{g,n}}\M_{g,n}(X).$$

Moreover, since $\MM_{g,n}$ is smooth, it is evident that the
virtual fundamental class on $\MM_{g,n}(X)$ is simply the pullback
of the virtual fundamental class of $\M_{g,n}(X)$
$$[\MM_{g,n}(X)]^{vir}=\pr^*_2 [\M_{g,n}(X)]^{vir},$$ and the
evaluation map $\MM_{g,n}(X)\rTo (\widehat{X})^n=(X \times G)^n$ is
simply the product of the evaluation maps $\be: \MM_{g,n}\rTo G^n$
and $\mathbf{ev}:\M_{g,n}(X)\rTo X^n.$

Thus in this special case, we have $$\Lambda^{G,X}_{g,n}=\st^*
\Lambda^X_{g,n},$$ where $\{\Lambda^X_{g,n}\}$ is the usual
Gromov-Witten \cft\ for $X$, and $\st:\MM_{g,n}\rTo \M_{g,n} $ is
the forgetful map ($\st:= \sth \circ \stt$).

Since $G$ acts trivially, we have
$$\ch(X)_{m_i}=H^{\bullet}(X^{m_i};\Theta) \cong H^{\bullet} (X;
\Theta)$$ for every $m_i \in G$.  So the state space $\ch(X)$ is
just $$H^{\bullet}(X;\Theta)\otimes \nc[G]$$ and
$$\Lambda^{G,X}_{g,n}((v_1 \otimes m_1),\dots,(v_n \otimes
m_n))=\st^*\Lambda^X_{g,n}(v_1,\dots, v_n) \cup \be^*(\vac),$$
which is clearly just the external tensor product of
$\Lambda^X_{g,n}$ with $\nc[G]$.  Thus we have proved the
following.

\begin{prop}
If $X$ is a smooth, projective variety with a trivial $G$ action,
then its associated \Gcft\ is isomorphic to the external tensor
product of the \cft\ of stable maps associated to $X$ (regarded as
an equivariant \cft\ for the trivial group) with $\nc[G]$, the
group ring \Gcft\  (see Example \ref{ex:grpring}).
\end{prop}

\begin{rem}
In the previous example, the induced \cft\ on the space of
$G$-coinvariants agrees with Proposition 3.7 in \cite{JK}.
\end{rem}

\bibliographystyle{amsplain}

\begin{thebibliography}{WDV}

\bibitem[ACV03]{ACV} D.~Abramovich, A~Corti, and A.~Vistoli,
\emph{Twisted bundles and admissible covers}.  Commun. Alg. \textbf{31} no 8 (2003), 3547--3618.
\url{math.AG/0106211}.

\bibitem[AGV02]{AGV} D.~Abramovich, T.~Graber, and A.~Vistoli,
\emph{Algebraic orbifold quantum products}. In A. Adem, J. Morava,
and Y. Ruan (eds.), {Orbifolds in Mathematics and Physics},
\emph{Contemp. Math.}, Amer. Math. Soc., Providence, RI.
\textbf{310}, (2002), 1--25. \url{math.AG/0112004}.

\bibitem[BK01]{BaKi} B.~Bakalov and A.~Kirillov, \emph{Lectures on tensor
categories and modular functors}.  Univ. Lecture Series
\textbf{21}, Amer. Math. Soc., Providence, RI, 2001.

\bibitem[Cos03]{Cos} K.~Costello, \emph{Higher-genus Gromov-Witten
invariants as genus 0 invariants of symmetric   products.}
\url{math.AG/0303387} 

\bibitem[CoKa99]{CoKa} D.~Cox and S.~Katz, \emph{Mirror symmetry and algebraic
geometry}. Math. Surv. \textbf{68}, Amer. Math. Soc., Providence,
RI, 1999.

\bibitem[CR00]{CR1} W.~Chen and Y.~Ruan, \emph{A new cohomology theory for
orbifold.} Commun. Math. Phys. {\bf 248} (2004), no. 1, 1--31.  \url{math.AG/0004129}.

\bibitem[CR02]{CR2} \bysame, \emph{Orbifold Gromov-Witten theory.} In A. Adem, J. Morava, and Y. Ruan
(eds.), {Orbifolds in Mathematics and Physics}, \emph{Contemp.
Math.}, Amer. Math. Soc., Providence, RI. \textbf{310}, (2002),
25--85. \url{math.AG/0103156}.

\bibitem[FG03]{FG} B.~Fantechi and L.~ G\"{o}ttsche, \emph{Orbifold
cohomology for global quotients}. Duke Math. J. \textbf{117} (2003),
no. 2, 197--227. \url{math.AG/0104207}.

\bibitem[FS03]{FrSz} E.~Frenkel and M.~Szczesny, \emph{Chiral de Rham
Complex and Orbifolds}. \url{math.AG/0307181}.

\bibitem[Fu95]{Fu} W.~Fulton, \emph{Algebraic topology}.  Springer, NY, 1995.

\bibitem[IP]{IP} In preparation.

\bibitem[JK02]{JK} T.~Jarvis and T.~Kimura, \emph{The orbifold quantum
cohomology of the classifying space of a finite group}.  In A.
Adem, J. Morava, and Y. Ruan (eds.), {Orbifolds in Mathematics and
Physics}, \emph{Contemp. Math.}, Amer. Math. Soc., Providence, RI.
\textbf{310}, (2002), 123--134. \url{math.AG/0112037}.

\bibitem[JKV01]{JKV} T.~Jarvis, T.~Kimura, and A.~Vaintrob, \emph{Moduli spaces of higher
spin curves and integrable hierarchies}. Compositio Math.
\textbf{126} (2001), no. 2, 157--212. \url{math.AG/9905034}.

\bibitem[Kas95]{Kas} C.~Kassel, \emph{Quantum groups}. Springer-Verlag, New
York, NY, 1995.

\bibitem[Kau02]{Ka} R.~Kaufmann, \emph{Orbifold Frobenius algebras, cobordisms,
and monodromies}. In A. Adem, J. Morava, and Y. Ruan (eds.),
{Orbifolds in Mathematics and Physics}, \emph{Contemp. Math.},
Amer. Math. Soc., Providence, RI. \textbf{310}, (2002),
135--162.

\bibitem[Kau03]{Ka2} R.~Kaufmann, \emph{Orbifolding Frobenius algebras}.  Int. J. of
Math. \textbf{14}, (2003), no. 6, 573--617. \url{math.AG/0107163}.

\bibitem[Ki02]{Ki} A.~Kirillov, Jr. \emph{Modular categories and orbifold
models}. Comm. Math. Phys. \textbf{229} (2002), no. 2, 309--335.

\bibitem[Kir03]{Ki2} \bysame, \emph{On the modular functor associated
with a finite group}. \url{math.QA/0310087} 

\bibitem[KM94]{KoMa}
M.~Kontsevich and Yu.~Manin, \emph{Gromov-Witten classes, quantum
cohomology, and enumerative geometry}. Commun. Math. Phys.
\textbf{164} (1994), 525--562.

\bibitem[Ma99]{Ma}
Yu.~Manin, \emph{Frobenius manifolds, quantum cohomology, and
moduli spaces}.  Colloquium Publ., \textbf{47},  Amer. Math. Soc.,
Providence, RI, 1999.

\bibitem[Mo01]{Mo}
G.~Moore, Lectures on \emph{D-branes, RR-fields, and K-theory}.
Miniprogram (June 18 - July 13, 2001) The Duality Workshop,
Institute for Theoretical Physics, Santa Barbara, CA.
\url{http://online.itp.ucsb.edu/online/mp01/moore2/}

\bibitem[PV01]{PoVa}
A.~Polishchuk and A.~Vaintrob, \emph{Algebraic construction of
Witten's top Chern class}.  In E. Previato (ed.), Contemp. Math.
\textbf{276}, Amer. Math. Soc., (2001), 229--249.
\url{math.AG/0011032}.

\bibitem[P02]{Po2}
A.~Polishchuk, \emph{Witten's top Chern class on the moduli space
of higher spin curves}. \url{math.AG/0208112}.
 
\bibitem[Szc04]{szc04} Matthew Szczesny,  
\emph{Orbifold conformal blocks and the stack of pointed G-covers}. Preprint. \url{math.AG/0408123}


\bibitem[T99]{Tu}
V.~ Turaev, \emph{Homotopy field theory in dimension 2 and
group-algebras}. \url{math.QA/9910010}.

\end{thebibliography}

\providecommand{\bysame}{\leavevmode\hbox
to3em{\hrulefill}\thinspace}

\end{document}